\documentclass{article}
\usepackage{amsmath,amssymb}
\usepackage{algorithm,algpseudocode}  
\usepackage{listings,url,verbatim}    
\usepackage{multicol}
\usepackage{multirow}         
\usepackage[pdftex]{graphicx} 
\usepackage{subfigure}        
\usepackage[usenames]{color}  

\graphicspath{{FIGURES/}}

\title{Communication-optimal parallel and sequential QR and LU factorizations: theory and practice}
\author{James Demmel, Laura Grigori, \\ Mark Hoemmen, and Julien Langou}
\date{\today}

\DeclareMathOperator{\Span}{span}

\DeclareMathOperator{\House}{House}

\newtheorem{theorem}{Theorem}
\newtheorem{lemma}{Lemma}
\newtheorem{corollary}{Corollary}
\newenvironment{proof}{\sl \noindent Proof: \rm}{$\Box$}

\begin{document}
\lstset{basicstyle=\small,showstringspaces=false,language=C} 

\maketitle

\begin{abstract}
  We present parallel and sequential dense QR factorization algorithms
  that are both \emph{optimal} (up to polylogarithmic factors) in the
  amount of communication they perform, and just as \emph{stable} as
  Householder QR.  Our first algorithm, Tall Skinny QR (TSQR), factors
  $m \times n$ matrices in a one-dimensional (1-D) block cyclic row
  layout, and is optimized for $m \gg n$.  Our second algorithm, CAQR
  (Communication-Avoiding QR), factors general rectangular matrices
  distributed in a two-dimensional block cyclic layout.  It invokes
  TSQR for each block column factorization.

  The new algorithms are superior in both theory and practice.  We
  have extended known lower bounds on communication for sequential and
  parallel matrix multiplication to provide latency lower bounds, and
  show these bounds apply to the LU and QR decompositions.  We not
  only show that our QR algorithms attain these lower bounds (up to
  polylogarithmic factors), but that existing LAPACK and ScaLAPACK
  algorithms perform asymptotically more communication.  We also point
  out recent LU algorithms in the literature that attain at least some
  of these lower bounds.

  Both TSQR and CAQR have asymptotically lower latency cost in the
  parallel case, and asymptotically lower latency and bandwidth costs
  in the sequential case.  In practice, we have implemented parallel
  TSQR on several machines, with speedups of up to $6.7\times$ on 16
  processors of a Pentium III cluster, and up to $4\times$ on 32
  processors of a BlueGene/L.  We have also implemented sequential
  TSQR on a laptop for matrices that do not fit in DRAM, so that slow
  memory is disk.  Our out-of-DRAM implementation was as little as
  $2\times$ slower than the predicted runtime as though DRAM were
  infinite.

  We have also modeled the performance of our parallel CAQR algorithm,
  yielding predicted speedups over ScaLAPACK's \lstinline!PDGEQRF! of
  up to $9.7\times$ on an IBM Power5, up to $22.9\times$ on a model
  Petascale machine, and up to $5.3\times$ on a model of the Grid.
\end{abstract}

\tableofcontents

\section{Introduction}\label{S:introduction}

The large and increasing costs of communication motivate redesigning
algorithms to avoid communication whenever possible.  Communication
matters for both parallel and sequential algorithms.  In the parallel
case, it refers to messages between processors, which may be sent over
a network or via a shared memory.  In the sequential case, it refers
to data movement between different levels of the memory hierarchy.
Many authors have pointed out the exponentially growing gaps between
floating-point arithmetic rate, bandwidth, and latency, for many
different storage devices and networks on modern high-performance
computers (see e.g., Graham et al.\ \cite{graham2005getting}).

We present parallel and sequential dense QR factorization algorithms
that are both \emph{optimal} (sometimes only up to polylogarithmic factors) 
in the amount of communication they perform, and just as \emph{stable} as
Householder QR.  Some of the algorithms are novel, and some extend
earlier work.  The first set of algorithms, ``Tall Skinny QR'' (TSQR),
are for matrices with many more rows than columns, and the second set,
``Communication-Avoiding QR'' (CAQR), are for general rectangular
matrices.  The algorithms have significantly lower latency cost in the
parallel case, and significantly lower latency and bandwidth costs in
the sequential case, than existing algorithms in LAPACK and ScaLAPACK.
Our algorithms are numerically stable in the same senses as in LAPACK
and ScaLAPACK.

The new algorithms are superior in both theory and practice.  We have
extended known lower bounds on communication for sequential and
parallel matrix multiplication (see Hong and Kung \cite{hong1981io}
and Irony, Toledo, and Tiskin \cite{irony2004communication}) to QR
decomposition, and shown both that the new algorithms attain these
lower bounds (sometimes only up to polylogarithmic factors), 
whereas existing LAPACK
and ScaLAPACK algorithms perform asymptotically more communication.
(LAPACK costs more in both latency and bandwidth, and ScaLAPACK in
latency; it turns out that ScaLAPACK already uses optimal bandwidth.)
Operation counts are shown in Tables
\ref{tbl:1-par-tsqr}--\ref{tbl:6-seq-caqr-square}, and will be
discussed below in more detail.

In practice, we have implemented parallel TSQR on several machines,
with significant speedups:
\begin{itemize}
\item up to $6.7\times$ on 16 processors of a Pentium III cluster, for
  a $100,000 \times 200$ matrix; and
\item up to $4\times$ on 32 processors of a BlueGene/L, for a
  $1,000,000 \times 50$ matrix.
\end{itemize}
Some of this speedup is enabled by TSQR being able to use a much
better local QR decomposition than ScaLAPACK can use, such as the
recursive variant by Elmroth and Gustavson (see
\cite{elmroth2000applying} and the performance results in Section
\ref{S:TSQR:perfres}).  We have also implemented sequential TSQR on a
laptop for matrices that do not fit in DRAM, so that slow memory is
disk.  This requires a special implementation in order to run at all,
since virtual memory does not accommodate matrices of the sizes we
tried.  By extrapolating runtime from matrices that do fit in DRAM, we
can say that our out-of-DRAM implementation was as little as $2\times$
slower than the predicted runtime as though DRAM were infinite.

We have also modeled the performance of our parallel CAQR algorithm
(whose actual implementation and measurement is future work), yielding
predicted speedups over ScaLAPACK's \lstinline!PDGEQRF! of up to
$9.7\times$ on an IBM Power5, up to $22.9\times$ on a model Petascale
machine, and up to $5.3\times$ on a model of the Grid.  The best
speedups occur for the largest number of processors used, and for
matrices that do not fill all of memory, since in this case latency
costs dominate. In general, when the largest possible matrices are
used, computation costs dominate the communication costs and improved
communication does not help.

Tables \ref{tbl:1-par-tsqr}--\ref{tbl:6-seq-caqr-square} summarize our
performance models for TSQR, CAQR, and ScaLAPACK's sequential and
parallel QR factorizations.  
We omit lower order terms.
In these tables, we make the optimal
choice of matrix layout for each algorithm.  In the parallel case,
that means choosing the block size $b$ as well as the processor grid
dimensions $P_r \times P_c$ in the 2-D block cyclic layout.  (See
Sections \ref{SS:CAQR:par:opt} and \ref{SS:PDGEQRF:par:opt} for
discussion of how to choose these parameters for parallel CAQR resp.\
ScaLAPACK.)  In case the matrix layout is fixed, 
Table~\ref{tbl:CAQR:par:model} in Section~\ref{S:CAQR} gives a general
performance model of parallel CAQR and \lstinline!PDGEQRF! as a
function of the block size $b$ (we assume square $b \times b$ blocks)
and the processor grid dimensions $P_r$ and $P_c$.  
(Table~\ref{tbl:CAQR:par:model} shows that for fixed $b$, $P_r$
and $P_c$, the number of flops and words transferred roughly match,
but the number of messages is about $b$ times lower for CAQR.)
In the sequential
case, choosing the optimal matrix layout means choosing the dimensions
of the matrix block in fast memory so as to minimize runtime with
respect to the fast memory size $W$.  (See Sections
\ref{SS:CAQR-seq-detailed:opt} and \ref{SS:seqLL:factor} for
discussion of how to choose these parameters for sequential CAQR
resp.\ (Sca)LAPACK QR.)
Equation \eqref{eq:CAQR:seq:modeltime:P} in Appendix
\ref{SS:CAQR-seq-detailed:factor} gives the performance model of
sequential CAQR as a function of the dimensions of the matrix block in
fast memory (or rather, as a function of the $P_r \times P_c$ block
layout, which uniquely determines the matrix block dimensions).

\begin{table}[h]
\small
\centering
\begin{tabular}{l | l | l | l}
         & TSQR & \lstinline!PDGEQRF! & Lower bound \\ \hline
\# flops & $\frac{2mn^2}{P} + \frac{2}{3}n^3 \log P$ 
         & $\frac{2mn^2}{P} - \frac{2n^3}{3P}$
         & $\Theta\left( \frac{mn^2}{P} \right)$ \\ 
\# words & $\frac{n^2}{2} \log P$ 
         & $\frac{n^2}{2} \log P$ 
         & $\frac{n^2}{2} \log P$ \\
\textbf{\# messages}
         & $\mathbf{\log P}$ 
         & $\mathbf{2n \log P}$ 
         & $\mathbf{\log P}$ \\ 
\end{tabular}
\caption{Performance models of parallel TSQR and ScaLAPACK's parallel
  QR factorization \lstinline!PDGEQRF! on an $m \times n$ matrix with 
  $P$ processors, along with lower bounds on the number of flops,
  words, and messages.  We assume $m/P \geq n$.  Everything (messages, 
  words, and flops) is counted along the critical path.  The boldface
  part of the table highlights TSQR's improvement over ScaLAPACK.}
\label{tbl:1-par-tsqr}
\end{table}

\begin{table}[h]
\tiny
\centering
\begin{tabular}{l | l | l | l}
         & Par.\ CAQR & \lstinline!PDGEQRF! & Lower bound \\ \hline
\# flops & $\frac{2mn^2}{P} + \frac{2n^3}{3}$
         & $\frac{2mn^2}{P} + \frac{2n^3}{3}$
         & $\Theta\left( \frac{mn^2}{P} \right)$ \\ 
\# words & $\sqrt{\frac{m n^3}{P}} \log P
            - \frac{1}{4} \sqrt{\frac{n^5}{m P}} 
              \log\left( \frac{n P}{m} \right)$    
         & $\sqrt{\frac{m n^3}{P}} \log P
            - \frac{1}{4} \sqrt{\frac{n^5}{m P}} 
              \log\left( \frac{n P}{m} \right)$    
         & $\Theta\left( \sqrt{\frac{m n^3}{P}} \right)$ \\
\textbf{\# messages} 
         & $\mathbf{ \frac{1}{4}
            \sqrt{\frac{n P}{m}}
            \log^2\left( 
                \frac{m P}{n} 
            \right) 
            \cdot \log\left( 
                P \sqrt{\frac{m P}{n}}
            \right) }$ 
         & $\mathbf{ \frac{n}{4} 
              \log\left( \frac{m P^5}{n} \right)
              \log\left( \frac{m P}{n} \right) }$
         & $\mathbf{ \Theta\left( \sqrt{\frac{nP}{m}} \right) }$ \\
\end{tabular}
\caption{Performance models of parallel CAQR and ScaLAPACK's parallel
  QR factorization \lstinline!PDGEQRF! on a $m \times n$ matrix with 
  $P$ processors, along with lower bounds on the number of flops,
  words, and messages.  The matrix is stored in a 2-D $P_r \times P_c$ 
  block cyclic layout with square $b \times b$ blocks.  We choose $b$,
  $P_r$, and $P_c$ optimally and independently for each algorithm.
  We assume $m \geq n$.  Everything (messages, words, and flops) is 
  counted along the critical path.  The boldface part of the table 
  highlights CAQR's improvement over ScaLAPACK.}
\label{tbl:2-par-caqr-general}
\end{table}

\begin{table}[h]
\small
\centering
\begin{tabular}{l | l | l | l}
         & Par.\ CAQR & \lstinline!PDGEQRF! & Lower bound \\ \hline
\# flops & $\frac{4n^3}{3P}$
         & $\frac{4n^3}{3P}$
         & $\Theta\left( \frac{n^3}{P} \right)$ \\ 
\# words & $\frac{3n^2}{4 \sqrt{P}} \log P$
         & $\frac{3n^2}{4 \sqrt{P}} \log P$  
         & $\Theta\left( \frac{n^2}{\sqrt{P}} \right)$ \\
\textbf{\# messages} 
         & $\mathbf{ \frac{3}{8} \sqrt{P} \log^3 P }$
         & $\mathbf{ \frac{5n}{4} \log^2 P }$
         & $\mathbf{ \Theta\left( \sqrt{P} \right) }$ \\
\end{tabular}
\caption{Performance models of parallel CAQR and ScaLAPACK's parallel
  QR factorization \lstinline!PDGEQRF! on a square $n \times n$ matrix with 
  $P$ processors, along with lower bounds on the number of flops,
  words, and messages.  The matrix is stored in a 2-D $P_r \times P_c$ 
  block cyclic layout with square $b \times b$ blocks.  We choose $b$,
  $P_r$, and $P_c$ optimally and independently for each algorithm.
  Everything (messages, words, and flops) is 
  counted along the critical path.  The boldface part of the table 
  highlights CAQR's improvement over ScaLAPACK.}
\label{tbl:3-par-caqr-square}
\end{table}

\begin{table}[h]
\small
\centering
\begin{tabular}{l | l | l | l}
         & Seq.\ TSQR & Householder QR & Lower bound \\ \hline
\# flops & $2mn^2$
         & $2mn^2$
         & $\Theta(mn^2)$ \\
\# words & $2mn$
         & $\frac{m^2 n^2}{2W}$
         & $2mn$ \\
\# messages & $\frac{2mn}{\widetilde{W}}$
            & $\frac{mn^2}{2W}$
            & $\frac{2mn}{W}$ \\
\end{tabular}
\caption{Performance models of sequential TSQR and blocked sequential
  Householder QR (either LAPACK's in-DRAM \lstinline!DGEQRF! or ScaLAPACK's
  out-of-DRAM \lstinline!PFDGEQRF!) on an $m \times n$ matrix with
  fast memory size $W$, along with lower bounds on the number of flops,
  words, and messages.  We assume $m \gg n$ and $W \geq 3n^2/2$.  The 
  boldface part of the table highlights TSQR's improvement over (Sca)LAPACK.
  $\widetilde{W} = W - n(n+1)/2$, which is at least about $\frac{2}{3}W$.}
\label{tbl:4-seq-tsqr}
\end{table}

\begin{table}[h]
\small
\centering
\begin{tabular}{l | l | l | l}
         & Seq.\ CAQR & Householder QR & Lower bound \\ \hline
\# flops & $2mn^2 - \frac{2n^3}{3}$
         & $2mn^2 - \frac{2n^3}{3}$
         & $\Theta(mn^2)$ \\
\textbf{\# words}
         & $\mathbf{ 3 \frac{mn^2}{\sqrt{W}} }$
         & $\mathbf{ \frac{m^2 n^2}{2W} - \frac{m n^3}{6W}
            + \frac{3mn}{2} - \frac{3n^2}{4} }$ 
         & $\mathbf{ \Theta(\frac{mn^2}{\sqrt{W}} )}$ \\
\textbf{\# messages}
         & $\mathbf{ 12 \frac{mn^2}{W^{3/2}} }$ 
         & $\mathbf{ \frac{mn^2}{2W} + \frac{2mn}{W} }$
         & $\mathbf{ \Theta(\frac{mn^2}{W^{3/2}} )}$ \\
\end{tabular}
\caption{Performance models of sequential CAQR and blocked sequential
  Householder QR (either LAPACK's in-DRAM \lstinline!DGEQRF! or ScaLAPACK's
  out-of-DRAM \lstinline!PFDGEQRF!) on an $m \times n$ matrix with
  fast memory size $W$, along with lower bounds on the number of flops,
  words, and messages.  
  The boldface
  part of the table highlights CAQR's improvement over (Sca)LAPACK.}
\label{tbl:5-seq-caqr-general}
\end{table}

\begin{table}[h]
\small
\centering
\begin{tabular}{l | l | l | l}
         & Seq.\ CAQR & Householder QR & Lower bound \\ \hline
\# flops & $\frac{4n^3}{3}$
         & $\frac{4n^3}{3}$
         & $\Theta(n^3)$ \\
\textbf{\# words}
         & $\mathbf{ 3 \frac{n^3}{\sqrt{W}} }$ 
         & $\mathbf{ \frac{n^4}{3W} + \frac{3n^2}{4} }$ 
         & $\mathbf{ \Theta(\frac{n^3}{\sqrt{W}} )}$ \\       
\textbf{\# messages}
         & $\mathbf{ 12 \frac{n^3}{W^{3/2}} }$ 
         & $\mathbf{ \frac{n^3}{2W} }$
         & $\mathbf{ \Theta(\frac{n^3}{W^{3/2}} )}$ \\       
\end{tabular}
\caption{Performance models of sequential CAQR and blocked sequential
  Householder QR (either LAPACK's in-DRAM \lstinline!DGEQRF! or ScaLAPACK's
  out-of-DRAM \lstinline!PFDGEQRF!) on a square $n \times n$ matrix with
  fast memory size $W$, along with lower bounds on the number of flops,
  words, and messages.  The boldface part of the table highlights
  CAQR's improvement over (Sca)LAPACK.}
\label{tbl:6-seq-caqr-square}
\end{table}

Here are highlights of the six tables in this section.  Tables
\ref{tbl:1-par-tsqr}--\ref{tbl:3-par-caqr-square} concern parallel
algorithms.  First, Table \ref{tbl:1-par-tsqr} compares parallel TSQR
and ScaLAPACK's parallel QR factorization \lstinline!PDGEQRF!.  TSQR
requires fewer messages: $\log P$, which is both optimal, and a factor
$2n$ fewer messages than ScaLAPACK.  Table
\ref{tbl:2-par-caqr-general} compares parallel CAQR and
\lstinline!PDGEQRF! on a general rectangular matrix.  Parallel CAQR
needs fewer messages: $\Theta(\sqrt{nP/m})$, which is both optimal
(modulo polylogarithmic factors), and a factor $\Theta(\sqrt{mn/P})$
fewer messages than ScaLAPACK.  Note that $\sqrt{mn/P}$ is the square
root of each processor's local memory size, up to a small constant
factor.  Table \ref{tbl:3-par-caqr-square} presents the same
comparison for the special case of a square $n \times n$ matrix.
There again, parallel CAQR requires fewer messages:
$\Theta(\sqrt{P})$, which is both optimal and a factor
$\Theta(n/\sqrt{P})$ fewer messages than \lstinline!PDGEQRF!.  This
factor is the square root of the local memory size, up to a small
constant factor.

Next, Tables \ref{tbl:4-seq-tsqr}--\ref{tbl:6-seq-caqr-square} concern
sequential QR factorization algorithms.  Table \ref{tbl:4-seq-tsqr}
compares sequential TSQR with sequential blocked Householder QR.  This
is LAPACK's QR factorization routine \lstinline!DGEQRF! when fast
memory is cache and slow memory is DRAM, and it is ScaLAPACK's
out-of-DRAM QR factorization routine \lstinline!PFDGEQRF! when fast
memory is DRAM and slow memory is disk.  Sequential TSQR transfers
fewer words between slow and fast memory: $2mn$, which is both optimal
and a factor $mn/(4W)$ fewer words than transferred by blocked
Householder QR.  Note that $mn/W$ is how many times larger the matrix
is than the fast memory size $W$.  Furthermore, TSQR requires fewer
messages: $2mn/\widetilde{W}$, which is close to optimal and $O(n)$
times lower than Householder QR.
Table~\ref{tbl:5-seq-caqr-general} compares sequential CAQR and sequential
blocked Householder QR on a general rectangular matrix.  Sequential
CAQR transfers fewer words between slow and fast memory:
$\Theta(mn^2/\sqrt{W})$, which is both optimal and a factor
$\Theta(m/\sqrt{W})$ fewer words transferred than blocked Householder
QR.  Note that $m/\sqrt{W}$ is $\sqrt{m^2/W}$ which is the square root
of how many times larger a square $m \times m$ matrix is than the fast
memory size.  Sequential CAQR also requires fewer messages: $12 mn^2 /
W^{3/2}$, which is optimal.  We note that our analysis of CAQR
applies for any $W$, whereas our analysis of the algorithms in
LAPACK and ScaLAPACK assume that at least 2 columns fit in fast memory,
that is $W \geq 2m$; otherwise they may communicate even more.
Finally, Table~\ref{tbl:6-seq-caqr-square} presents the same comparison for the
special case of a square $n \times n$ matrix.  There again, sequential
CAQR transfers fewer words between slow and fast memory:
$\Theta(n^3/\sqrt{W})$, which is both optimal and a factor
$\Theta(n/\sqrt{W})$ fewer words transferred than blocked Householder
QR.  Sequential CAQR also requires fewer messages: $12 n^3 / W^{3/2}$,
which is optimal.

We expect parallel CAQR to outperform ScaLAPACK's current parallel QR
factorization especially well in the strong scaling regime, i.e., when
the matrix dimensions are constant and the number of processors $P$
varies.  Table \ref{tbl:3-par-caqr-square} shows that the number of
floating-point operations for both algorithms scales as $1/P$, and the
number of words transferred scales as $n^2 \log P / \sqrt{P}$.
However, for ScaLAPACK, the number of messages is proportional to $n
\log^2 P$, whereas for parallel CAQR, the number of messages is
proportional to $\sqrt{P} \log^3 P$, a factor of $n / \sqrt{P}$ fewer
messages.  In either case, the number of messages grows with the
number of processors and also with the data size, if we assume a
limited amount of memory per processor, so reducing communication
costs is important to achieving strong scalability.

We have concentrated on the cases of a homogeneous parallel computer
and a sequential computer with a two-level memory hierarchy. But real
computers are obviously more complicated, combining many levels of
parallelism and memory hierarchy, perhaps heterogeneously. So we have
shown that our parallel and sequential TSQR designs correspond to the
two simplest cases of reduction trees (binary and flat, respectively),
and that different choices of reduction trees will let us optimize
TSQR for more general architectures.

Now we briefly describe related work and our contributions.
The tree-based QR idea itself is not novel (see for example,
\cite{buttari2007class,cunha2002new,golub1988parallel,gunter2005parallel,kurzak2008qr,pothen1989distributed,quintana-orti2008scheduling,rabani2001outcore}),
but we have a number of optimizations and generalizations:
\begin{itemize}
\item Our algorithm can perform almost all its floating-point
  operations using any fast sequential QR factorization routine.  In
  particular, we can achieve significant speedups by invoking Elmroth
  and Gustavson's recursive QR (see
  \cite{elmroth1998new,elmroth2000applying}).

\item We apply TSQR to the parallel factorization of arbitrary
  rectangular matrices in a two-dimensional block cyclic layout.

\item We adapt TSQR to work on general reduction trees.  This
  flexibility lets schedulers overlap communication and computation,
  and minimize communication for more complicated and realistic
  computers with multiple levels of parallelism and memory hierarchy
  (e.g., a system with disk, DRAM, and cache on multiple boards each
  containing one or more multicore chips of different clock speeds,
  along with compute accelerator hardware like GPUs).

\item We prove optimality for both our parallel and sequential algorithms, with
  a 1-D layout for TSQR and 2-D block layout for CAQR, i.e., that they
  minimize bandwidth and latency costs.  This assumes $O(n^3)$
  (non-Strassen-like algorithms), and is done in a Big-Oh sense, sometimes
  modulo polylogarithmic terms.

\item We describe special cases in which existing sequential algorithms 
by Elmroth and Gustavson \cite{elmroth2000applying} and also LAPACK's DGEQRF
attain minimum bandwidth. In particular, with the correct choice of
block size, Elmroth's and Gustavson's RGEQRF algorithm attains minimum
bandwidth and flop count, though not minimum latency.

\item We observe that there are alternative LU algorithms in 
  the literature that attain at least some of these communication
  lower bounds: \cite{grigori2008calu} describes a parallel LU algorithm
  attaining both bandwidth and latency lower bounds, and
  \cite{toledo1997locality} describes a sequential LU algorithm that
  at least attains the bandwidth lower bound.

\item We outline how to extend both algorithms and optimality results
  to certain kinds of hierarchical architectures, either with multiple
  levels of memory hierarchy, or multiple levels of parallelism
  (e.g., where each node in a parallel machine consists of other parallel
   machines, such as multicore).
\end{itemize}

We note that the $Q$ factor is represented as a tree of smaller $Q$
factors, which differs from the traditional layout.  Many previous
authors did not explain in detail how to apply a stored TSQR $Q$
factor, quite possibly because this is not needed for solving least
squares problems.  Adjoining the right-hand side(s) to the matrix $A$,
and taking the QR factorization of the result, requires only the $R$
factor.  Previous authors discuss this optimization.  However, many of
our applications require storing and working with the implicit
representation of the $Q$ factor.  
Our performance models show that applying this tree-structured $Q$
has about the same cost as the traditionally represented $Q$.

\subsection{Outline}\label{SS:intro:outline}

The rest of this report is organized as follows.  
Section~\ref{S:abbrev} first gives a list of terms and abbreviations.  
We then begin the discussion of Tall Skinny QR by Section~\ref{S:motivation},
which motivates the algorithm, giving a variety of applications where
it is used, beyond as a building block for general QR.  
Section~\ref{S:TSQR:algebra} introduces the TSQR algorithm and shows how the
parallel and sequential versions correspond to different reduction or
all-reduction trees.  
After that, Section~\ref{S:reduction}
illustrates how TSQR is actually a reduction, introduces corresponding
terminology, and discusses some design choices.  
Section~\ref{S:TSQR:localQR} shows how the local QR decompositions in TSQR can
be further optimized, including ways that current ScaLAPACK cannot
exploit.  We also explain how to apply the $Q$ factor from TSQR
efficiently, which is needed both for general QR and other
applications.  
Section~\ref{S:perfmodel} explains about our parallel and 
sequential machine models, and what parameters we use to describe them.  
Next, Sections~\ref{S:TSQR:perfcomp} and \ref{S:TSQR:stability}
describe other "tall skinny QR" algorithms, such as CholeskyQR and
Gram-Schmidt, and compare their cost (Section~\ref{S:TSQR:perfcomp})
and numerical stability (Section~\ref{S:TSQR:stability}) to that of
TSQR.  These sections show that TSQR is the only algorithm that
simultaneously minimizes communication and is numerically stable.
Section~\ref{S:TSQR:platforms} describes the platforms used for
testing TSQR, and Section~\ref{S:TSQR:perfres} concludes the
discussion of TSQR proper by describing the TSQR performance results.

Our discussion of CAQR presents both the parallel and the sequential
CAQR algorithms for the QR factorization of general rectangular matrices.  
Section~\ref{S:CAQR} describes the parallel CAQR algorithm
and constructs a performance model.  
Section~\ref{S:CAQR-seq} does the same for sequential CAQR.  
Subsection~\ref{sec:seq_qr_other} analyzes other sequential
QR algorithms including those of Elmroth and Gustavson.
Next, Section~\ref{S:CAQR-counts} compares
the performance of parallel CAQR and ScaLAPACK's \lstinline!PDGEQRF!,
showing CAQR to be superior, for the same choices of block sizes and
data layout parameters, as well as when these parameters are chosen
optimally and independently for CAQR and \lstinline!PDGEQRF!.  
After that, Section~\ref{S:CAQR:perfest} presents performance predictions
comparing CAQR to \lstinline!PDGEQRF!. Future work includes actual
implementation and measurements.

The next two sections in the body of the text concern theoretical
results about CAQR and other parallel and sequential QR factorizations.  
Section~\ref{S:lowerbounds} describes how to extend
known lower bounds on communication for matrix multiplication to QR,
and shows that these are attained (modulo polylogarithmic factors) by
TSQR and CAQR.  
Section~\ref{S:limits-to-par} reviews known lower
bounds on parallelism for QR, using a PRAM model of parallel
computation.

The final section, Section~\ref{S:hierarchies} briefly outlines how 
to extend the algorithms and optimality results to hierarchical architectures, 
either with several levels of memory hierarchy, or several levels
of parallelism.

The Appendices provide details of operation counts and other results
summarized in previous sections.  Appendix~\ref{S:localQR-flops}
presents flop counts for optimizations of local QR decompositions
described in Section~\ref{S:localQR}.  
Appendices~\ref{S:TSQR-seq-detailed}, \ref{S:CAQR-seq-detailed},
\ref{S:TSQR-par-detailed}, and \ref{S:CAQR-par-detailed} give details
of performance models for sequential TSQR, sequential CAQR, parallel
TSQR and parallel CAQR, respectively.  
Appendix~\ref{S:PFDGEQRF} models sequential QR based on ScaLAPACK's out-of-DRAM
routine \lstinline!PFDGEQRF!.
Finally, Appendix~\ref{S:CommLowerBoundsFromCalculus} 
proves communication lower bounds
needed in Section~\ref{S:lowerbounds}.

\subsection{Future work}

Implementations of sequential and parallel CAQR are currently
underway.  Optimization of the TSQR reduction tree for more general,
practical architectures (such as multicore, multisocket, or GPUs) is
future work, as well as optimization of the rest of CAQR to the
most general architectures, with proofs of optimality.

It is natural to ask to how much of dense linear algebra one
can extend the results of this paper, that is finding algorithms that
attain communication lower bounds.
In the case of parallel LU with
pivoting, refer to the technical report by Grigori, Demmel, and Xiang
\cite{grigori2008calu}, and in the case of sequential LU, refer to the
paper by Toledo \cite{toledo1997locality} (at least for minimizing bandwidth).
More broadly, we hope to extend the
results of this paper to the rest of linear algebra, including 
two-sided factorizations (such as reduction to symmetric tridiagonal,
bidiagonal, or (generalized) upper Hessenberg forms).  Once 
a matrix is symmetric tridiagonal (or bidiagonal) and so takes
little memory,  fast algorithms for the eigenproblem (or SVD)
are available. Most challenging is likely to be find
eigenvalues of a matrix in upper Hessenberg form (or of 
a matrix pencil).


\section{List of terms and abbreviations}\label{S:abbrev}

\begin{description}
  \item[alpha-beta model] A simple model for communication time,
    involving a latency parameter $\alpha$ and an inverse bandwidth
    parameter $\beta$:  the time to transfer a single message
    containing $n$ words is $\alpha + \beta n$.

  \item[CAQR] Communication-Avoiding QR -- a parallel and/or
    \emph{explicitly swapping} QR factorization algorithm, intended
    for input matrices of general shape.  Invokes \emph{TSQR} for
    panel factorizations.

  \item[CholeskyQR] A fast but numerically unstable QR factorization
    algorithm for tall and skinny matrices, based on the Cholesky
    factorization of $A^T A$.

  \item[\texttt{DGEQRF}] LAPACK QR factorization routine for general
    dense matrices of double-precision floating-point numbers.  May or
    may not exploit shared-memory parallelism via a multithreaded BLAS
    implementation.

  \item[GPU] Graphics processing unit.

  \item[Explicitly swapping] Refers to algorithms explicitly written
    to save space in one level of the memory hierarchy (``fast
    memory'') by using the next level (``slow memory'') as swap space.
    Explicitly swapping algorithms can solve problems too large to fit
    in fast memory.  Special cases include \emph{out-of-DRAM} (a.k.a.\
    out-of-core), \emph{out-of-cache} (which is a performance
    optimization that manages cache space explicitly in the
    algorithm), and algorithms written for processors with
    non-cache-coherent local scratch memory and global DRAM (such as
    Cell).

  \item[Flash drive] A persistent storage device that uses nonvolatile
    flash memory, rather than the spinning magnetic disks used in hard
    drives.  These are increasingly being used as replacements for
    traditional hard disks for certain applications.  Flash drives are
    a specific kind of solid-state drive (SSD), which uses solid-state
    (not liquid, gas, or plasma) electronics with no moving parts to
    store data.

  \item[Local store] A user-managed storage area which functions like
    a cache (in that it is smaller and faster than main memory), but
    has no hardware support for cache coherency.

  \item[Out-of-cache] Refers to algorithms explicitly written to save
    space in cache (or local store), by using the next larger level of
    cache (or local store), or main memory (DRAM), as swap space.

  \item[Out-of-DRAM] Refers to algorithms explicitly written to save
    space in main memory (DRAM), by using disk as swap space.
    (``Core'' used to mean ``main memory,'' as main memories were once
    constructed of many small solenoid cores.)  See \emph{explicitly
      swapping}.

  \item[\texttt{PDGEQRF}] ScaLAPACK parallel QR factorization
    routine for general dense matrices of double-precision
    floating-point numbers.

  \item[\texttt{PFDGEQRF}] ScaLAPACK parallel out-of-core QR
    factorization routine for general dense matrices of
    double-precision floating-point numbers.

  \item[TSQR] Tall Skinny QR -- our reduction-based QR factorization
    algorithm, intended for ``tall and skinny'' input matrices (i.e.,
    those with many more rows than columns).
\end{description}


\section{Motivation for TSQR}\label{S:motivation}

\subsection{Block iterative methods}

Block iterative methods frequently compute the QR factorization of a
tall and skinny dense matrix.  This includes algorithms for solving
linear systems $Ax = B$ with multiple right-hand sides (such as
variants of GMRES, QMR, or CG
\cite{vital:phdthesis:90,Freund:1997:BQA,oleary:80}), as well as block
iterative eigensolvers (for a summary of such methods, see
\cite{templatesEigenBai,templatesEigenLehoucq}).  Many of these
methods have widely used implementations, on which a large community
of scientists and engineers depends for their computational tasks.
Examples include TRLAN (Thick Restart Lanczos), BLZPACK (Block
Lanczos), Anasazi (various block methods), and PRIMME (block
Jacobi-Davidson methods)
\cite{TRLANwebpage,BLZPACKwebpage,BLOPEXwebpage,irbleigs,TRILINOSwebpage,PRIMMEwebpage}.
Eigenvalue computation is particularly sensitive to the accuracy of
the orthogonalization; two recent papers suggest that large-scale
eigenvalue applications require a stable QR factorization
\cite{lehoucqORTH,andrewORTH}.

\subsection{$s$-step Krylov methods}

Recent research has reawakened an interest in alternate formulations
of Krylov subspace methods, called \emph{$s$-step Krylov methods}, in
which some number $s$ steps of the algorithm are performed all at
once, in order to reduce communication.  Demmel et al.\ review the
existing literature and discuss new advances in this area
\cite{demmel2008comm}.  Such a method begins with an $n \times n$
matrix $A$ and a starting vector $v$, and generates some basis for the
Krylov subspace $\Span\{v$, $Av$, $A^2 v$, $\dots$, $A^s v\}$, using a
small number of communication steps that is independent of $s$.  Then,
a QR factorization is used to orthogonalize the basis vectors.

The goal of combining $s$ steps into one is to leverage existing basis
generation algorithms that reduce the number of messages and/or the
volume of communication between different levels of the memory
hierarchy and/or different processors.  These algorithms make the
resulting number of messages independent of $s$, rather than growing
with $s$ (as in standard Krylov methods).  However, this means that
the QR factorization is now the communications bottleneck, at least in
the parallel case: the current \lstinline!PDGEQRF!  algorithm in
ScaLAPACK takes $2s \log_2 P$ messages (in which $P$ is the
number of processors), compared to $\log_2 P$ messages for TSQR.
Numerical stability considerations limit $s$, so that it is
essentially a constant with respect to the matrix size $m$.
Furthermore, a stable QR factorization is necessary in order to
restrict the loss of stability caused by generating $s$ steps of the
basis without intermediate orthogonalization.  This is an ideal
application for TSQR, and in fact inspired its (re-)discovery.

\subsection{Panel factorization in general QR}

Householder QR decompositions of tall and skinny matrices also
comprise the panel factorization step for typical QR factorizations of
matrices in a more general, two-dimensional layout.  This includes the
current parallel QR factorization routine \lstinline!PDGEQRF! in
ScaLAPACK, as well as ScaLAPACK's out-of-DRAM QR factorization
\lstinline!PFDGEQRF!.  Both algorithms use a standard column-based
Householder QR for the panel factorizations, but in the parallel case
this is a latency bottleneck, and in the out-of-DRAM case it is a
bandwidth bottleneck.  Replacing the existing panel factorization with
TSQR would reduce this cost by a factor equal to the number of columns
in a panel, thus removing the bottleneck.  TSQR requires more
floating-point operations, though some of this computation can be
overlapped with communication.  Section \ref{S:CAQR} will discuss the
advantages of this approach in detail.



\section{TSQR matrix algebra}\label{S:TSQR:algebra}

In this section, we illustrate the insight behind the TSQR algorithm.
TSQR uses a reduction-like operation to compute the QR factorization
of an $m \times n$ matrix $A$, stored in a 1-D block row
layout.\footnote{The ScaLAPACK Users' Guide has a good explanation of
  1-D and 2-D block and block cyclic layouts of dense matrices
  \cite{scalapackusersguide}.  In particular, refer to the section
  entitled ``Details of Example Program \#1.''}  We begin with
parallel TSQR on a binary tree of four processors ($P = 4$), and later
show sequential TSQR on a linear tree with four blocks.

\subsection{Parallel TSQR on a binary tree}

\begin{figure}
  \begin{center}
    \includegraphics[angle=0,scale=0.7]{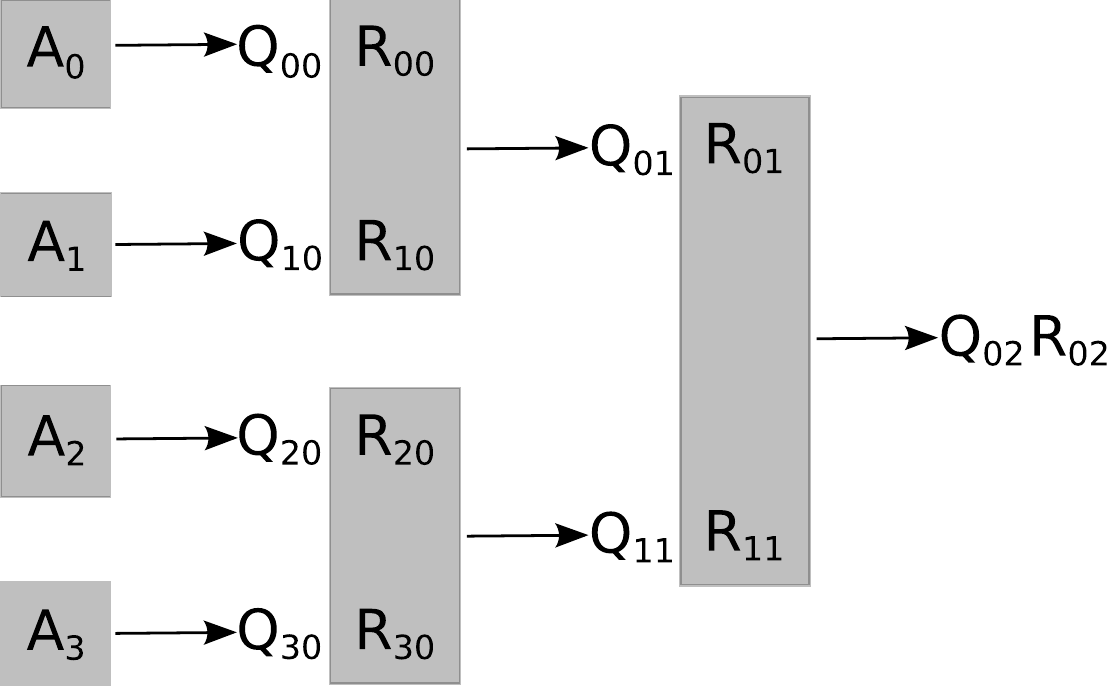}
    \caption{Execution of the parallel TSQR factorization on a binary
      tree of four processors.  The gray boxes indicate where local QR
      factorizations take place.  The $Q$ and $R$ factors each have
      two subscripts: the first is the sequence number within that stage,
      and the second is the stage number.}
    \label{Fi:TSQR:algebra:par4}
  \end{center}
\end{figure}

The basic idea of using a reduction on a binary tree to compute a tall
skinny QR factorization has been rediscovered more than once (see
e.g., \cite{cunha2002new,pothen1989distributed}).  (TSQR was also
suggested by Golub et al.\ \cite{golub1988parallel}, but they did not
reduce the number of messages from $n \log P$ to $\log P$.)  We repeat
it here in order to show its generalization to a whole space of
algorithms.  First, we decompose the $m \times n$ matrix $A$ into four
$m/4 \times n$ block rows:
\[
A = 
\begin{pmatrix}
A_0 \\
A_1 \\
A_2 \\
A_3 \\
\end{pmatrix}.
\]
Then, we independently compute the QR factorization of each block row:
\[
\begin{pmatrix}
A_0 \\
A_1 \\
A_2 \\
A_3 \\
\end{pmatrix}
=
\begin{pmatrix}
  Q_{00} R_{00} \\
  Q_{10} R_{10} \\
  Q_{20} R_{20} \\
  Q_{30} R_{30} \\
\end{pmatrix}.
\]
This is ``stage 0'' of the computation, hence the second subscript 0
of the $Q$ and $R$ factors.  The first subscript indicates the block
index at that stage.  (Abstractly, we use the Fortran convention that
the first index changes ``more frequently'' than the second index.)
Stage 0 operates on the $P = 4$ leaves of the tree.  We can write this
decomposition instead as a block diagonal orthogonal matrix times a
column of blocks:
\[
A =
\begin{pmatrix}
  Q_{00} R_{00} \\
  Q_{10} R_{10} \\
  Q_{20} R_{20} \\
  Q_{30} R_{30} \\
\end{pmatrix}
=
\left(
  \begin{array}{c | c | c | c}
    Q_{00} & & & \\ \hline
    & Q_{10} & & \\ \hline
    & & Q_{20} & \\ \hline
    & & & Q_{30} \\ 
  \end{array}
\right)
\cdot
\begin{pmatrix}
  R_{00} \\ \hline
  R_{10} \\ \hline
  R_{20} \\ \hline
  R_{30} \\ 
\end{pmatrix},
\]
although we do not have to store it this way.  After this stage 0,
there are $P = 4$ of the $R$ factors.  We group them into successive
pairs $R_{i,0}$ and $R_{i+1,0}$, and do the QR factorizations of
grouped pairs in parallel:
\[
\begin{pmatrix}
  R_{00} \\
  R_{10} \\ \hline
  R_{20} \\ 
  R_{30} \\
\end{pmatrix}
= 
\begin{pmatrix}
  \begin{pmatrix} 
    R_{00} \\
    R_{10} \\
  \end{pmatrix} \\ \hline
  \begin{pmatrix}
    R_{20} \\
    R_{30} \\
  \end{pmatrix}
\end{pmatrix}
=
\begin{pmatrix}
  Q_{01} R_{01} \\ \hline
  Q_{11} R_{11} \\ 
\end{pmatrix}.
\]
As before, we can rewrite the last term as a block diagonal orthogonal
matrix times a column of blocks:
\[
\begin{pmatrix}
  Q_{01} R_{01} \\ \hline
  Q_{11} R_{11} \\ 
\end{pmatrix}
=
\left(
  \begin{array}{c | c}
    Q_{01} &        \\ \hline
          &  Q_{11} \\ 
  \end{array}
\right)
\cdot
\begin{pmatrix}
  R_{01} \\ \hline
  R_{11} \\
\end{pmatrix}.
\]
This is stage 1, as the second subscript of the $Q$ and $R$ factors
indicates.  We iteratively perform stages until there is only one $R$
factor left, which is the root of the tree:
\[
\begin{pmatrix}
  R_{01} \\
  R_{11} \\
\end{pmatrix}
=
Q_{02} R_{02}.
\]
Equation \eqref{eq:TSQR:algebra:par4:final} shows the whole 
factorization:
\begin{equation}\label{eq:TSQR:algebra:par4:final}
A = 
\begin{pmatrix}
A_0 \\
A_1 \\
A_2 \\
A_3 \\
\end{pmatrix}
=
\left(
  \begin{array}{c | c | c | c}
    Q_{00} & & & \\ \hline
    & Q_{10} & & \\ \hline
    & & Q_{20} & \\ \hline
    & & & Q_{30} \\ 
  \end{array}
\right)
\cdot
\left(
  \begin{array}{c | c}
    Q_{01} &       \\ \hline
          & Q_{11} \\ 
  \end{array}
\right)
\cdot
Q_{02} \cdot R_{02},
\end{equation}
in which the product of the first three matrices has orthogonal
columns, since each of these three matrices does.  Note the binary
tree structure in the nested pairs of $R$ factors.

Figure \ref{Fi:TSQR:algebra:par4} illustrates the binary tree on which
the above factorization executes.  Gray boxes highlight where local QR
factorizations take place.  By ``local,'' we refer to a factorization
performed by any one processor at one node of the tree; it may involve
one or more than one block row.  If we were to compute all the above
$Q$ factors explicitly as square matrices, each of the $Q_{i0}$ would
be $m/P \times m/P$, and $Q_{ij}$ for $j > 0$ would be $2n \times 2n$.
The final $R$ factor would be upper triangular and $m \times n$, with
$m - n$ rows of zeros.  In a ``thin QR'' factorization, in which the
final $Q$ factor has the same dimensions as $A$, the final $R$ factor
would be upper triangular and $n \times n$.  In practice, we prefer to
store all the local $Q$ factors implicitly until the factorization is
complete.  In that case, the implicit representation of $Q_{i0}$ fits
in an $m/P \times n$ lower triangular matrix, and the implicit
representation of $Q_{ij}$ (for $j > 0$) fits in an $n \times n$ lower
triangular matrix (due to optimizations that will be discussed in
Section \ref{S:TSQR:localQR}).

Note that the maximum per-processor memory requirement is $\max\{mn/P,
n^2 + O(n)\}$, since any one processor need only factor two $n \times n$
upper triangular matrices at once, or a single $m/P \times n$ matrix.

\subsection{Sequential TSQR on a flat tree}

\begin{figure}
  \begin{center}
    \includegraphics[angle=0,scale=0.7]{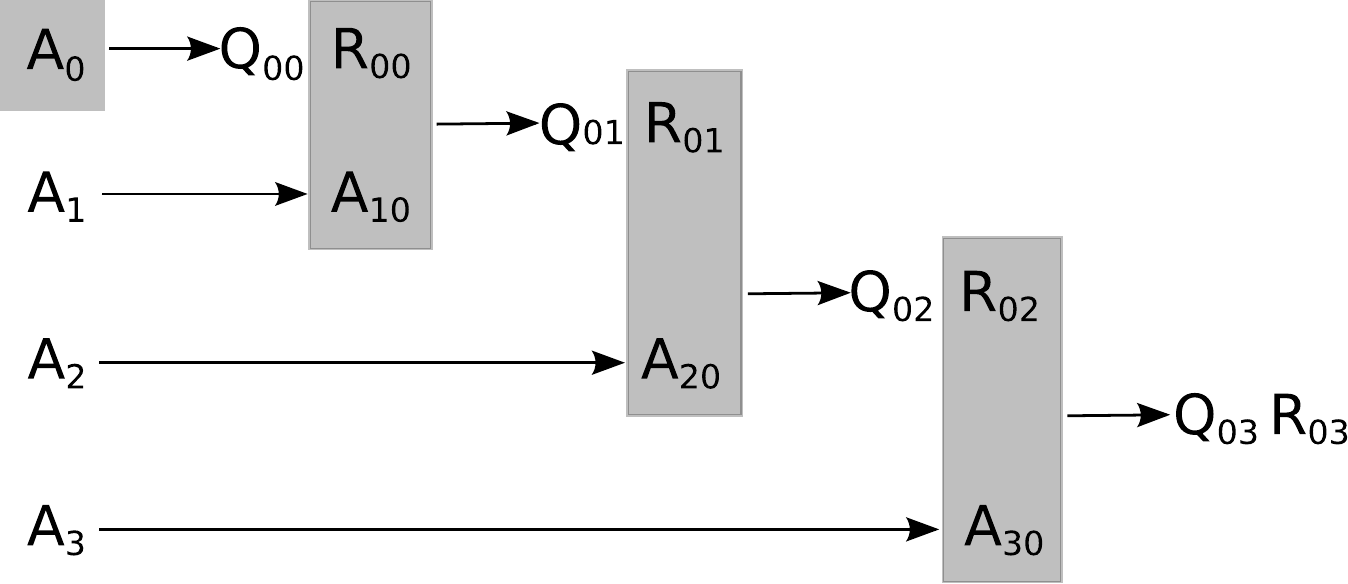}
    \caption{Execution of the sequential TSQR factorization on a flat
      tree with four submatrices.  The gray boxes indicate where local
      QR factorizations take place The $Q$ and $R$ factors each have
      two subscripts: the first is the sequence number for that stage,
      and the second is the stage number.}
    \label{Fi:TSQR:algebra:seq4}
  \end{center}
\end{figure}

Sequential TSQR uses a similar factorization process, but with a
``flat tree'' (a linear chain).  It may also handle the leaf nodes of
the tree slightly differently, as we will show below.  Again, the
basic idea is not new; see e.g.,
\cite{buttari2007class,buttari2007parallel,gunter2005parallel,kurzak2008qr,quintana-orti2008scheduling,rabani2001outcore}.
(Some authors (e.g.,
\cite{buttari2007class,kurzak2008qr,quintana-orti2008scheduling})
refer to sequential TSQR as ``tiled QR.''  We use the phrase
``sequential TSQR'' because both our parallel and sequential
algorithms could be said to use tiles.)  In particular, Gunter and van
de Geijn develop a parallel out-of-DRAM QR factorization algorithm
that uses a flat tree for the panel factorizations
\cite{gunter2005parallel}.  Buttari et al.\ suggest using a QR
factorization of this type to improve performance of parallel QR on
commodity multicore processors \cite{buttari2007class}.  Quintana-Orti
et al.\ develop two variations on block QR factorization algorithms,
and use them with a dynamic task scheduling system to parallelize the
QR factorization on shared-memory machines
\cite{quintana-orti2008scheduling}.  Kurzak and Dongarra use similar
algorithms, but with static task scheduling, to parallelize the QR
factorization on Cell processors \cite{kurzak2008qr}.
The reason these authors use what we call sequential TSQR in a parallel
context ...

We will show that the basic idea of sequential TSQR fits into the same
general framework as the parallel QR decomposition illustrated above,
and also how this generalization expands the tuning space of QR
factorization algorithms.  In addition, we will develop detailed
performance models of sequential TSQR and the current sequential QR
factorization implemented in LAPACK.

We start with the same block row decomposition as with parallel TSQR
above:
\[
A = 
\begin{pmatrix}
A_0 \\
A_1 \\
A_2 \\
A_3 \\
\end{pmatrix}
\]
but begin with a QR factorization of $A_0$, rather than of all the
block rows:
\[
\begin{pmatrix}
A_0 \\
A_1 \\
A_2 \\
A_3 \\
\end{pmatrix}
=
\begin{pmatrix}
Q_{00} R_{00} \\
A_1 \\
A_2 \\
A_3 \\
\end{pmatrix}.
\]
This is ``stage 0'' of the computation, hence the second subscript 0
of the $Q$ and $R$ factor.  We retain the first subscript for
generality, though in this example it is always zero.  We can write
this decomposition instead as a block diagonal matrix times a column
of blocks:
\[
\begin{pmatrix}
  Q_{00} R_{00} \\
  A_1 \\
  A_2 \\
  A_3 \\
\end{pmatrix}
=
\left(
  \begin{array}{c | c | c | c}
    Q_{00} & & & \\ \hline
    & I & & \\ \hline
    & & I & \\ \hline
    & & & I \\ 
  \end{array}
\right)
\cdot
\begin{pmatrix}
  R_{00} \\ \hline
  A_1    \\ \hline
  A_2    \\ \hline
  A_3    \\
\end{pmatrix}.
\]
We then combine $R_{00}$ and $A_1$ using a QR factorization:
\[
\begin{pmatrix}
R_{00} \\
A_1 \\
A_2 \\
A_3 \\
\end{pmatrix}
=
\begin{pmatrix}
R_{00} \\
A_1 \\ \hline
A_2 \\
A_3 \\
\end{pmatrix}
=
\begin{pmatrix}
  Q_{01} R_{01} \\ \hline
  A_2 \\
  A_3 \\
\end{pmatrix}
\]
This can be rewritten as a block diagonal matrix times a column of blocks:
\[
\begin{pmatrix}
  Q_{01} R_{01} \\ \hline
  A_2 \\
  A_3 \\
\end{pmatrix}
=
\left(
\begin{array}{c | c | c}
  Q_{01} &   &    \\ \hline
        & I &    \\ \hline
        &   & I  \\
\end{array}
\right)
\cdot
\begin{pmatrix}
  R_{01} \\ \hline
  A_2   \\ \hline
  A_3   \\ 
\end{pmatrix}.
\]
We continue this process until we run out of $A_i$ factors.  The
resulting factorization has the following structure:
\begin{equation}\label{eq:TSQR:algebra:seq4:final}
\begin{pmatrix}
A_0 \\
A_1 \\
A_2 \\
A_3 \\
\end{pmatrix}
=
\left(
  \begin{array}{c | c | c | c}
    Q_{00} & & & \\ \hline
    & I & & \\ \hline
    & & I & \\ \hline
    & & & I \\ 
  \end{array}
\right)
\cdot
\left(
\begin{array}{c | c | c}
  Q_{01} &   &    \\ \hline
        & I &    \\ \hline
        &   & I  \\
\end{array}
\right)
\cdot
\left(
  \begin{array}{c | c | c}
    I & & \\ \hline
    & Q_{02} & \\ \hline
    & & I \\
  \end{array}
\right)
\cdot
\left(
  \begin{array}{c | c | c}
    I & & \\ \hline
    & I & \\ \hline
    & & Q_{03} \\ 
  \end{array}
\right)
R_{30}.
\end{equation}
Here, the $A_i$ blocks are $m/P \times n$.  If we were to compute all
the above $Q$ factors explicitly as square matrices, then $Q_{00}$
would be $m/P \times m/P$ and $Q_{0j}$ for $j > 0$ would be $2m/P
\times 2m/P$.  The above $I$ factors would be $m/P \times m/P$.  The
final $R$ factor, as in the parallel case, would be upper triangular
and $m \times n$, with $m - n$ rows of zeros.  In a ``thin QR''
factorization, in which the final $Q$ factor has the same dimensions
as $A$, the final $R$ factor would be upper triangular and $n \times
n$.  In practice, we prefer to store all the local $Q$ factors
implicitly until the factorization is complete.  In that case, the
implicit representation of $Q_{00}$ fits in an $m/P \times n$ lower
triangular matrix, and the implicit representation of $Q_{0j}$ (for $j
> 0$) fits in an $m/P \times n$ lower triangular matrix as well (due
to optimizations that will be discussed in Section
\ref{S:TSQR:localQR}).

Figure \ref{Fi:TSQR:algebra:seq4} illustrates the flat tree on which
the above factorization executes.  Gray boxes highlight where
``local'' QR factorizations take place.  

The sequential algorithm differs from the parallel one in that it does
not factor the individual blocks of the input matrix $A$, excepting
$A_0$.  This is because in the sequential case, the input matrix has
not yet been loaded into working memory.  In the fully parallel case,
each block of $A$ resides in some processor's working memory.  It then
pays to factor all the blocks before combining them, as this reduces
the volume of communication (only the triangular $R$ factors need to
be exchanged) and reduces the amount of arithmetic performed at the
next level of the tree.  In contrast, the sequential algorithm never
writes out the intermediate $R$ factors, so it does not need to
convert the individual $A_i$ into upper triangular factors.  Factoring
each $A_i$ separately would require writing out an additional $Q$
factor for each block of $A$.  It would also add another level to the
tree, corresponding to the first block $A_0$.  

Note that the maximum per-processor memory requirement is $mn/P +
n^2/2 + O(n)$, since only an $m/P \times n$ block and an $n \times n$
upper triangular block reside in fast memory at one time.  We could
save some fast memory by factoring each $A_i$ block separately before
combining it with the next block's $R$ factor, as long as each block's
$Q$ and $R$ factors are written back to slow memory before the next
block is loaded.  One would then only need to fit no more than two $n
\times n$ upper triangular factors in fast memory at once.  However,
this would result in more writes, as each $R$ factor (except the last)
would need to be written to slow memory and read back into fact
memory, rather than just left in fast memory for the next step.

In both the parallel and sequential algorithms, a vector or matrix is
multiplied by $Q$ or $Q^T$ by using the implicit representation of the
$Q$ factor, as shown in Equation \eqref{eq:TSQR:algebra:par4:final}
for the parallel case, and Equation \eqref{eq:TSQR:algebra:seq4:final}
for the sequential case.  This is analogous to using the Householder
vectors computed by Householder QR as an implicit representation of
the $Q$ factor.

\subsection{TSQR on general trees}
\label{SS:TSQR:GeneralTrees}

\begin{figure}
  \begin{center}
    \includegraphics[angle=0,scale=0.35]{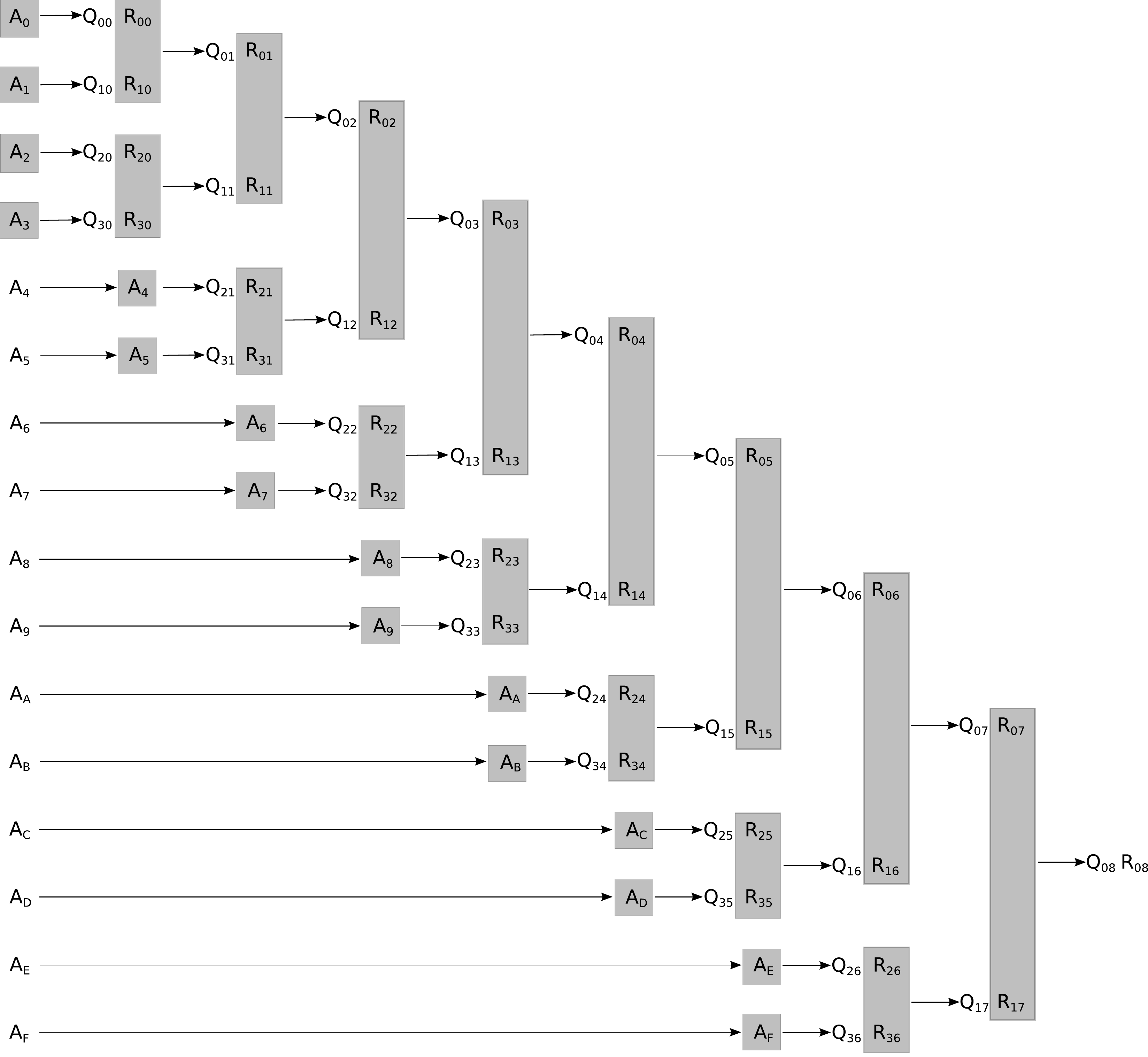}
    \caption{Execution of a hybrid parallel / out-of-core TSQR
      factorization.  The matrix has 16 blocks, and four processors
      can execute local QR factorizations simultaneously.  The gray
      boxes indicate where local QR factorizations take place.  We
      number the blocks of the input matrix $A$ in hexadecimal to save
      space (which means that the subscript letter A is the number
      $10_{10}$, but the non-subscript letter $A$ is a matrix block).
      The $Q$ and $R$ factors each have two subscripts: the first is
      the sequence number for that stage, and the second is the stage
      number.}
    \label{Fi:TSQR:algebra:pooc16}
  \end{center}
\end{figure}

The above two algorithms are extreme points in a large set of possible
QR factorization methods, parametrized by the tree structure.  Our
version of TSQR is novel because it works on any tree.  In general,
the optimal tree may depend on both the architecture and the matrix
dimensions.  This is because TSQR is a reduction (as we will discuss
further in Section \ref{S:reduction}).  Trees of types other than
binary often result in better reduction performance, depending on the
architecture (see e.g., \cite{nishtala2008performance}).  Throughout
this paper, we discuss two examples -- the binary tree and the flat
tree -- as easy extremes for illustration.  We will show that the
binary tree minimizes the number of stages and messages in the
parallel case, and that the flat tree minimizes the number and volume
of input matrix reads and writes in the sequential case.  Section
\ref{S:reduction} shows how to perform TSQR on any tree.  Methods for
finding the best tree in the case of TSQR are future work.
Nevertheless, we can identify two regimes in which a ``nonstandard''
tree could improve performance significantly: parallel memory-limited
CPUs, and large distributed-memory supercomputers.

The advent of desktop and even laptop multicore processors suggests a
revival of parallel out-of-DRAM algorithms, for solving cluster-sized
problems while saving power and avoiding the hassle of debugging on a
cluster.  TSQR could execute efficiently on a parallel memory-limited
device if a sequential flat tree were used to bring blocks into
memory, and a parallel tree (with a structure that reflects the
multicore memory hierarchy) were used to factor the blocks.  Figure
\ref{Fi:TSQR:algebra:pooc16} shows an example with 16 blocks executing
on four processors, in which the factorizations are pipelined for
maximum utilization of the processors.  The algorithm itself needs no
modification, since the tree structure itself encodes the pipelining.
This is, we believe, a novel extension of the parallel out-of-core QR
factorization of Gunter et al.\ \cite{gunter2005parallel}.

TSQR's choice of tree shape can also be optimized for modern
supercomputers.  A tree with different branching factors at different
levels could naturally accommodate the heterogeneous communication
network of a cluster of multicores.  The subtrees at the lowest level
may have the same branching factor as the number of cores per node (or
per socket, for a multisocket shared-memory architecture).

Note that the maximum per-processor memory requirement of all TSQR
variations is bounded above by
\[
\frac{q n(n+1)}{2} + \frac{mn}{P},
\]
in which $q$ is the maximum branching factor in the tree.


\section{TSQR as a reduction}\label{S:reduction}

Section \ref{S:TSQR:algebra} explained the algebra of the TSQR
factorization.  It outlined how to reorganize the parallel QR
factorization as a tree-structured computation, in which groups of
neighboring processors combine their $R$ factors, perform (possibly
redundant) QR factorizations, and continue the process by
communicating their $R$ factors to the next set of neighbors.
Sequential TSQR works in a similar way, except that communication
consists of moving matrix factors between slow and fast memory.  This
tree structure uses the same pattern of communication found in a
reduction or all-reduction.  Thus, effective optimization of TSQR
requires understanding these operations.

\subsection{Reductions and all-reductions}\label{S:reduction:defs}

Reductions and all-reductions are operations that take a collection as
input, and combine the collection using some (ideally) associative
function into a single item.  The result is a function of all the
items in the input.  Usually, one speaks of (all-) reductions in the
parallel case, where ownership of the input collection is distributed
across some number $P$ of processors.  A reduction leaves the final
result on exactly one of the $P$ processors; an all-reduction leaves a
copy of the final result on all the processors.  See, for example,
\cite{gropp1999using}.

In the sequential case, there is an analogous operation.  Imagine that
there are $P$ ``virtual processors.''  To each one is assigned a
certain amount of fast memory.  Virtual processors communicate by
sending messages via slow memory, just as the ``real processors'' in
the parallel case communicate via the (relatively slow) network. Each
virtual processor owns a particular subset of the input data, just as
each real processor does in a parallel implementation.  A virtual
processor can read any other virtual processor's subset by reading
from slow memory (this is a ``receive'').  It can also write some data
to slow memory (a ``send''), for another virtual processor to read.
We can run programs for this virtual parallel machine on an actual
machine with only one processor and its associated fast memory by
scheduling the virtual processors' tasks on the real processor(s) in a
way that respects task dependencies.  Note that all-reductions and
reductions produce the same result when there is only one actual
processor, because if the final result ends up in fast memory on any
of the virtual processors, it is also in fast memory on the one actual
processor.

The ``virtual processors'' argument may also have practical use when
implementing (all-) reductions on clusters of SMPs or vector
processors, multicore out-of-core, or some other combination
consisting of tightly-coupled parallel units with slow communication
links between the units.  A good mapping of virtual processors to real
processors, along with the right scheduling of the ``virtual''
algorithm on the real machine, can exploit multiple levels of
parallelism and the memory hierarchy.

\subsection{(All-) reduction trees}
\lstset{language=Lisp}

\begin{figure}
  \centering
  \scalebox{0.4}{\includegraphics{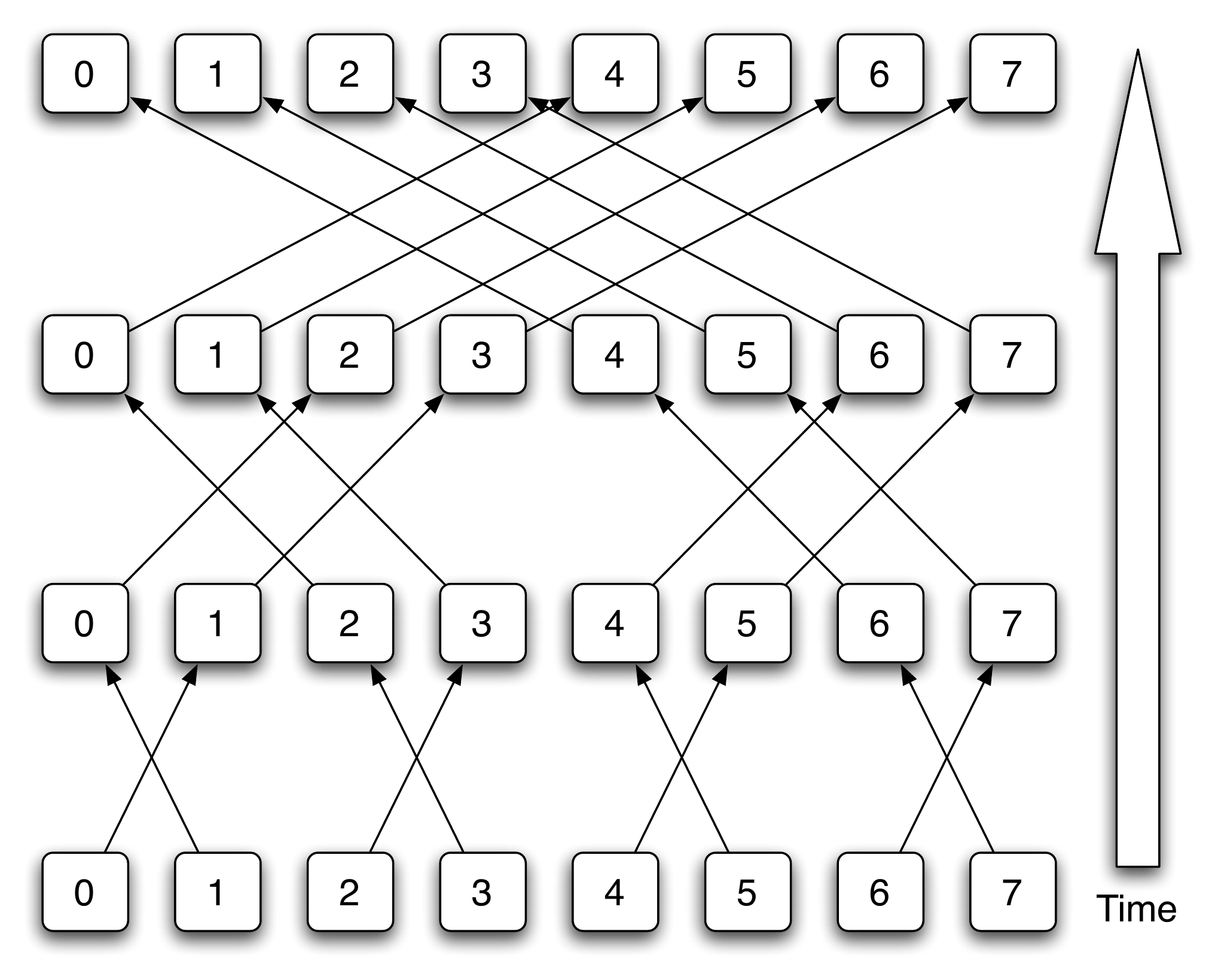}}
  \caption{Diagram of a parallel butterfly all-reduction on a binary
    tree of eight processors.  Each arrow represents a message from
    one processor to another.  Time moves
    upwards.}\label{Fi:allreduction:butterfly}
\end{figure}

Reductions and all-reductions are performed on directed trees.  In a
reduction, each node represents a processor, and each edge a message
passed from one processor to another.  All-reductions have two
different implementation strategies:
\begin{itemize}
\item ``Reduce-broadcast'': Perform a standard reduction to one
  processor, followed by a broadcast (a reduction run backwards) of
  the result to all processors.

\item ``Butterfly'' method, with a communication pattern like that of
  a fast Fourier transform.
\end{itemize}  
The butterfly method uses a tree with the following recursive
structure:
\begin{itemize}
\item Each leaf node corresponds to a single processor.

\item Each interior node is an ordered tuple whose members are the
  node's children.

\item Each edge from a child to a parent represents a complete
  exchange of information between all individual processors at the
  same positions in the sibling tuples.
\end{itemize}
We call the processors that communicate at a particular stage
\emph{neighbors}.  For example, in a a binary tree with eight
processors numbered 0 to 7, processors 0 and 1 are neighbors at the
first stage, processors 0 and 2 are neighbors at the second stage, and
processors 0 and 4 are neighbors at the third (and final) stage.  At
any stage, each neighbor sends its current reduction value to all the
other neighbors.  The neighbors combine the values redundantly, and
the all-reduction continues.  Figure \ref{Fi:allreduction:butterfly}
illustrates this process.  The butterfly all-reduction can be extended
to any number of processors, not just powers of two.

The reduce-broadcast implementation requires about twice as many
stages as the butterfly pattern (in the case of a binary tree) and
thus as much as twice the latency.  However, it reduces the total
number of messages communicated per level of the tree (not just the
messages on the critical path).  In the case of a binary tree,
reduce-broadcast requires at most $P/2$ messages at any one level, and
$P \log(P) / 2$ total messages.  A butterfly always generates $P$
messages at every level, and requires $P \log(P)$ total messages.  The
choice between reduce-broadcast and butterfly depends on the
properties of the communication network.

\subsection{TSQR-specific (all-) reduction requirements}\label{SS:reduction:requirements}

TSQR uses an (all-) reduction communication pattern, but has
requirements that differ from the standard (all-) reduction.  For
example, if the $Q$ factor is desired, then TSQR must store
intermediate results (the local $Q$ factor from each level's
computation with neighbors) at interior nodes of the tree.  This
requires reifying and preserving the (all-) reduction tree for later
invocation by users.  Typical (all-) reduction interfaces, such as
those provided by MPI or OpenMP, do not allow this (see e.g.,
\cite{gropp1999using}).  They may not even guarantee that the same
tree will be used upon successive invocations of the same (all-)
reduction operation, or that the inputs to a node of the (all-)
reduction tree will always be in the same order.


\section{Optimizations for local QR factorizations}\label{S:TSQR:localQR} 
\label{S:localQR}

Although TSQR achieves its performance gains because it optimizes
communication, the local QR factorizations lie along the critical path
of the algorithm.  The parallel cluster benchmark results in Section
\ref{S:TSQR:perfres} show that optimizing the local QR factorizations
can improve performance significantly.  In this section, we outline a
few of these optimizations, and hint at how they affect the
formulation of the general CAQR algorithm in Section \ref{S:CAQR}.

\subsection{Structured QR factorizations}\label{SS:TSQR:localQR:structured}

Many of the inputs to the local QR factorizations have a particular
structure.  In the parallel case, they are vertical stacks of $n
\times n$ upper triangular matrices, and in the sequential case, at
least one of the blocks is upper triangular.  In this section, we show
how to modify a standard dense Householder QR factorization in order
to exploit this structure.  This can save a factor of $5 \times$ flops
and (at least) $3 \times$ storage, in the parallel case, and a factor
of $2 \times$ flops and (up to) $2 \times$ storage in the sequential
case.  We also show how to perform the trailing matrix update with
these \emph{structured QR factorizations}, as it will be useful for
Section \ref{S:CAQR}.

Suppose that we have two upper triangular matrices $R_0$ and $R_1$,
each of size $5 \times 5$.  (The notation here is generic and not
meant to correspond to a specific stage of TSQR.  This is extended
easily enough to the case of $q$ upper triangular matrices, for $q =
2, 3, \dots$.)  Then, we can write their vertical concatenation as
follows, in which an $x$ denotes a structural nonzero of the matrix,
and empty spaces denote zeros:
\begin{equation}\label{eq:fact_2rs}
\begin{pmatrix}
R_0 \\
R_1\\
\end{pmatrix}
=
\begin{pmatrix}
x & x & x & x & x \\
  & x & x & x & x \\
  &   & x & x & x \\
  &   &   & x & x \\
  &   &   &   & x \\
x & x & x & x & x \\
  & x & x & x & x \\
  &   & x & x & x \\
  &   &   & x & x \\
  &   &   &   & x \\
\end{pmatrix}.
\end{equation}
Note that we do not need to store the ones on the diagonal
explicitly.  The $Q$ factor effectively overwrites $R_1$ and the $R$
factor overwrites $R_0$.

The approach used for performing the QR factorization of the first
block column affects the storage for the Householder vectors as well
as the update of any trailing matrices that may exist.  In general,
Householder transformations have the form $I - \tau_j v_j v_j^T$, in
which the Householder vector $v_i$ is normalized so that $v_i(1) = 1$.
This means that $v_i(1)$ need not be stored explicitly.  Furthermore,
if we use \emph{structured Householder transformations}, we can avoid
storing and computing with the zeros in Equation \eqref{eq:fact_2rs}.
As the Householder vector always has the same nonzero pattern as the
vector from which it is calculated, the nonzero structure of the
Householder vector is trivial to determine.

For a $2n \times n$ rectangular matrix composed of $n \times n$ upper
triangular matrices, the $i$-th Householder vector $v_i$ in the QR
factorization of the matrix is a vector of length $2n$ with nonzeros
in entries $n+1$ through $n+i$, a one in entry $i$, and zeros
elsewhere.  If we stack all $n$ Householder vectors into a $2n \times
n$ matrix, we obtain the following representation of the $Q$ factor
(not including the $\tau$ array of multipliers):
\begin{equation}\label{eq:house:2nxn}
\begin{pmatrix}
1 &   &   &   &   \\
  & 1 &   &   &   \\
  &   & 1 &   &   \\
  &   &   & 1 &   \\
  &   &   &   & 1 \\
x & x & x & x & x  \\
  & x & x & x & x \\
  &   & x & x & x \\
  &   &   & x & x \\
  &   &   &   & x \\
\end{pmatrix}.
\end{equation}

\begin{algorithm}[h]
  \caption{Sequential QR factorization of $qn \times n$ matrix
    $A$, with structure as in Equation \eqref{eq:fact_2rs}}
\label{Alg:QR:qnxn}
\begin{algorithmic}[1]
\For{$j = 1$ to $n$}
  \State{Let $\mathcal{I}_j$ be the index set $\{j$, $n+1 : n+j$,
    $\dots$, $(q-1)n + 1 : (q-1)n + j\}$}
  \State{$w := A(\mathcal{I}_j, j)$}\Comment{Gather pivot column of $A$ into $w$}
  \State{$[\tau_j, v] := \House(w)$}\Comment{Compute Householder
    reflection, normalized so that $v(1) = 1$}
  \State{$X := A(\mathcal{I}_j, j+1:n)$}\Comment{Gather from $A$ into
    $X$.  One would normally perform the update in place; we use a
    copy to improve clarity.}
  \State{$X := (I - \tau_j v v^T) X$}\Comment{Apply Householder
    reflection}
  \State{$A(\mathcal{I}_j \setminus \{j\}, j) :=
    v(2:end)$}\Comment{Scatter $v(2:end)$ back into $A$}
  \State{$A(\mathcal{I}_j, j+1:n) := X$}\Comment{Scatter $X$ back into $A$}
\EndFor
\end{algorithmic}
\end{algorithm}

Algorithm \ref{Alg:QR:qnxn} shows a standard, column-by-column
sequential QR factorization of the $qn \times n$ matrix of upper
triangular $n \times n$ blocks, using structured Householder
reflectors.  To analyze the cost, consider the components:
\begin{enumerate}
\item $\House(w)$: the cost of this is dominated by finding the norm of
  the vector $w$ and scaling it.  
  
\item Applying a length $n$ Householder reflector, whose vector
  contains $k$ nonzeros, to an $n \times b$ matrix $A$.  This is an
  operation $(I - \tau v v^T) A = A - v (\tau (v^T A))$.  
\end{enumerate}
Appendix \ref{S:localQR-flops} counts the arithmetic operations in
detail.  There, we find that the total cost is about
\[
\frac{2}{3}(q-1)n^3 
\]
flops, to factor a $qn \times n$ matrix (we showed the specific case
$q = 2$ above).  The flop count increases by about a factor of $3
\times$ if we ignore the structure of the inputs.


\subsection{BLAS 3 structured Householder QR}\label{SS:TSQR:localQR:BLAS3structured}

\begin{algorithm}[h]
  \caption{Computing $Y$ and $T$ in the $(Y,T)$ representation
    of a collection of $n$ Householder reflectors.  Modification
    of an algorithm in \cite{schreiber1989storage} so that
    $P_j = I - \tau_j v_j v_j^T$.}
\label{Alg:YT:scaled}
\begin{algorithmic}[1]
  \Require{$n$ Householder reflectors $\rho_j = I - \tau_j v_j v_j^T$}
  \For{$j = 1$ to $n$}
    \If{$j = 1$}
      \State{$Y := [ v_1 ]$}
      \State{$T := [ -\tau_j ]$}
    \Else
      \State{$z := -\tau_j (T (Y^T v_j))$}
      \State{$Y := \begin{pmatrix} Y & v_j \end{pmatrix}$}
      \State{$T := \begin{pmatrix} T & z \\ 0 & -\tau_j \\ \end{pmatrix}$}
    \EndIf
  \EndFor
  \Ensure{$Y$ and $T$ satisfy $\rho_1 \cdot \rho_2 \cdot \dots \rho_n
    = I + Y T Y^T$}
\end{algorithmic}
\end{algorithm}

Representing the local $Q$ factor as a collection of Householder
transforms means that the local QR factorization is dominated by BLAS
2 operations (dense matrix-vector products).  A number of authors have
shown how to reformulate the standard Householder QR factorization so
as to coalesce multiple Householder reflectors into a block, so that
the factorization is dominated by BLAS 3 operations.  For example,
Schreiber and Van Loan describe a so-called YT representation of a
collection of Householder reflectors \cite{schreiber1989storage}.
BLAS 3 transformations like this are now standard in LAPACK and
ScaLAPACK.

We can adapt these techniques in a straightforward way in order to
exploit the structured Householder vectors depicted in Equation
\eqref{eq:house:2nxn}.  Schreiber and Van Loan use a slightly
different definition of Householder reflectors: $\rho_j = I - 2v_j
v_j^T$, rather than LAPACK's $\rho_j = I - \tau_j v_j v_j^T$.
Schreiber and Van Loan's $Y$ matrix is the matrix of Householder
vectors $Y = [v_1\, v_2\, \dots \, v_n]$; its construction requires no
additional computation as compared with the usual approach.  However,
the $T$ matrix must be computed, which increases the flop count by a
constant factor.  The cost of computing the $T$ factor for the $qn
\times n$ factorization above is about $qn^3 / 3$.  Algorithm
\ref{Alg:YT:scaled} shows the resulting computation.  Note that the
$T$ factor requires $n(n-1)/2$ additional storage per processor on
which the $T$ factor is required.

\subsection{Recursive Householder QR}\label{SS:TSQR:localQR:recursive}

In Section \ref{S:TSQR:perfres}, we show large performance gains
obtained by using Elmroth and Gustavson's recursive algorithm for the
local QR factorizations \cite{elmroth2000applying}.  The authors
themselves observed that their approach works especially well with
``tall thin'' matrices, and others have exploited this effect in their
applications (see e.g., \cite{rabani2001outcore}).  The recursive
approach outperforms LAPACK because it makes the panel factorization a
BLAS 3 operation.  In LAPACK, the panel QR factorization consists only
of matrix-vector and vector-vector operations.  This suggests why
recursion helps especially well with tall, thin matrices.  Elmroth and
Gustavson's basic recursive QR does not perform well when $n$ is
large, as the flop count grows cubically in $n$, so they opt for a
hybrid approach that divides the matrix into panels of columns, and
performs the panel QR factorizations using the recursive method.

Elmroth and Gustavson use exactly the same representation of the $Q$
factor as Schreiber and Van Loan \cite{schreiber1989storage}, so the
arguments of the previous section still apply.

\subsection{Trailing matrix update}\label{SS:TSQR:localQR:trailing}

Section \ref{S:CAQR} will describe how to use TSQR to factor matrices
in general 2-D layouts.  For these layouts, once the current panel
(block column) has been factored, the panels to the right of the
current panel cannot be factored until the transpose of the current
panel's $Q$ factor has been applied to them.  This is called a
\emph{trailing matrix update}.  The update lies along the critical
path of the algorithm, and consumes most of the floating-point
operations in general.  This holds regardless of whether the
factorization is left-looking, right-looking, or some hybrid of the
two.\footnote{For descriptions and illustrations of the difference
  between left-looking and right-looking factorizations, see e.g.,
  \cite{dongarra1996key}.}  Thus, it's important to make the updates
efficient.

The trailing matrix update consists of a sequence of applications of
local $Q^T$ factors to groups of ``neighboring'' trailing matrix
blocks.  (Section \ref{S:reduction} explains the meaning of the word
``neighbor'' here.)  We now explain how to do one of these local $Q^T$
applications.  (Do not confuse the local $Q$ factor, which we label
generically as $Q$, with the entire input matrix's $Q$ factor.)

Let the number of rows in a block be $M$, and the number of columns
in a block be $N$.  We assume $M \geq N$.  Suppose that we want to
apply the local $Q^T$ factor from the above $qN \times N$ matrix
factorization, to two blocks $C_0$ and $C_1$ of a trailing matrix
panel.  (This is the case $q = 2$, which we assume for simplicity.)
We divide each of the $C_i$ into a top part and a bottom part:
\[
C_i = 
\begin{pmatrix}
  C_i(1:N, :) \\
  C_i(N+1 : M, :)
\end{pmatrix} =
\begin{pmatrix}
  C_i' \\
  C_i''
\end{pmatrix}.
\]
Our goal is to perform the operation
\[
\begin{pmatrix}
 R_0 & C_0' \\
 R_1 & C_1' \\
\end{pmatrix}
= 
\begin{pmatrix}
QR & C_0' \\
   & C_1' \\
\end{pmatrix}
=
Q \cdot
\begin{pmatrix}
R & \hat{C}_0' \\
  & \hat{C}_1' \\
\end{pmatrix},
\]
in which $Q$ is the local $Q$ factor and $R$ is the local $R$ factor
of $[R_0; R_1]$.  Implicitly, the local $Q$ factor has the dimensions
$2M \times 2M$, as Section \ref{S:TSQR:algebra} explains.  However, it
is not stored explicitly, and the implicit operator that is stored has
the dimensions $2N \times 2N$.  We assume that processors $P_0$ and
$P_1$ each store a redundant copy of $Q$, that processor $P_2$ has
$C_0$, and that processor $P_3$ has $C_1$.  We want to apply $Q^T$ to
the matrix
\[
C = 
\begin{pmatrix}
C_0 \\
C_1 \\
\end{pmatrix}.
\]
First, note that $Q$ has a specific structure.  If stored explicitly,
it would have the form
\[
Q = 
\begin{pmatrix}
  \begin{matrix}
    U_{00} & \\
          & I_{M - N} 
  \end{matrix} & 
  \begin{matrix}
    U_{01} &  \\
          & \mathbf{0}_{M - N}
  \end{matrix} \\
  \begin{matrix}
    U_{10} &  \\
          & \mathbf{0}_{M - N}
  \end{matrix} &
  \begin{matrix}
    U_{11} & \\
          & I_{M - N} 
  \end{matrix} \\
\end{pmatrix},
\]
in which the $U_{ij}$ blocks are each $N \times N$.  This makes the
only nontrivial computation when applying $Q^T$ the following:
\begin{equation}\label{eq:localQTupdate}
\begin{pmatrix}
  \hat{C}_0' \\
  \hat{C}_1' \\
\end{pmatrix}
:=
\begin{pmatrix}
  U_{00}^T & U_{10}^T \\
  U_{01}^T & U_{11}^T
\end{pmatrix}
\cdot
\begin{pmatrix}
C_0' \\
C_1' \\
\end{pmatrix}.
\end{equation}
We see, in particular, that only the uppermost $N$ rows of each block
of the trailing matrix need to be read or written.  Note that it is
not necessary to construct the $U_{ij}$ factors explicitly; we need
only operate on $C_0'$ and $C_1'$ with $Q^T$.

If we are using a standard Householder QR factorization (without BLAS
3 optimizations), then computing Equation \eqref{eq:localQTupdate} is
straightforward.  When one wishes to exploit structure (as in Section
\ref{SS:TSQR:localQR:structured}) and use a local QR factorization
that exploits BLAS 3 operations (as in Section
\ref{SS:TSQR:localQR:BLAS3structured}), more interesting load balance
issues arise.  We will discuss these in the following section.

\subsubsection{Trailing matrix update with structured BLAS 3 QR}

An interesting attribute of the YT representation is that the $T$
factor can be constructed using only the $Y$ factor and the $\tau$
multipliers.  This means that it is unnecessary to send the $T$ factor
for updating the trailing matrix; the receiving processors can each
compute it themselves.  However, one cannot compute $Y$ from $T$
and $\tau$ in general.  

When the YT representation is used, the update of the trailing
matrices takes the following form:
\[
\begin{pmatrix}
\hat{C_0}' \\
\hat{C_1}' \\
\end{pmatrix}
:=
\begin{pmatrix}
I -
\begin{pmatrix}
 I \\
Y_1 \\
\end{pmatrix}
\cdot 
\begin{array}{cc} 	
T^T \\
\\
\end{array}
\cdot
\begin{pmatrix}
 I \\
Y_1 \\
\end{pmatrix}^T
\end{pmatrix}
\begin{pmatrix}
C_0' \\
C_1' \\
\end{pmatrix}.
\]
Here, $Y_1$ starts on processor $P_1$, $C_0'$ on processor $P_2$, and
$C_1'$ on processor $P_3$.  The matrix $T$ must be computed from
$\tau$ and $Y_1$; we can assume that $\tau$ is on processor $P_1$.
The updated matrices $\hat{C_0}'$ and $\hat{C_1}'$ are on processors
$P_2$ resp.\ $P_3$.  

There are many different ways to perform this parallel update.  The
data dependencies impose a directed acyclic graph (DAG) on the flow of
data between processors.  One can find the the best way to do the
update by realizing an optimal computation schedule on the DAG.  Our
performance models can be used to estimate the cost of a particular
schedule.

Here is a straightforward but possibly suboptimal schedule.  First,
assume that $Y_1$ and $\tau$ have already been sent to $P_3$.  Then,
\begin{multicols}{2}
$P_2$'s tasks:
\begin{itemize}
  \item Send $C_0'$ to $P_3$
  \item Receive $W$ from $P_3$
  \item Compute $\hat{C_0}' = C_0' - W$
\end{itemize}
\vspace{1cm}

$P_3$'s tasks:
\begin{itemize}
  \item Compute the $T$ factor and $W := T^T (C_0' + Y_1^T C_1')$
  \item Send $W$ to $P_2$
  \item Compute $\hat{C_1}' := C_1' - Y_1 W$
\end{itemize}
\end{multicols}
However, this leads to some load imbalance, since $P_3$ performs more
computation than $P_2$.  It does not help to compute $T$ on $P_0$ or
$P_1$ before sending it to $P_3$, because the computation of $T$ lies
on the critical path in any case.  We will see in Section \ref{S:CAQR}
that part of this computation can be overlapped with the
communication.  

For $q \geq 2$, we can write the update operation as
\[
\begin{pmatrix}
\hat{C_0}' \\
\hat{C_1}' \\
\vdots \\
\hat{C_{q-1}}' \\
\end{pmatrix}
:=
\left( I - 
\begin{pmatrix}
I_{N \times N} \\
Y_1 \\
\vdots \\
Y_{q-1} \\
\end{pmatrix}
T^T
\begin{pmatrix}
I_{N \times N} & Y_1^T & \dots & Y_{q-1}^T \\
\end{pmatrix}
\right)
\begin{pmatrix}
C_0' \\
C_1' \\
\vdots \\
C_{q-1}' \\
\end{pmatrix}.
\]
If we let
\[
D := C_0' + Y_1^T C_1' + Y_2^T C_2' + \dots + Y_{q-1}^T C_{q-1}'
\]
be the ``inner product'' part of the update operation formulas, then
we can rewrite the update formulas as
\[
\begin{aligned}
\hat{C_0}' &:= C_0' - T^T D, \\
\hat{C_1}' &:= C_1' - Y_1 T^T D, \\
\vdots  & \\
\hat{C_{q-1}}' &:= C_{q-1}' - Y_{q-1} T^T D.
\end{aligned}
\]
As the branching factor $q$ gets larger, the load imbalance becomes
less of an issue.  The inner product $D$ should be computed as an
all-reduce in which the processor owning $C_i$ receives $Y_i$ and
$T$.  Thus, all the processors but one will have the same
computational load.

\section{Machine model}\label{S:perfmodel}

\subsection{Parallel machine model}

Throughout this work, we use the ``alpha-beta'' or latency-bandwidth
model of communication, in which a message of size $n$ floating-point
words takes time $\alpha + \beta n$ seconds.  The $\alpha$ term
represents message latency (seconds per message), and the $\beta$ term
inverse bandwidth (seconds per floating-point word communicated).  Our
algorithms only need to communicate floating-point words, all of the
same size.  We make no attempt to model overlap of communication and
computation, but we do mention the possibility of overlap when it
exists.  Exploiting overlap could potentially speed up our algorithms
(or any algorithm) by a factor of two.

We predict floating-point performance by counting floating-point
operations and multiplying them by $\gamma$, the inverse peak
floating-point performance, also known as the floating-point
throughput.  The quantity $\gamma$ has units of seconds per flop (so
it can be said to measure the bandwidth of the floating-point
hardware).  If we need to distinguish between adds and multiplies on
one hand, and divides on the other, we use $\gamma$ for the throughput
of adds and multiplies, and $\gamma_d$ for the throughput of divides.

When appropriate, we may scale the peak floating-point performance
prediction of a particular matrix operation by a factor, in order to
account for the measured best floating-point performance of local QR
factorizations.  This generally gives the advantage to competing
algorithms rather than our own, as our algorithms are designed to
perform better when communication is much slower than arithmetic.

\subsection{Sequential machine model}

We also apply the alpha-beta model to communication between levels of
the memory hierarchy in the sequential case.  We restrict our model to
describe only two levels at one time: fast memory (which is smaller)
and slow memory (which is larger).  The terms ``fast'' and ``slow''
are always relative.  For example, DRAM may be considered fast if the
slow memory is disk, but DRAM may be considered slow if the fast
memory is cache.  As in the parallel case, the time to complete a
transfer between two levels is modeled as $\alpha + \beta n$.  We
assume that user has explicit control over data movement (reads and
writes) between fast and slow memory.  This offers an upper bound when
control is implicit (as with caches), and also allows our model as
well as our algorithms to extend to systems like the Cell processor
(in which case fast memory is an individual local store, and slow
memory is DRAM).

We assume that the fast memory can hold $W$ floating-point words.  For
any QR factorization operating on an $m \times n$ matrix, the quantity
\[
\frac{mn}{W}
\]
bounds from below the number of loads from slow memory into fast
memory (as the method must read each entry of the matrix at least
once).  It is also a lower bound on the number of stores from fast
memory to slow memory (as we assume that the algorithm must write the
computed $Q$ and $R$ factors back to slow memory).  Sometimes we may
refer to the block size $P$.  In the case of TSQR, we usually choose
\[
P = \frac{mn}{3W},
\]
since at most three blocks of size $P$ must be in fast memory at one
time when applying the $Q$ or $Q^T$ factor in sequential TSQR (see
Section \ref{S:TSQR:algebra}).

In the sequential case, just as in the parallel case, we assume all
memory transfers are nonoverlapped.  Overlapping communication and
computation may provide up to a twofold performance improvement.
However, some implementations may consume fast memory space in order
to do buffering correctly.  This matters because the main goal of our
sequential algorithms is to control fast memory usage, often to solve
problems that do not fit in fast memory.  We usually want to use as
much of fast memory as possible, in order to avoid expensive transfers
to and from slow memory.


\section{TSQR implementation}\label{S:TSQR:impl}

In this section, we describe the TSQR factorization algorithm in
detail.  We also build a performance model of the algorithm, based on
the machine model in Section \ref{S:perfmodel} and the operation
counts of the local QR factorizations in Section \ref{S:TSQR:localQR}.
Parallel TSQR performs $2mn^2/P + \frac{2n^3}{3}\log P$ flops,
compared to the $2mn^2/P - 2n^3/(3P)$ flops performed by ScaLAPACK's
parallel QR factorization \lstinline!PDGEQRF!, but requires $2n$ times
fewer messages.  The sequential TSQR factorization performs the same
number of flops as sequential blocked Householder QR, but requires
$\Theta(n)$ times fewer transfers between slow and fast memory, and a
factor of $\Theta(m/\sqrt{W})$ times fewer words transferred, in which
$W$ is the fast memory size.

\subsection{Reductions and all-reductions}

In Section \ref{S:reduction}, we gave a detailed description of
(all-)reductions.  We did so because the TSQR factorization is itself an
(all-)reduction, in which additional data (the components of the $Q$
factor) is stored at each node of the (all-)reduction tree.  Applying
the $Q$ or $Q^T$ factor is also a(n) (all-)reduction.

If we implement TSQR with an all-reduction, then we get the final $R$
factor replicated over all the processors.  This is especially useful
for Krylov subspace methods.  If we implement TSQR with a reduction,
then the final $R$ factor is stored only on one processor.  This
avoids redundant computation, and is useful both for block column
factorizations for 2-D block (cyclic) matrix layouts, and for solving
least squares problems when the $Q$ factor is not needed.

\subsection{Factorization}

We now describe the parallel and sequential TSQR factorizations for
the 1-D block row layout.  (We omit the obvious generalization to a
1-D block cyclic row layout.)

Parallel TSQR computes an $R$ factor which is duplicated over all the
processors, and a $Q$ factor which is stored implicitly in a
distributed way.  The algorithm overwrites the lower trapezoid of
$A_{i}$ with the set of Householder reflectors for that block, and the
$\tau$ array of scaling factors for these reflectors is stored
separately.  The matrix $R_{i,k}$ is stored as an $n \times n$ upper
triangular matrix for all stages $k$.  Algorithm \ref{Alg:TSQR:par}
shows an implementation of parallel TSQR, based on an all-reduction.
(Note that running Algorithm \ref{Alg:TSQR:par} on a matrix stored in
a 1-D block cyclic layout still works, though it performs an implicit
block row permutation on the $Q$ factor.)

\begin{algorithm}[h]
\caption{Parallel TSQR}
\label{Alg:TSQR:allred:blkrow}
\label{Alg:TSQR:par}
\begin{algorithmic}[1]
\Require{$\Pi$ is the set of $P$ processors}
\Require{All-reduction tree with height $L$.  If $P$ is a power of two
  and we want a binary all-reduction tree, then $L = \log_2 P$.}
\Require{$i \in \Pi$: my processor's index}
\Require{The $m \times n$ input matrix $A$ is distributed in a 1-D 
  block row layout over the processors; $A_{i}$ is the block of rows 
  belonging to processor $i$}.
\State{Compute $[Q_{i,0}, R_{i,0}] := qr(A_{i})$ using sequential 
  Householder QR}\label{Alg:TSQR:allred:blkrow:QR1}
\For{$k$ from 1 to $L$}
  \If{I have any neighbors in the all-reduction tree at this level}
    \State{Send (non-blocking) $R_{i,k-1}$ to each neighbor not myself}
    \State{Receive (non-blocking) $R_{j,k-1}$ from each neighbor $j$ not myself}
    \State{Wait until the above sends and receives
      complete}\Comment{Note: \emph{not} a global barrier.}
    \State{Stack the upper triangular $R_{j,k-1}$ from all neighbors 
      (including my own $R_{i,k-1}$), by order of processor ids, into 
      a $qn \times n$ array $C$, in which $q$ is the number of 
      neighbors.}
    \State{Compute $[Q_{i,k}, R_{i,k}] := qr(C)$ using Algorithm
      \ref{Alg:QR:qnxn} in Section \ref{SS:TSQR:localQR:structured}}
  \Else
    \State{$R_{i,k} := R_{i,k-1}$}
    \State{$Q_{i,k} := I_{n \times n}$}\Comment{Stored implicitly}
  \EndIf
  \State{Processor $i$ has an implicit representation of its block column
    of $Q_{i,k}$.  The blocks in the block column are $n \times n$
    each and there are as many of them as there are neighbors at stage
    $k$ (including $i$ itself).  We don't need to compute the blocks
    explicitly here.}
\EndFor
\Ensure{$R_{i,L}$ is the $R$ factor of $A$, for all processors $i \in \Pi$.}
\Ensure{The $Q$ factor is implicitly represented by $\{Q_{i,k}\}$:
  $i \in \Pi$, $k \in \{0, 1, \dots, L\}\}$.}
\end{algorithmic}
\end{algorithm}

Sequential TSQR begins with an $m \times n$ matrix $A$ stored in slow
memory.  The matrix $A$ is divided into $P$ blocks $A_0$, $A_1$,
$\dots$, $A_{P-1}$, each of size $m/P \times n$.  (Here, $P$ has
nothing to do with the number of processors.)  Each block of $A$ is
loaded into fast memory in turn, combined with the $R$ factor from the
previous step using a QR factorization, and the resulting $Q$ factor
written back to slow memory.  Thus, only one $m/P \times n$ block of
$A$ resides in fast memory at one time, along with an $n \times n$
upper triangular $R$ factor.  Sequential TSQR computes an $n \times n$
$R$ factor which ends up in fast memory, and a $Q$ factor which is
stored implicitly in slow memory as a set of blocks of Householder
reflectors.  Algorithm \ref{Alg:TSQR:seq} shows an implementation of
sequential TSQR.

\begin{algorithm}[h]
\caption{Sequential TSQR}\label{Alg:TSQR:seq}
\begin{algorithmic}[1]
\Require{The $m \times n$ input matrix $A$, stored in slow memory, 
  is divided into $P$ row blocks $A_0$, $A_1$, $\dots$, $A_{P-1}$}
\State{Load $A_0$ into fast memory}
\State{Compute $[Q_{00}, R_{00}] := qr(A_{0})$ using standard
  sequential QR.   Here, the $Q$ factor is represented implicitly by
  an $m/P \times n$ lower triangular array of Householder reflectors 
  $Y_{00}$ and their $n$ associated scaling factors $\tau_{00}$}
\State{Write $Y_{00}$ and $\tau_{00}$ back to slow memory; keep
  $R_{00}$ in fast memory}
\For{$k = 1$ to $P - 1$}
    \State{Load $A_k$}
    \State{Compute $[Q_{01}, R_{01}] = qr([R_{0,k-1}; A_k])$ 
      using the structured method analyzed in Appendix 
      \ref{SSS:localQR-flops:seq:2blocks}.   Here, the $Q$ factor 
      is represented implicitly by a full $m/P \times n$ array of 
      Householder reflectors $Y_{0k}$ and their $n$ associated 
      scaling factors $\tau_{0k}$.}
    \State{Write $Y_{0k}$ and $\tau_{0k}$ back to slow memory; 
      keep $R_{0k}$ in fast memory}
\EndFor
\Ensure{$R_{0,P-1}$ is the $R$ factor in the QR factorization of $A$,
  and is in fast memory}
\Ensure{The $Q$ factor is implicitly represented by $Q_{00}$,
  $Q_{01}$, $\dots$, $Q_{0,P-1}$, and is in slow memory}
\end{algorithmic}
\end{algorithm}

\subsubsection{Performance model}

In Appendix \ref{S:TSQR-par-detailed}, we develop a performance model
for parallel TSQR on a binary tree.  Appendix \ref{S:TSQR-seq-detailed}
does the same for sequential TSQR on a flat tree.

A parallel TSQR factorization on a binary reduction tree performs the
following computations along the critical path: One local QR
factorization of a fully dense $m/P \times n$ matrix, and $\log P$
factorizations, each of a $2n \times n$ matrix consisting of two $n
\times n$ upper triangular matrices.  The factorization requires
\[
\frac{2mn^2}{P} + \frac{2n^3}{3} \log P
\]
flops and $\log P$ messages, and transfers a total of $(1/2) n^2 \log
P$ words between processors.  In contrast, parallel Householder QR
requires
\[
\frac{2mn^2}{P} - \frac{2n^3}{3}
\]
flops and $2n \log P$ messages, but also transfers $(1/2) n^2 \log P$
words between processors.  For details, see Table
\ref{tbl:QR:perfcomp:par} in Section \ref{S:TSQR:perfcomp}.

Sequential TSQR on a flat tree performs the same number of flops as
sequential Householder QR, namely 
\[
2mn^2 - \frac{2n^3}{3}
\]
flops.  However, sequential TSQR only transfers 
\[
2mn - \frac{n(n+1)}{2} + \frac{mn^2}{\tilde{W}}
\]
words between slow and fast memory, in which $\tilde{W} = W -
n(n+1)/2$, and only performs
\[
\frac{2mn}{\tilde{W}}
\]
transfers between slow and fast memory.  In contrast, blocked
sequential Householder QR transfers 
\[
\frac{m^2 n^2}{2W} 
- \frac{mn^3}{6W}
+ \frac{3mn}{2} 
- \frac{3n^2}{4}
\]
words between slow and fast memory, and only performs
\[
\frac{2mn}{W} + \frac{mn^2}{2W}
\]
transfers between slow and fast memory.  For details, see Table
\ref{tbl:QR:perfcomp:seq} in Section \ref{S:TSQR:perfcomp}.

\subsection{Applying  $Q$ or $Q^T$ to vector(s)}\label{SS:TSQR:application}

Just like Householder QR, TSQR computes an implicit representation of
the $Q$ factor.  One need not generate an explicit representation of
$Q$ in order to apply the $Q$ or $Q^T$ operators to one or more
vectors.  In fact, generating an explicit $Q$ matrix requires just as
many messages as applying $Q$ or $Q^T$.  (The performance model for
applying $Q$ or $Q^T$ is an obvious extension of the factorization
performance model; the parallel performance model is developed in
Appendix \ref{SS:TSQR-par-detailed:apply} and the sequential
performance model in Appendix \ref{SS:TSQR-seq-detailed:apply}.)
Furthermore, the implicit representation can be updated or downdated,
by using standard techniques (see e.g., \cite{govl:96}) on the local
QR factorizations recursively.  The $s$-step Krylov methods mentioned
in Section \ref{S:motivation} employ updating and downdating
extensively.

In the case of the ``thin'' $Q$ factor (in which the vector input is
of length $n$), applying $Q$ involves a kind of broadcast operation
(which is the opposite of a reduction).  If the ``full'' $Q$ factor is
desired, then applying $Q$ or $Q^T$ is a kind of all-to-all (like the
fast Fourier transform).  Computing $Q \cdot x$ runs through the nodes
of the (all-)reduction tree from leaves to root, whereas computing
$Q^T \cdot y$ runs from root to leaves.



\section{Other ``tall skinny'' QR algorithms}\label{S:TSQR:perfcomp}

There are many other algorithms besides TSQR for computing the QR
factorization of a tall skinny matrix.  They differ in terms of
performance and accuracy, and may store the $Q$ factor in different
ways that favor certain applications over others.  In this section, we
model the performance of the following competitors to TSQR:
\begin{itemize}
\item Four different Gram-Schmidt variants
\item CholeskyQR (see \cite{stwu:02})
\item Householder QR, with a block row layout
\end{itemize}
Each includes parallel and sequential versions.  For Householder QR,
we base our parallel model on the ScaLAPACK routine
\lstinline!PDGEQRF!, and the sequential model on left-looking blocked
Householder.  Our left-looking blocked Householder implementation is
modeled on the out-of-core ScaLAPACK routine \lstinline!PFDGEQRF!,
which is left-looking instead of right-looking in order to minimize
the number of writes to slow memory (the total amount of data moved
between slow and fast memory is the same for both left-looking and
right-looking blocked Householder QR).  See Appendix \ref{S:PFDGEQRF}
for details.  In the subsequent Section \ref{S:TSQR:stability}, we
summarize the numerical accuracy of these QR factorization methods,
and discuss their suitability for different applications.

In the parallel case, CholeskyQR and TSQR have comparable numbers of
messages and communicate comparable numbers of words, but CholeskyQR
requires a constant factor fewer flops along the critical path.
However, the $Q$ factor computed by TSQR is always numerically
orthogonal, whereas the $Q$ factor computed by CholeskyQR loses
orthogonality proportionally to $\kappa_2(A)^2$.  The variants of
Gram-Schmidt require at best a factor $n$ more messages than these two
algorithms, and lose orthogonality at best proportionally to
$\kappa_2(A)$.

\subsection{Gram-Schmidt orthogonalization}\label{SS:TSQR:perfcomp:GS}

Gram-Schmidt has two commonly used variations: ``classical'' (CGS) and
``modified'' (MGS).  Both versions have the same floating-point
operation count, but MGS performs them in a different order to improve
stability.  We will show that a parallel implementation of MGS uses at
best $2n \log P$ messages, in which $P$ is the number of processors,
and a blocked sequential implementation requires at least
\[
\frac{m n^2}{ 2W - n(n+1) }
\]
transfers between slow and fast memory, in which $W$ is the fast
memory capacity.  In contrast, parallel TSQR requires only $\log P$
messages, and sequential TSQR only requires
\[
\frac{4mn}{2W - n(n+1)}
\]
transfers between slow and fast memory, a factor of about $n/4$ less.
See Tables \ref{tbl:gram-schmidt-variants:par} and
\ref{tbl:gram-schmidt-variants:seq} for details.

\subsubsection{Left- and right-looking}

Just like many matrix factorizations, both MGS and CGS come in
left-looking and right-looking variants.  To distinguish between the
variants, we append ``\_L'' resp. ``\_R'' to the algorithm name to
denote left- resp.\ right-looking.  We show all four combinations as
Algorithms \ref{Alg:MGS:RL}--\ref{Alg:CGS:LL}.  Both right-looking and
left-looking variants loop from left to right over the columns of the
matrix $A$.  At iteration $k$ of this loop, the left-looking version
only accesses columns $1$ to $k$ inclusive, whereas the right-looking
version only accesses columns $k$ to $n$ inclusive.  Thus,
right-looking algorithms require the entire matrix to be available,
which forbids their use when the matrix is to be generated and
orthogonalized one column at a time.  (In this case, only left-looking
algorithms may be used.)  We assume here that the entire matrix is
available at the start of the algorithm.

Right-looking Gram-Schmidt is usually called ``row-oriented
Gram-Schmidt,'' and by analogy, left-looking Gram-Schmidt is usually
called ``column-oriented Gram-Schmidt.''  We use the terms
``right-looking'' resp.\ ``left-looking'' for consistency with the
other QR factorization algorithms in this paper.

\begin{algorithm}[h]
\caption{Modified Gram-Schmidt, right-looking}\label{Alg:MGS:RL}
\begin{algorithmic}[1]
\Require{$A$: $m \times n$ matrix with $m \geq n$}
\For{$k = 1$ to $n$}
  \State{$R(k,k) := \|A(:,k)\|_2$}  
  \State{$Q(:,k) := A(:,k) / R(k,k)$} 
  \State{$R(k, k+1 : n) := Q(:,k)^T \cdot A(:, k+1 : n)$} 
  \State{$A(:, k+1 : n) := A(:, k+1:n) - R(k, k+1:n) \cdot Q(:,k)$} 
\EndFor
\end{algorithmic}
\end{algorithm}

\begin{algorithm}[h]
\caption{Modified Gram-Schmidt, left-looking}\label{Alg:MGS:LL}
\begin{algorithmic}[1]
\Require{$A$: $m \times n$ matrix with $m \geq n$}
\For{$k = 1$ to $n$}
  \State{$v := A(:, k)$} 
  \For{$j = 1$ to $k - 1$}\Comment{Data dependencies hinder vectorization}
    \State{$R(j,k) := Q(:,j)^T \cdot v$}\Comment{Change $v$ to $A(:,k)$ to get CGS}
    \State{$v := v - R(j,k) \cdot Q(:, j)$}  
  \EndFor
  \State{$R(k,k) := \|v\|_2$}
  \State{$Q(:, k) := v / R(k,k)$}
\EndFor  
\end{algorithmic}
\end{algorithm}

\begin{algorithm}[h]
\caption{Classical Gram-Schmidt, right-looking}\label{Alg:CGS:RL}
\begin{algorithmic}[1]
\Require{$A$: $m \times n$ matrix with $m \geq n$}
\State{$V := A$}\Comment{Not copying $A$ would give us right-looking MGS.}
\For{$k = 1$ to $n$}
  \State{$R(k,k) := \|V(:,k)\|_2$} 
  \State{$Q(:,k) := V(:,k) / R(k,k)$}
  \State{$R(k, k+1 : n) := Q(:,k)^T \cdot A(:, k+1 : n)$}
  \State{$V(:, k+1 : n) := V(:, k+1:n) - R(k, k+1:n) \cdot Q(:,k)$}
\EndFor
\end{algorithmic}
\end{algorithm}

\begin{algorithm}[h]
\caption{Classical Gram-Schmidt, left-looking}\label{Alg:CGS:LL}
\begin{algorithmic}[1]
\Require{$A$: $m \times n$ matrix with $m \geq n$}
\For{$k = 1$ to $n$}
  \State{$R(1:k-1, k) := Q(:, 1:k-1)^T \cdot A(:,k)$}\Comment{This and
    the next statement are not vectorized in left-looking MGS.}
  \State{$A(:,k) := A(:,k) - R(1:k-1, k) \cdot Q(:,
    1:k-1)$}\Comment{In the sequential case, one can coalesce the read
    of each block of $A(:,k)$ in this statement with the read of each
    block of $A(:,k)$ in the next statement.}
  \State{$R(k,k) := \|A(:,k)\|_2$}
  \State{$Q(:, k) := A(:,k) / R(k,k)$}
\EndFor  
\end{algorithmic}
\end{algorithm}

\subsubsection{Reorthogonalization}

One can improve the stability of CGS by reorthogonalizing the vectors.
The simplest way is to make two orthogonalization passes per column,
that is, to orthogonalize the current column against all the previous
columns twice.  We call this ``CGS2.''  This method only makes sense
for left-looking Gram-Schmidt, when there is a clear definition of
``previous columns.''  Normally one would orthogonalize the column
against all previous columns once, and then use some orthogonality
criterion to decide whether to reorthogonalize the column.  As a
result, the performance of CGS2 is data-dependent, so we do not model
its performance here.  In the worst case, it can cost twice as much as
CGS\_L.  Section \ref{S:TSQR:stability} discusses the numerical
stability of CGS2 and why ``twice is enough.''

\subsubsection{Parallel Gram-Schmidt}

MGS\_L (Algorithm \ref{Alg:MGS:LL}) requires about $n/4$ times more
messages than MGS\_R (Algorithm \ref{Alg:MGS:RL}), since the
left-looking algorithm's data dependencies prevent the use of
matrix-vector products.  CGS\_R (Algorithm \ref{Alg:CGS:RL}) requires
copying the entire input matrix; not doing so results in MGS\_R
(Algorithm \ref{Alg:MGS:RL}), which is more numerically stable in any
case.  Thus, for the parallel case, we favor MGS\_R and CGS\_L for a
fair comparison with TSQR.  

In the parallel case, all four variants of MGS and CGS listed here
require
\[
\frac{2mn^2}{P} + O\left( \frac{mn}{P} \right)
\]
arithmetic operations, and involve communicating 
\[
\frac{n^2}{2}\log(P) + O(n\log(P))
\]
floating-point words in total.  MGS\_L requires
\[
\frac{n^2}{2} \log(P) + O(n \log(P)) 
\]
messages, whereas the other versions only need $2 n \log(P)$ messages.
Table \ref{tbl:gram-schmidt-variants:par} shows all four performance
models.  

\begin{table}
  \centering
  \begin{tabular}{l|c|c|c}
    Parallel algorithm & \# flops & \# messages & \# words \\ \hline
    Right-looking MGS & $2mn^2/P$ 
                      & $2n \log(P)$ 
                      & $\frac{n^2}{2}\log(P)$ \\
    Left-looking MGS  & $2mn^2/P$ 
                      & $\frac{n^2}{2}\log(P)$ 
                      & $\frac{n^2}{2}\log(P)$ \\
    Right-looking CGS & $2mn^2/P$  
                      & $2n \log(P)$ 
                      & $\frac{n^2}{2}\log(P)$ \\
    Left-looking CGS  & $2mn^2/P$ 
                      & $2n \log(P)$ 
                      & $\frac{n^2}{2}\log(P)$ \\
  \end{tabular}
  \caption{Arithmetic operation counts, number of messages,
    and total communication volume (in number of words 
    transferred) for parallel left-looking and right-looking variants 
    of CGS and MGS.  Reductions are performed using a binary tree on
    $P$ processors.  Lower-order terms omitted.}
  \label{tbl:gram-schmidt-variants:par}
\end{table}

\subsubsection{Sequential Gram-Schmidt}

For one-sided factorizations in the out-of-slow-memory regime,
left-looking algorithms require fewer writes than their right-looking
analogues (see e.g., \cite{toledo99survey}).  We will see this in the
results below, which is why we spend more effort analyzing the
left-looking variants.

Both MGS and CGS can be reorganized into blocked variants that work on
panels.  These variants perform the same floating-point operations as
their unblocked counterparts, but save some communication.  In the
parallel case, the blocked algorithms encounter the same latency
bottleneck as ScaLAPACK's parallel QR factorization
\lstinline!PDGEQRF!, so we do not analyze them here.  The sequential
case offers more potential for communication savings.  

The analysis of blocked sequential Gram-Schmidt's communication costs
resembles that of blocked left-looking Householder QR (see Appendix
\ref{S:PFDGEQRF}), except that Gram-Schmidt computes and stores the
$Q$ factor explicitly.  This means that Gram-Schmidt stores the upper
triangle of the matrix twice: once for the $R$ factor, and once for
the orthogonalized vectors.  Left-looking MGS and CGS would use a left
panel of width $b$ and a current panel of width $c$, just like blocked
left-looking Householder QR.  Right-looking Gram-Schmidt would use a
current panel of width $c$ and a right panel of width $b$.  Unlike
Householder QR, however, Gram-Schmidt requires storing the $R$ factor
separately, rather than overwriting the original matrix's upper
triangle.  We assume here that $W \gg n(n+1)/2$, so that the entire
$R$ factor can be stored in fast memory.  This need not be the case,
but it is a reasonable assumption for the ``tall skinny'' regime.  If
$m \approx n$, then Gram-Schmidt's additional bandwidth requirements,
due to working with the upper triangle twice (once for the $Q$ factor
and once for the $R$ factor), make Householder QR more competitive
than Gram-Schmidt.

For sequential MGS\_L, the number of words transferred between slow
and fast memory is about
\begin{equation}\label{eq:MGSL:seq:bw:bc}
\sum_{j=1}{\frac{n}{c}} \left( 
    2cm + \sum{k=1}{\frac{c(j-1)}{b}} bm
\right) = 
\frac{3mn}{2} + \frac{m n^2}{2c},
\end{equation}
and the number of messages is about
\begin{equation}\label{eq:MGSL:seq:lat:bc}
\frac{n^2}{2 b c} 
+ \frac{2n}{c}
- \frac{n}{2b}.
\end{equation}
The fast memory usage is about $(b+c)m + n(n+1)/2$, so if we optimize
for bandwidth and take $b = 1$ and 
\[
c \approx \frac{W - n(n+1)/2}{m},
\]
the number of words transferred between slow and fast memory is about
\begin{equation}\label{eq:MGSL:seq:bw:W}
\frac{3mn}{2} 
+ \frac{m^2 n^2}{2W - n(n+1)},
\end{equation}
and the number of messages is about (using the highest-order term only)
\begin{equation}\label{eq:MGSL:seq:lat:W}
\frac{m n^2}{ 2W - n(n+1) }.
\end{equation}

For sequential MGS\_R, the number of words transferred between slow
and fast memory is about
\begin{equation}\label{eq:MGSR:seq:bw:bc}
\sum_{j=1}^{\frac{n}{c}} \left(
    2cm + \sum_{k=\frac{c(j-1)}{b}}^{\frac{n}{b}} 2bm
\right) =
3mn + \frac{m n^2}{c} + \frac{2bmn}{c}.
\end{equation}
This is always greater than the number of words transferred by
MGS\_L.  The number of messages is about 
\begin{equation}\label{eq:MGSR:seq:lat:bc}
\frac{n^2}{bc} + \frac{4n}{c} + \frac{n}{b},
\end{equation}
which is also always greater than the number of messages transferred
by MGS\_L.  Further analysis is therefore unnecessary; we should
always use the left-looking version.

Table \ref{tbl:gram-schmidt-variants:seq} shows performance models for
blocked versions of left-looking and right-looking sequential MGS and
CGS.  We omit CGS\_R as it requires extra storage and provides no
benefits over MGS\_R.

\begin{table}
  \centering
  \begin{tabular}{l|c|c|c}
    Sequential algorithm & \# flops & \# messages & \# words \\ \hline
    Right-looking MGS
      & $2mn^2$ 
      & $\frac{2 m n^2}{2W - n(n+1)}$
      & $3mn + \frac{2 m n^2}{2W - n(n+1)}$ \\
    Left-looking MGS
      & $2mn^2$
      & $\frac{m n^2}{ 2W - n(n+1) }$
      & $\frac{3mn}{2} + \frac{m^2 n^2}{2W - n(n+1)}$ \\
    Left-looking CGS
      & $2mn^2$ 
      & $\frac{m n^2}{ 2W - n(n+1) }$
      & $\frac{3mn}{2} + \frac{m^2 n^2}{2W - n(n+1)}$ \\
  \end{tabular}
  \caption{Arithmetic operation counts, number of reads and writes, 
    and total communication volume (in number of words read and 
    written) for sequential left-looking CGS and MGS.  $W$ is the fast
    memory capacity in number of floating-point words.  Lower-order
    terms omitted.}
  \label{tbl:gram-schmidt-variants:seq}
\end{table}


\subsection{CholeskyQR}\label{SS:TSQR:perfcomp:CholeskyQR}

\begin{algorithm}[h]
\caption{CholeskyQR factorization}\label{Alg:CholeskyQR}
\begin{algorithmic}[1]
  \Require{$A$:  $m \times n$ matrix with $m \geq n$}
  \State{$W := A^T A$}\Comment{(All-)reduction}
  \State{Compute the Cholesky factorization $L \cdot L^T$ of $W$}
  \State{$Q := A L^{-T}$}
  \Ensure{$[Q, L^T]$ is the QR factorization of $A$}
\end{algorithmic}
\end{algorithm}

CholeskyQR (Algorithm \ref{Alg:CholeskyQR}) is a QR factorization that
requires only one all-reduction \cite{stwu:02}.  In the parallel case, it
requires $\log_2 P$ messages, where $P$ is the number of processors.
In the sequential case, it reads the input matrix only once.  Thus, it
is optimal in the same sense that TSQR is optimal.  Furthermore, the
reduction operator is matrix-matrix addition rather than a QR
factorization of a matrix with comparable dimensions, so CholeskyQR
should always be faster than TSQR.  Section \ref{S:TSQR:perfres}
supports this claim with performance data on a cluster.  Note that in
the sequential case, $P$ is the number of blocks, and we assume
conservatively that fast memory must hold $2mn/P$ words at once (so
that $W = 2mn/P$).

\begin{table}
\centering
\begin{tabular}{l|c|c|c}
Algorithm  & \# flops & \# messages & \# words \\ \hline
Parallel CholeskyQR   & $\frac{2mn^2}{P} + \frac{n^3}{3}$ 
           & $\log(P)$ 
           & $\frac{n^2}{2}\log(P)$ \\ 
Sequential CholeskyQR & $2mn^2 + \frac{n^3}{3}$  
           & $\frac{6mn}{W}$
           & $3mn$ \\ 
\end{tabular}
\caption{Performance model of the parallel and sequential 
  CholeskyQR factorization.  We assume $W = 2mn/P$ in the sequential
  case, where $P$ is the number of blocks and $W$ is the number of
  floating-point words that fit in fast memory.  Lower-order terms
  omitted.  All parallel terms are counted along the critical path.}
\label{tbl:CholeskyQR:counts}
\end{table}

CholeskyQR begins by computing half of the symmetric matrix $A^T A$.
In the parallel case, each processor $i$ computes half of its
component $A_i^T A_i$ locally.  In the sequential case, this happens
one block at a time.  Since this result is a symmetric $n \times n$
matrix, the operation takes only $mn^2/P + O(mn/P)$ flops.  These local
components are then summed using a(n) (all-)reduction, which can also
exploit symmetry.  The final operation, the Cholesky factorization,
requires $n^3/3 + O(n^2)$ flops.  (Choosing a more stable or robust
factorization does not improve the accuracy bound, as the accuracy has
already been lost by computing $A^T A$.)  Finally, the $Q := A L^{-T}$
operation costs $mn^2/P + O(mn/P)$ flops per block of $A$.  Table
\ref{tbl:CholeskyQR:counts} summarizes both the parallel and
sequential performance models.  In Section \ref{S:TSQR:stability}, we
compare the accuracy of CholeskyQR to that of TSQR and other ``tall
skinny'' QR factorization algorithms.

\subsection{Householder QR}\label{SS:TSQR:perfcomp:HQR}
\label{SS:TSQR:perfcomp:ScaLAPACK} 

Householder QR uses orthogonal reflectors to reduce a matrix to upper
tridiagonal form, one column at a time (see e.g., \cite{govl:96}).  In
the current version of LAPACK and ScaLAPACK, the reflectors are
coalesced into block columns (see e.g., \cite{schreiber1989storage}).
This makes trailing matrix updates more efficient, but the panel
factorization is still standard Householder QR, which works one column
at a time.  These panel factorizations are an asymptotic latency
bottleneck in the parallel case, especially for tall and skinny
matrices.  Thus, we model parallel Householder QR without considering
block updates.  In contrast, we will see that operating on blocks of
columns can offer asymptotic bandwidth savings in sequential
Householder QR, so it pays to model a block column version.

\subsubsection{Parallel Householder QR}

ScaLAPACK's parallel QR factorization routine, \lstinline!PDGEQRF!,
uses a right-looking Householder QR approach \cite{lawn80}.  The cost
of \lstinline!PDGEQRF! depends on how the original matrix $A$ is
distributed across the processors.  For comparison with TSQR, we
assume the same block row layout on $P$ processors.

\lstinline!PDGEQRF! computes an explicit representation of the $R$
factor, and an implicit representation of the $Q$ factor as a sequence
of Householder reflectors.  The algorithm overwrites the upper
triangle of the input matrix with the $R$ factor.  Thus, in our case,
the $R$ factor is stored only on processor zero, as long as $m/P \geq
n$.  We assume $m/P \geq n$ in order to simplify the performance
analysis.

Section \ref{SS:TSQR:localQR:BLAS3structured} describes BLAS 3
optimizations for Householder QR.  \lstinline!PDGEQRF! exploits these
techniques in general, as they accelerate the trailing matrix updates.
We do not count floating-point operations for these optimizations
here, since they do nothing to improve the latency bottleneck in the
panel factorizations.

In \lstinline!PDGEQRF!, some processors may need to perform fewer
flops than other processors, because the number of rows in the current
working column and the current trailing matrix of $A$ decrease by one
with each iteration.  With the assumption that $m/P \geq n$, however,
all but the first processor must do the same amount of work at each
iteration.  In the tall skinny regime, ``flops on the critical path''
(which is what we count) is a good approximation of ``flops on each
processor.''  We count floating-point operations, messages, and words
transferred by parallel Householder QR on general matrix layouts in
Section \ref{S:CAQR-counts}; in particular, Equation
\eqref{Eq:ScaLAPACK:time} in that section gives a performance model.

\begin{table}
  \centering
  \begin{tabular}{l|c|c|c}
    Parallel algorithm & \# flops & \# messages & \# words \\ \hline
    TSQR & $\frac{2mn^2}{P} + \frac{2n^3}{3} \log(P)$
         & $\log(P)$
         & $\frac{n^2}{2} \log(P)$ \\
    \lstinline!PDGEQRF!
        & $\frac{2mn^2}{P} - \frac{2n^3}{3P}$ 
        & $2n \log(P)$
        & $\frac{n^2}{2} \log(P)$ \\
    MGS\_R & $\frac{2mn^2}{P}$
          & $2n \log(P)$
          & $\frac{n^2}{2}\log(P)$ \\
    CGS\_L & $\frac{2mn^2}{P}$
          & $2n \log(P)$ 
          & $\frac{n^2}{2}\log(P)$ \\ 
    CholeskyQR & $\frac{2mn^2}{P} + \frac{n^3}{3}$ 
               & $\log(P)$ 
               & $\frac{n^2}{2}\log(P)$ \\ 
  \end{tabular}
  \caption{Performance model of various parallel QR factorization
    algorithms.  ``CGS2'' means CGS with one reorthogonalization pass.
    Lower-order terms omitted.  All parallel terms are counted along
    the critical path.  We show only the best-performing versions of
    MGS and CGS.  We omit CGS2 because it is no slower than applying 
    CGS twice, but the number of orthogonalization steps may vary 
    based on the numerical properties of the input, so it is hard to 
    predict performance \emph{a priori}.}
  \label{tbl:QR:perfcomp:par}
\end{table}


Table \ref{tbl:QR:perfcomp:par} compares the performance of all the
parallel QR factorizations discussed here.  We see that 1-D TSQR and
CholeskyQR save both messages and bandwidth over MGS\_R and
ScaLAPACK's \lstinline!PDGEQRF!, but at the expense of a higher-order
$n^3$ flops term.

\subsubsection{Sequential Householder QR}

LAPACK Working Note \#118 describes a left-looking out-of-DRAM QR
factorization \lstinline!PFDGEQRF!, which is implemented as an
extension of ScaLAPACK \cite{dazevedo1997design}.  It uses ScaLAPACK's
parallel QR factorization \lstinline!PDGEQRF! to perform the current
panel factorization in DRAM.  Thus, it is able to exploit parallelism.
We assume here, though, that it is running sequentially, since we are
only interested in modeling the traffic between slow and fast memory.
\lstinline!PFDGEQRF! is a left-looking method, as usual with
out-of-DRAM algorithms.  The code keeps two panels in memory: a left
panel of fixed width $b$, and the current panel being factored, whose
width $c$ can expand to fill the available memory.  Appendix
\ref{S:PFDGEQRF} describes the method in more detail with performance
counts, and Algorithm \ref{Alg:PFDGEQRF:outline} in the Appendix gives
an outline of the code.

See Equation \eqref{eq:PFDGEQRF:runtime:W} in Appendix
\ref{S:PFDGEQRF} for the following counts.  The \lstinline!PFDGEQRF!
algorithm performs
\[
    2mn^2 - \frac{2n^3}{3}
\]
floating-point arithmetic operations, just like any sequential
Householder QR factorization.  (Here and elsewhere, we omit
lower-order terms.)  It transfers a total of about
\[
\frac{m^2 n^2}{2W}
- \frac{m n^3}{6W}
+ \frac{3mn}{2} 
- \frac{3n^2}{4} 
\]
floating-point words between slow and fast memory, and accesses slow
memory (counting both reads and writes) about
\[
\frac{mn^2}{2W} 
+ \frac{2mn}{W} 
- \frac{n}{2}
\]
times.  In contrast, sequential TSQR only requires 
\[
\frac{2mn}{\tilde{W}}
\]
slow memory accesses, where $\tilde{W} = W - n(n+1)/2$, and only
transfers
\[
2mn 
- \frac{n(n+1)}{2} 
+ \frac{mn^2}{\tilde{W}}
\]
words between slow and fast memory (see Equation
\eqref{eq:TSQR:seq:modeltimeW:factor} in Appendix
\ref{S:TSQR-seq-detailed}).  We note that we expect $W$ to be a
reasonably large multiple of $n^2$, so that $\tilde{W} \approx W$.

Table \ref{tbl:QR:perfcomp:seq} compares the performance of the
sequential QR factorizations discussed in this section, including our
modeled version of \lstinline!PFDGEQRF!.

\begin{table}
\small
  \centering
  \begin{tabular}{l|c|c|c}
    Sequential algorithm & \# flops & \# messages & \# words  \\\hline
    TSQR & $2mn^2 - \frac{2n^3}{3}$ 
         & $\frac{2mn}{\tilde{W}}$
         & $2mn - \frac{n(n+1)}{2} 
           + \frac{mn^2}{\tilde{W}}$ \\
    \lstinline!PFDGEQRF!
         & $2mn^2 - \frac{2 n^3}{3}$
         & $\frac{2mn}{W} + \frac{mn^2}{2W}$
         & $\frac{m^2 n^2}{2W} - \frac{mn^3}{6W}
            + \frac{3mn}{2} - \frac{3n^2}{4}$ \\
    MGS & $2mn^2$ 
        & $\frac{2mn^2}{\tilde{W}}$
        & $\frac{3mn}{2} + \frac{m^2 n^2}{2 \tilde{W}}$ \\
    CholeskyQR & $2mn^2 + \frac{n^3}{3}$
               & $\frac{6mn}{W}$ 
               & $3mn$ \\
  \end{tabular}
  \caption{Performance model of various sequential QR factorization
    algorithms.  \lstinline!PFDGEQRF! is our model of ScaLAPACK's
    out-of-DRAM QR factorization; $W$ is the fast memory size, and
    $\tilde{W} = W - n(n+1)/2$.  Lower-order terms omitted.  We omit
    CGS2 because it is no slower than applying CGS twice, but the
    number of orthogonalization steps may vary based on the numerical
    properties of the input, so it is hard to predict performance
    \emph{a priori}.}
  \label{tbl:QR:perfcomp:seq}
\end{table}

\section{Numerical stability of TSQR and other QR
  factorizations}\label{TSQR:stability}\label{S:TSQR:stability}

In the previous section, we modeled the performance of various QR
factorization algorithms for tall and skinny matrices on a block row
layout.  Our models show that CholeskyQR should have better
performance than all the other methods.  However, numerical accuracy
is also an important consideration for many users.  For example, in
CholeskyQR, the loss of orthogonality of the computed $Q$ factor
depends quadratically on the condition number of the input matrix (see
Table \ref{tbl:TSQR:stability}).  This is because computing the Gram
matrix $A^T A$ squares the condition number of $A$.  One can avoid
this stability loss by computing and storing $A^T A$ in doubled
precision.  However, this doubles the communication volume.  It also
increases the cost of arithmetic operations by a hardware-dependent
factor.

\begin{table}[h]
\centering
\begin{tabular}{l|c|c|c} 
  Algorithm       & $\| I - Q^T Q \|_2$ bound 
                  & Assumption on $\kappa(A)$ & Reference(s) \\ \hline
  Householder QR  & $O(\varepsilon)$ & None & \cite{govl:96} \\
  TSQR            & $O(\varepsilon)$ & None & \cite{govl:96} \\
  CGS2            & $O(\varepsilon)$ 
                  & $O(\varepsilon \kappa(A)) < 1$ 
                  & \cite{abde:71,kiel:74} \\ 
  MGS             & $O(\varepsilon \kappa(A))$   
                  & None
                  & \cite{bjor:67} \\ 
  CholeskyQR      & $O(\varepsilon \kappa(A)^2)$ 
                  & None
                  & \cite{stwu:02} \\
  CGS             & $O(\varepsilon \kappa(A)^{n-1})$ 
                  & None
                  & \cite{kiel:74,smbl:06} \\ 
\end{tabular}
\caption{Upper bounds on deviation from orthogonality of
  the $Q$ factor from various QR algorithms.  Machine precision is 
  $\varepsilon$.  ``Assumption on $\kappa(A)$'' refers to any 
  constraints which $\kappa(A)$ must satisfy in order for the bound 
  in the previous column to hold.}
\label{tbl:TSQR:stability}
\end{table}

Unlike CholeskyQR, CGS, or MGS, Householder QR is
\emph{unconditionally stable}.  That is, the computed $Q$ factors are
always orthogonal to machine precision, regardless of the properties
of the input matrix \cite{govl:96}.  This also holds for TSQR, because
the algorithm is composed entirely of no more than $P$ Householder QR
factorizations, in which $P$ is the number of input blocks.  Each of
these factorizations is itself unconditionally stable.  In contrast,
the orthogonality of the $Q$ factor computed by CGS, MGS, or
CholeskyQR depends on the condition number of the input matrix.
Reorthogonalization in MGS and CGS can make the computed $Q$ factor
orthogonal to machine precision, but only if the input matrix $A$ is
numerically full rank, i.e., if $O(\varepsilon \kappa(A)) < 1$.
Reorthogonalization also doubles the cost of the algorithm.

However, sometimes some loss of accuracy can be tolerated, either to
improve performance, or for the algorithm to have a desirable
property.  For example, in some cases the input vectors are
sufficiently well-conditioned to allow using CholeskyQR, and the
accuracy of the orthogonalization is not so important.  Another
example is GMRES.  Its backward stability was proven first for
Householder QR orthogonalization, and only later for modified
Gram-Schmidt orthogonalization \cite{greenbaum1997numerical}.  Users
traditionally prefer the latter formulation, mainly because the
Householder QR version requires about twice as many floating-point
operations (as the $Q$ matrix must be computed explicitly).  Another
reason is that most GMRES descriptions make the vectors available for
orthogonalization one at a time, rather than all at once, as
Householder QR would require (see e.g.,
\cite{walker1985implementation}).  (Demmel et al.\ review existing
techniques and present new methods for rearranging GMRES and other
Krylov subspace methods for use with Householder QR and TSQR
\cite{demmel2008comm}.)

We care about stability for two reasons.  First, an important
application of TSQR is the orthogonalization of basis vectors in
Krylov methods.  When using Krylov methods to compute eigenvalues of
large, ill-conditioned matrices, the whole solver can fail to converge
or have a considerably slower convergence when the orthogonality of
the Ritz vectors is poor \cite{lehoucqORTH,andrewORTH}.  Second, we
will use TSQR in Section \ref{S:CAQR} as the panel factorization in a
QR decomposition algorithm for matrices of general shape.  Users who
ask for a QR factorization generally expect it to be numerically
stable.  This is because of their experience with Householder QR,
which does more work than LU or Cholesky, but produces more accurate
results.  Users who are not willing to spend this additional work
already favor faster but less stable algorithms.

Table \ref{tbl:TSQR:stability} summarizes known upper bounds on the
deviation from orthogonality $\|I - Q^T Q\|_2$ of the computed $Q$
factor, as a function of the machine precision $\varepsilon$ and the
input matrix's two-norm condition number $\kappa(A)$, for various QR
factorization algorithms.  Except for CGS, all these bounds are sharp.
Smoktunowicz et al.\ demonstrate a matrix satisfying $O(\varepsilon
\kappa(A)^2) < 1$ for which $\|I - Q^T Q\|_2$ is not $O(\varepsilon
\kappa(A)^2)$, but as far as we know, no matrix has yet been found for
which the $\|I - Q^T Q\|_2$ is $O(\varepsilon \kappa(A)^{n-1})$ bound
is sharp \cite{smbl:06}.

In the table, ``CGS2'' refers to classical Gram-Schmidt with one
reorthogonalization pass.  A single reorthgonalization pass suffices
to make the $Q$ factor orthogonal to machine precision, as long as the
input matrix is numerically full rank, i.e., if $O(\varepsilon
\kappa(A)) < 1$.  This is the source of Kahan's maxim, ``Twice is
enough'' \cite{parlett1998symmetric}: the accuracy reaches its
theoretical best after one reorthogonalization pass (see also
\cite{abde:71}), and further reorthogonalizations do not improve
orthogonality.  However, TSQR needs only half as many messages to do
just as well as CGS2.  In terms of communication, TSQR's stability
comes for free.


\section{Platforms of interest for TSQR experiments and models}\label{S:TSQR:platforms}

\subsection{A large, but composable tuning space}

TSQR is not a single algorithm, but a space of possible algorithms.
It encompasses all possible reduction tree shapes, including:
\begin{enumerate}
  \item Binary (to minimize number of messages in the parallel case)
  \item Flat (to minimize communication volume in the sequential
    case)
  \item Hybrid (to account for network topology, and/or to balance
    bandwidth demands with maximum parallelism)
\end{enumerate}
as well as all possible ways to perform the local QR factorizations,
including:
\begin{enumerate}
  \item (Possibly multithreaded) standard LAPACK (DGEQRF)
  \item An existing parallel QR factorization, such as ScaLAPACK's
    PDGEQRF 
  \item A ``divide-and-conquer'' QR factorization (e.g.,
    \cite{elmroth1998new})
  \item Recursive (invoke another form of TSQR)
\end{enumerate}
Choosing the right combination of parameters can help minimize
communication between any or all of the levels of the memory
hierarchy, from cache and shared-memory bus, to DRAM and local disk,
to parallel filesystem and distributed-memory network interconnects,
to wide-area networks.

The huge tuning space makes it a challenge to pick the right platforms
for experiments.  Luckily, TSQR's hierarchical structure makes tunings
\emph{composable}.  For example, once we have a good choice of
parameters for TSQR on a single multicore node, we don't need to
change them when we tune TSQR for a cluster of these nodes.  From the
cluster perspective, it's as if the performance of the individual
nodes improved.  This means that we can benchmark TSQR on a small,
carefully chosen set of scenarios, with confidence that they represent
many platforms of interest.

\subsection{Platforms of  interest}\label{SS:TSQR:platforms:interest}
\label{TSQR:perfres:platforms-of-interest}

Here we survey a wide variety of interesting platforms for TSQR, and
explain the key features of each that we will distill into a small
collection of experiments.

\subsubsection{Single-node parallel, and explicitly swapping}

The ``cluster of shared-memory parallel (SMP) nodes'' continues to
provide a good price-performance point for many large-scale
applications.  This alone would justify optimizing the single-node
case.  Perhaps more importantly, the ``multicore revolution'' seeks to
push traditionally HPC applications into wider markets, which favor
the single-node workstation or even the laptop over the expensive,
power-hungry, space-consuming, difficult-to-maintain cluster.  A large
and expanding class of users may never run their jobs on a cluster.

%
%
Multicore SMPs can help reduce communication costs, but cannot
eliminate them.  TSQR can exploit locality by sizing individual
subproblems to fit within any level of the memory hierarchy.  This
gives programmers explicit control over management of communication
between levels, much like a traditional ``out-of-core''
algorithm.\footnote{We avoid this label because it's an anachronism
  (``core'' refers to main system memory, constructed of solenoids
  rather than transistors or DRAM), and because people now easily
  confuse ``core'' with ``processing unit'' (in the sense of
  ``multicore'').  We prefer the more precise term \emph{explicitly
    swapping}, or ``out-of-X'' for a memory hierarchy level X.  For
  example, ``out-of-DRAM'' means using a disk, flash drive, or other
  storage device as swap space for problems too large to solve in
  DRAM.}  TSQR's hierarchical structure makes explicit swap management
easy; it's just another form of communication.  It gives us an
optimized implementation for ``free'' on platforms like Cell or GPUs,
which require explicitly moving data into separate storage areas for
processing.  Also, it lets us easily and efficiently solve problems
too large to fit in DRAM.  This seems like an old-fashioned issue,
since an individual node nowadays can accommodate as much as 16 GB of
DRAM.  Explicitly swapping variants of libraries like ScaLAPACK tend
to be ill-maintained, due to lack of interest.  However, we predict a
resurgence of interest in explicitly-swapping algorithms, for the
following reasons:
\begin{itemize}
  \item Single-node workstations will become more popular than
    multinode clusters, as the number of cores per node increases.

  \item The amount of DRAM per node cannot scale linearly with the
    number of cores per node, because of DRAM's power requirements.
    Trying to scale DRAM will wipe out the power savings promised by
    multicore parallelism.  
  \item The rise to prominence of mobile computing -- e.g., more
    laptops than desktops were sold in U.S.\ retail in 2006 --
    drives increasing concern for total-system power use.

  \item Most operating systems do not treat ``virtual memory'' as
    another level of the memory hierarchy.  Default and recommended
    configurations for Linux, Windows XP, Solaris, and AIX on modern
    machines assign only 1--3 times as much swap space as DRAM, so
    it's not accurate to think of DRAM as a cache for disk.  Few
    operating systems expand swap space on demand, and expanding it
    manually generally requires administrative access to the machine.
    It's better for security and usability to ask applications to adapt
    to the machine settings, rather than force users to change their
    machine for their applications.

  \item \emph{Unexpected} use of virtual memory swap space generally
    slows down applications by orders of magnitude.  HPC programmers
    running batch jobs consider this a performance problem serious
    enough to warrant terminating the job early and sizing down the
    problem.  Users of interactive systems typically experience large
    (and often frustrating) delays in whole-system responsiveness when
    extensive swapping occurs.

  \item In practice, a single application need not consume all memory
    in order to trigger the virtual memory system to swap extensively.

  \item Explicitly swapping software does not stress the OS's virtual
    memory system, and can control the amount of memory and disk
    bandwidth resources that it uses.

  \item Alternate storage media such as solid-state drives offer more
    bandwidth than traditional magnetic hard disks.  Typical hard disk
    read or write bandwidth as of this work's publication date is
    around 60 MB/s, whereas Samsung announced in May 2008 the upcoming
    release of a 256 GB capacity solid-state drive with 200 MB/s read
    bandwidth and 160 MB/s write bandwidth \cite{samsung2008ssd}.
    Solid-state drives are finding increasing use in mobile devices
    and laptops, due to their lower power requirements.  This will
    make out-of-DRAM applications more attractive by widening any
    bandwidth bottlenecks.
\end{itemize}

\subsubsection{Distributed-memory machines}

Avoiding communication is a performance-enhancing strategy for
distributed-memory architectures as well.  TSQR can improve
performance on traditional clusters as well as other networked
systems, such as grid and perhaps even volunteer computing.  Avoiding
communication also makes improving network reliability less of a
performance burden, as software-based reliability schemes use some
combination of redundant and/or longer messages.  Many
distributed-memory supercomputers have high-performance parallel
filesystems, which increase the bandwidth available for out-of-DRAM
TSQR.  This enables reducing per-node memory requirements without
increasing the number of nodes needed to solve the problem.

\subsection{Pruning the platform space}\label{SS:TSQR:platforms:pruning}

For single-node platforms, we think it pays to investigate both
problems that fit in DRAM (perhaps with explicit cache management),
and problems too large to fit in DRAM, that call for explicit swapping
to a local disk.  High-performance parallel filesystems offer
potentially much more bandwidth, but we chose not to use them for our
experiments for the following reasons:
\begin{itemize}
  \item Lack of availability of a single-node machine with exclusive
    access to a parallel filesystem
  \item On clusters, parallel filesystems are usually shared among all
    cluster users, which would make it difficult to collect repeatable
    timings.
\end{itemize}
For multinode benchmarks, we opted for traditional clusters rather
than volunteer computing, due to the difficulty of obtaining
repeatable timings in the latter case.

\subsection{Platforms for experiments}\label{SS:TSQR:platforms:exp}

We selected the following experiments as characteristic of the space
of platforms:
\begin{itemize}
\item Single node, sequential, out-of-DRAM, and
\item Distributed memory, in-DRAM on each node.
\end{itemize}

We ran sequential TSQR on a laptop with a single PowerPC CPU.  It
represents the embedded and mobile space, with its tighter power and
heat requirements.  Details of the platform are as follows:
\begin{itemize}
\item Single-core PowerPC G4 (1.5 GHz) 
\item 512 KB of L2 cache
\item 512 MB of DRAM on a 167 MHz bus
\item One Fujitsu MHT2080AH 80 HB hard drive (5400 RPM) 
\item MacOS X 10.4.11
\item GNU C Compiler (gcc), version 4.0.1
\item vecLib (Apple's optimized dense linear algebra library),
  version 3.2.2
\end{itemize}

We ran parallel TSQR on the following distributed-memory machines:
\begin{enumerate}
\item Pentium III cluster (``Beowulf'')
\begin{itemize}
\item Operated by the University of Colorado at Denver and the Health
  Sciences Center
\item 35 dual-socket 900 MHz Pentium III nodes with Dolphin
  interconnect
\item Floating-point rate:  900 Mflop/s per processor, peak
\item Network latency: less than 2.7 $\mu$s, benchmarked\footnote{See
  \url{http://www.dolphinics.com/products/benchmarks.html}.}
\item Network bandwidth: 350 MB/s, benchmarked upper bound
\end{itemize}

\item IBM BlueGene/L (``Frost'')
\begin{itemize}
\item Operated by the National Center for Atmospheric Research
\item One BlueGene/L rack with 1024 700 MHz compute CPUs
\item Floating-point rate:  2.8 Gflop/s per processor, peak
\item Network\footnote{The BlueGene/L has two separate networks -- a
    torus for nearest-neighbor communication and a tree for
    collectives.  The latency and bandwidth figures here are for the
    collectives network.} latency:  1.5 $\mu$s, hardware
\item Network one-way bandwidth:  350 MB/s, hardware
\end{itemize}
\end{enumerate}

\subsection{Platforms for performance models}\label{SS:TSQR:platforms:models}

In Section \ref{S:CAQR:perfest}, we estimate performance of CAQR, our
QR factorization algorithm on a 2-D matrix layout, on three different
parallel machines: an existing IBM POWER5 cluster with a total of 888
processors (``IBM POWER5''), a future proposed petascale machine with
8192 processors (``Peta''), and a collection of 128 processors linked
together by the Internet (``Grid'').  Here are the parameters we use
in our models for the three parallel machines:
\begin{itemize}
\item IBM POWER5
\begin{itemize}
  \item 888 processors
  \item Floating-point rate:  7.6 Gflop/s per processor, peak
  \item Network latency:  5 $\mu$s
  \item Network bandwidth:  3.2 GB/s
\end{itemize}

\item Peta
\begin{itemize}
  \item 8192 processors
  \item Floating-point rate:  500 Gflop/s per processor, peak
  \item Network latency:  10 $\mu$s
  \item Network bandwidth:  4 GB/s
\end{itemize}

\item Grid
\begin{itemize}
  \item 128 processors
  \item Floating-point rate:  10 Tflop/s, peak
  \item Network latency:  0.1 s
  \item Network bandwidth:  0.32 GB/s
\end{itemize}
\end{itemize}
Peta is our projection of a future high-performance computing cluster,
and Grid is our projection of a collection of geographically separated
high-performance clusters, linked over a TeraGrid-like backbone.  Each
``processor'' of Peta may itself be a parallel multicore node, but we
consider it as a single fast sequential processor for the sake of our
model.  Similarly, each ``processor'' of Grid is itself a cluster, but
we consider it as a single very fast sequential processor.


\section{TSQR performance results}\label{S:TSQR:perfres}


\subsection{Scenarios used in experiments}

Previous work covers some parts of the tuning space mentioned in
Section \ref{SS:TSQR:platforms:pruning}.  Gunter et al.\ implemented
an out-of-DRAM version of TSQR on a flat tree, and used a parallel
distributed QR factorization routine to factor in-DRAM blocks
\cite{gunter2005parallel}.  Pothen and Raghavan
\cite{pothen1989distributed} and Cunha et al. \cite{cunha2002new} both
benchmarked parallel TSQR using a binary tree on a distributed-memory
cluster, and implemented the local QR factorizations with a
single-threaded version of \lstinline!DGEQRF!.  All these researchers
observed significant performance improvements over previous QR
factorization algorithms.

We chose to run two sets of experiments.  The first set covers the
out-of-DRAM case on a single CPU.  The second set is like the parallel
experiments of previous authors in that it uses a binary tree on a
distributed-memory cluster, but it improves on their approach by using
a better local QR factorization (the divide-and-conquer approach --
see \cite{elmroth2000applying}).

\subsection{Sequential out-of-DRAM tests}

We developed an out-of-DRAM version of TSQR that uses a flat reduction
tree.  It invokes the system vendor's native BLAS and LAPACK
libraries.  Thus, it can exploit a multithreaded BLAS on a machine
with multiple CPUs, but the parallelism is limited to operations on a
single block of the matrix.  We used standard POSIX blocking file
operations, and made no attempt to overlap communication and
computation.  Exploiting overlap could at best double the performance.

We ran sequential tests on a laptop with a single PowerPC CPU, as
described in Section \ref{SS:TSQR:platforms:exp}.  In our experiments,
we first used both out-of-DRAM TSQR and standard LAPACK QR to factor a
collection of matrices that use only slightly more than half of the
total DRAM for the factorization.  This was so that we could collect
comparison timings.  Then, we ran only out-of-DRAM TSQR on matrices
too large to fit in DRAM or swap space, so that an out-of-DRAM
algorithm is necessary to solve the problem at all.  For the latter
timings, we extrapolated the standard LAPACK QR timings up to the
larger problem sizes, in order to estimate the runtime if memory were
unbounded.  LAPACK's QR factorization swaps so much for out-of-DRAM
problem sizes that its actual runtimes are many times larger than
these extrapolated unbounded-memory runtime estimates.  As mentioned
in Section \ref{SS:TSQR:platforms:interest}, once an in-DRAM algorithm
begins swapping, it becomes so much slower that most users prefer to
abort the computation and try solving a smaller problem.  No attempt
to optimize by overlapping communication and computation was made.

We used the following power law for the extrapolation:
\[
t = A_1 b m^{A_2} n^{A_3},
\]
in which $t$ is the time spent in computation, $b$ is the number of
input matrix blocks, $m$ is the number of rows per block, and $n$ is
the number of columns in the matrix.  After taking logarithms of both
sides, we performed a least squares fit of $\log(A_1)$, $A_2$, and
$A_3$.  The value of $A_2$ was 1, as expected.  The value of $A_3$ was
about 1.6.  This is less than 2 as expected, given that increasing the
number of columns increases the computational intensity and thus the
potential for exploitation of locality (a BLAS 3 effect).  We expect
around two digits of accuracy in the parameters, which in themselves
are not as interesting as the extrapolated runtimes; the parameter
values mainly serve as a sanity check.

\subsubsection{Results}

\begin{figure}
  \begin{center}        
    \includegraphics[angle=0,scale=0.6]{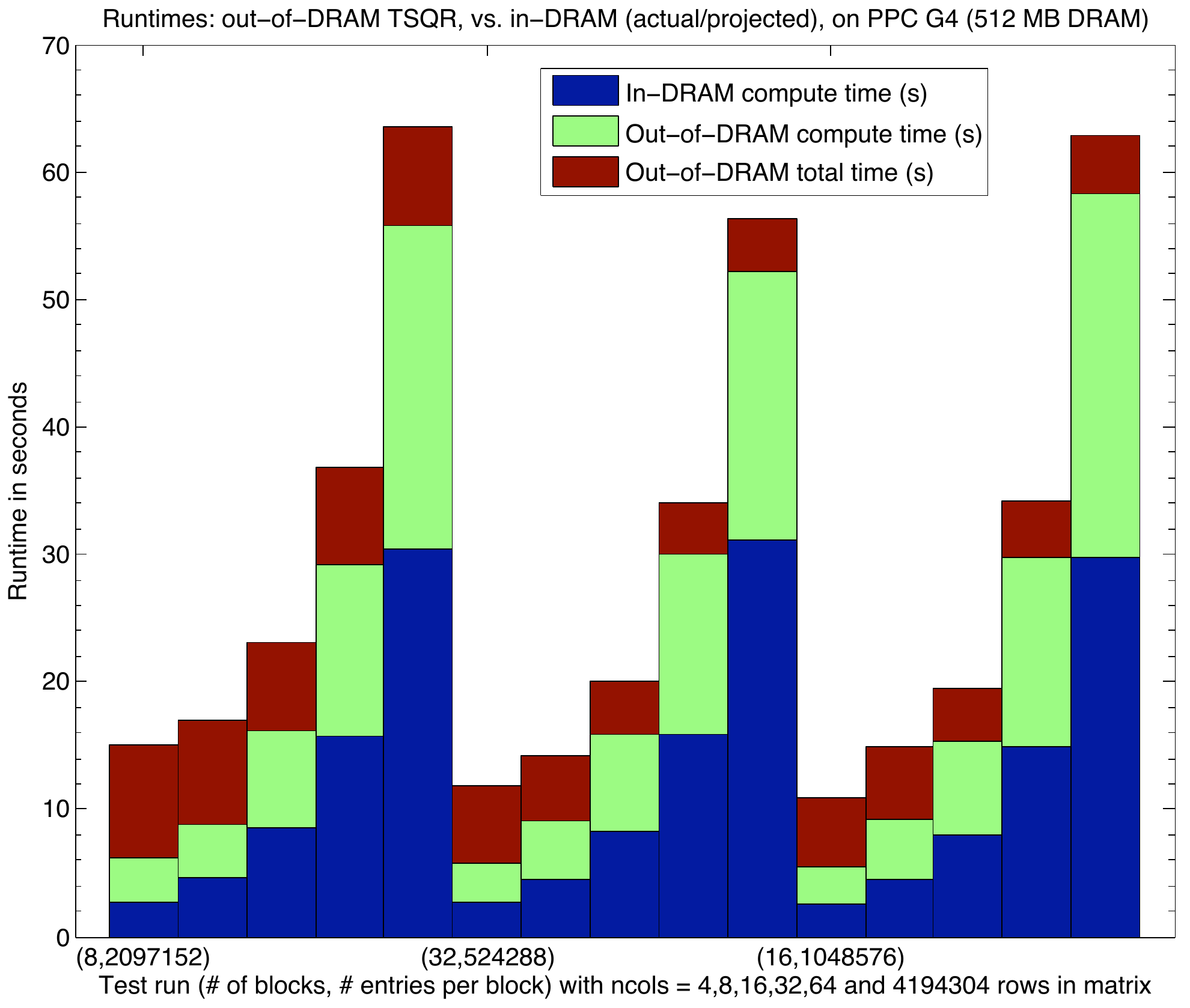}
    \caption{\small Runtimes (in seconds) of out-of-DRAM TSQR and
      standard QR (LAPACK's \lstinline!DGEQRF!) on a single-processor
      laptop.  All data is measured.  We limit memory usage to 256 MB,
      which is half of the laptop's total system memory, so that we
      can collect performance data for DGEQRF.  The graphs show
      different choices of block dimensions and number of blocks.  The
      top of the blue bar is the benchmarked total runtime for
      \lstinline!DGEQRF!, the top of the green bar is the benchmarked
      compute time for TSQR, and the top of the brown bar is the
      benchmarked total time for TSQR.  Thus the height of the brown
      bar alone is the I/O time.  Note that LAPACK starts and ends in
      DRAM, and TSQR starts and ends on disk.}
    \label{fig:tsqr-laptop-no-extrap}
  \end{center}
\end{figure}

\begin{figure}
  \begin{center}        
    \includegraphics[angle=0,scale=0.6]{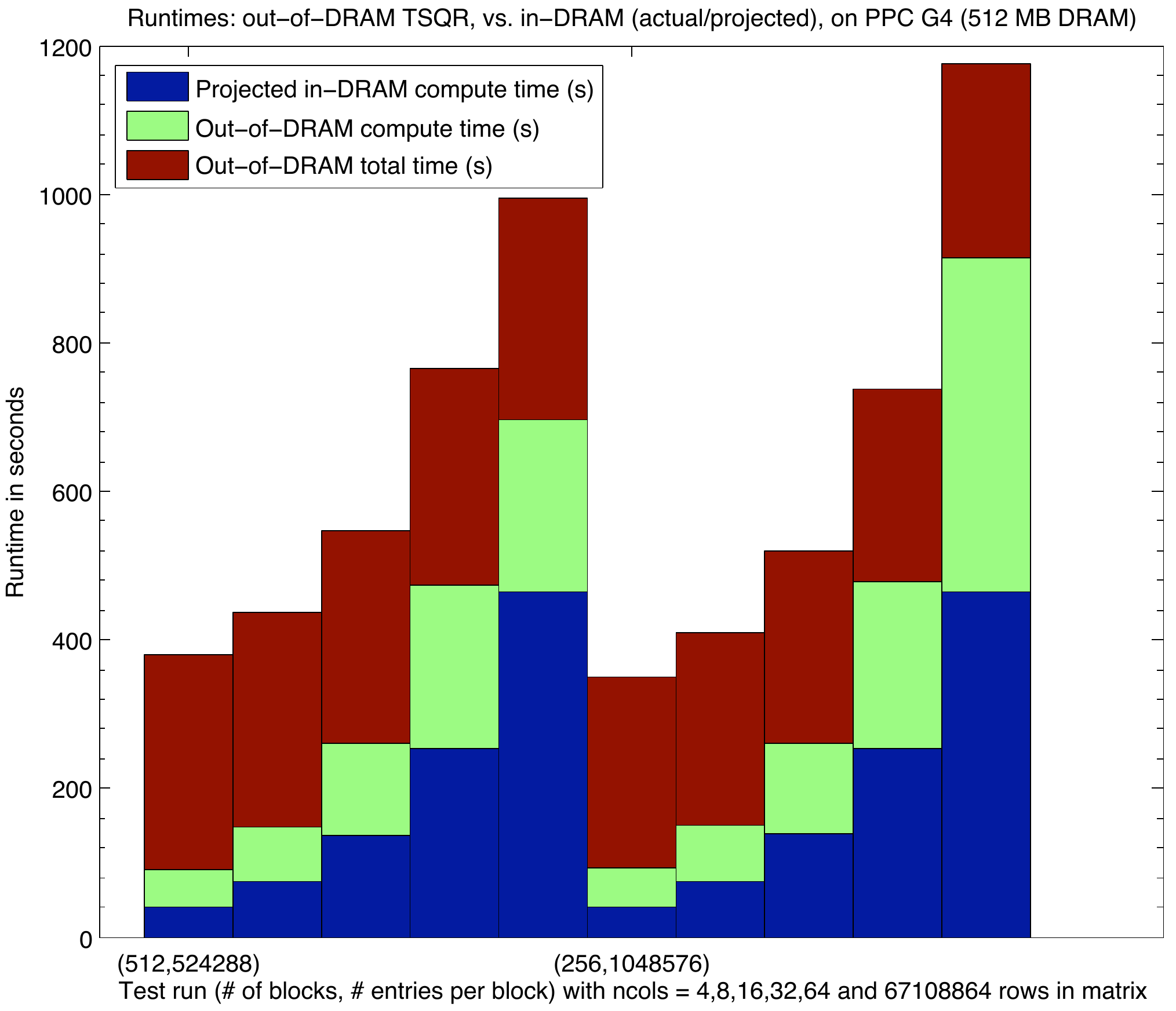}
    \caption{\small Measured runtime (in seconds) of out-of-DRAM TSQR,
      compared against extrapolated runtime (in seconds) of standard
      QR (LAPACK's DGEQRF) on a single-processor laptop.  We use the
      data in Figure \ref{fig:tsqr-laptop-no-extrap} to construct a
      power-law performance extrapolation.  The graphs show different
      choices of block dimensions and number of blocks.  The top of
      the blue bar is the extrapolated total runtime for
      \lstinline!DGEQRF!, the top of the green bar is the benchmarked
      compute time for TSQR, and the top of the brown bar is the
      benchmarked total time for TSQR.  Thus the height of the brown
      bar alone is the I/O time.  Note that LAPACK starts and ends in
      DRAM (if it could fit in DRAM), and TSQR starts and ends on
      disk.}
    \label{fig:tsqr-laptop-extrap}
  \end{center}
\end{figure}

Figure \ref{fig:tsqr-laptop-no-extrap} shows the measured in-DRAM
results on the laptop platform, and Figure
\ref{fig:tsqr-laptop-extrap} shows the (measured TSQR, extrapolated
LAPACK) out-of-DRAM results on the same platform.  In these figures,
the number of blocks used, as well as the number of elements in the
input matrix (and thus the total volume of communication), is the same
for each group of five bars.  We only varied the number of blocks and
the number of columns in the matrix.  For each graph, the total number
of rows in the matrix is constant for all groups of bars.  Note that
we have not tried to overlap I/O and computation in this
implementation.  The trends in Figure \ref{fig:tsqr-laptop-no-extrap}
suggest that the extrapolation is reasonable: TSQR takes about twice
as much time for computation as does standard LAPACK QR, and the
fraction of time spent in I/O is reasonable and decreases with problem
size.  

TSQR assumes that the matrix starts and ends on disk, whereas LAPACK
starts and ends in DRAM.  Thus, to compare the two, one could also
estimate LAPACK performance with infinite DRAM but where the data
starts and ends on disk.  The height of the reddish-brown bars in
Figures \ref{fig:tsqr-laptop-no-extrap} and
\ref{fig:tsqr-laptop-extrap} is the I/O time for TSQR, which can be
used to estimate the LAPACK I/O time.  Add this to the blue bar (the
LAPACK compute time) to estimate the runtime when the LAPACK QR
routine must load the matrix from disk and store the results back to
disk.

\subsubsection{Conclusions}

The main purpose of our out-of-DRAM code is not to outperform existing
in-DRAM algorithms, but to be able to solve classes of problems which
the existing algorithms cannot solve.  The above graphs show that the
penalty of an explicitly swapping approach is about 2x, which is small
enough to warrant its practical use.  This holds even for problems
with a relatively low computational intensity, such as when the input
matrix has very few columns.  Furthermore, picking the number of
columns sufficiently large may allow complete overlap of file I/O by
computation.

\subsection{Parallel cluster tests}

\begin{table}
{\small
\begin{tabular}{c|r|r|r|r|r|r} 
\# procs & CholeskyQR & TSQR & CGS & MGS & TSQR & ScaLAPACK \\
 & & (\lstinline!DGEQR3!) & & & (\lstinline!DGEQRF!) & \\ \hline
  1 & 1.02 & 4.14 &  3.73 &  7.17 &  9.68 & 12.63 \\ 
  2 & 0.99 & 4.00 &  6.41 & 12.56 & 15.71 & 19.88 \\ 
  4 & 0.92 & 3.35 &  6.62 & 12.78 & 16.07 & 19.59 \\ 
  8 & 0.92 & 2.86 &  6.87 & 12.89 & 11.41 & 17.85 \\ 
 16 & 1.00 & 2.56 &  7.48 & 13.02 &  9.75 & 17.29 \\ 
 32 & 1.32 & 2.82 &  8.37 & 13.84 &  8.15 & 16.95 \\ 
 64 & 1.88 & 5.96 & 15.46 & 13.84 &  9.46 & 17.74 \\ 
\end{tabular}
}
\caption{Runtime in seconds of various parallel QR factorizations on
  the Beowulf machine.  The total number of rows $m = 100000$ and the 
  ratio $\lceil n / \sqrt{P} \rceil = 50$ (with $P$ being the number of 
  processors) were kept constant as $P$ varied from 1 to 64.  This
  illustrates weak scaling with respect to the square of the number of 
  columns $n$ in the matrix, which is of interest because the number
  of floating-point operations in sequential QR is $\Theta(mn^2)$.  If 
  an algorithm scales perfectly, then all the numbers in that
  algorithm's column should be constant.  Both the $Q$ 
  and $R$ factors were computed explicitly; in particular, for those 
  codes which form an implicit representation of $Q$, the conversion
  to an explicit representation was included in the runtime measurement.}
\label{tbl:TSQR:cluster:weak:n}
\end{table}

\begin{table}
{\small
\begin{tabular}{c|r|r|r|r|r|r}
\# procs & CholeskyQR & TSQR & CGS & MGS & TSQR & ScaLAPACK \\
 & & (\lstinline!DGEQR3!) & & & (\lstinline!DGEQRF!) & \\ \hline
  1 & 0.45 & 3.43 & 3.61 &  7.13 &  7.07 &  7.26 \\ 
  2 & 0.47 & 4.02 & 7.11 & 14.04 & 11.59 & 13.95 \\ 
  4 & 0.47 & 4.29 & 6.09 & 12.09 & 13.94 & 13.74 \\ 
  8 & 0.50 & 4.30 & 7.53 & 15.06 & 14.21 & 14.05 \\ 
 16 & 0.54 & 4.33 & 7.79 & 15.04 & 14.66 & 14.94 \\ 
 32 & 0.52 & 4.42 & 7.85 & 15.38 & 14.95 & 15.01 \\ 
 64 & 0.65 & 4.45 & 7.96 & 15.46 & 14.66 & 15.33 \\ 
\end{tabular}
}
\caption{Runtime in seconds of various parallel QR factorizations on
  the Beowulf machine, illustrating weak scaling with respect to the 
  total number of rows $m$ in the matrix.  The ratio $\lceil m/P
  \rceil = 100000$ and the total number of columns $n = 50$ were kept 
  constant as the number of processors $P$ varied from 1 to 64.  If 
  an algorithm scales perfectly, then all the numbers in that 
  algorithm's column should be constant.  For those algorithms which
  compute an implicit representation of the $Q$ factor, that
  representation was left implicit.}
\label{tbl:TSQR:cluster:weak:m}
\end{table}

\begin{table}
{\small
\begin{tabular}{c|r|r|r|r|r}
\# procs & \multicolumn{2}{c|}{TSQR} & \multicolumn{2}{c|}{ScaLAPACK} \\
         & (\lstinline!DGEQR3!) & (\lstinline!DGEQRF!) & (\lstinline!PDGEQRF!) & (\lstinline!PDGEQR2!) \\ \hline
 32  & 690 & 276 & 172 & 206 \\
 64  & 666 & 274 & 172 & 206 \\
 128 & 662 & 316 & 196 & 232 \\
 256 & 610 & 322 & 184 & 218 \\
\end{tabular}
}
\caption{Performance per processor (Mflop / s / (\# processors)) on a
  $10^6 \times 50$ matrix, on the Frost machine.  This metric
  illustrates strong scaling (constant problem size, but number of
  processors increases).  If an algorithm scales perfectly, than all 
  the numbers in that algorithm's column should be constant.  
  \lstinline!DGEQR3! is a recursive local QR factorization, and 
  \lstinline!DGEQRF! LAPACK's standard local QR factorization.}
\label{tbl:TSQR:BGL:strong}
\end{table}

We also have results from a parallel MPI implementation of TSQR on a
binary tree.  Rather than LAPACK's \lstinline!DGEQRF!, the code uses a
custom local QR factorization, \lstinline!DGEQR3!, based on the
recursive approach of Elmroth and Gustavson
\cite{elmroth2000applying}.  Tests show that \lstinline!DGEQR3!
consistently outperforms LAPACK's \lstinline!DGEQRF! by a large margin
for matrix dimensions of interest.  

We ran our experiments on two platforms: a Pentium III cluster
(``Beowulf'') and on a BlueGene/L (``Frost''), both described in
detail in Section \ref{SS:TSQR:platforms:exp}.  The experiments
compare many different implementations of a parallel QR factorization.
``CholeskyQR'' first computes the product $A^T A$ using a reduction,
then performs a QR factorization of the product.  It is less stable
than TSQR, as it squares the condition number of the original input
matrix (see Table \ref{tbl:TSQR:stability} in Section
\ref{S:TSQR:stability} for a stability comparison of various QR
factorization methods).  TSQR was tested both with the recursive local
QR factorization \lstinline!DGEQR3!, and the standard LAPACK routine
\lstinline!DGEQRF!.  Both CGS and MGS were timed.

\subsubsection{Results}

Tables \ref{tbl:TSQR:cluster:weak:n} and \ref{tbl:TSQR:cluster:weak:m}
show the results of two different performance experiments on the
Pentium III cluster.  In the first of these, the total number of rows
$m = 100000$ and the ratio $\lceil n / \sqrt{P} \rceil = 50$ (with $P$
being the number of processors) were kept constant as $P$ varied from
1 to 64.  This was meant to illustrate weak scaling with respect to
$n^2$ (the square of the number of columns in the matrix), which is of
interest because the number of floating-point operations in sequential
QR is $\Theta(mn^2)$.  If an algorithm scales perfectly, then all the
numbers in that algorithm's column should be constant.  Both the $Q$
and $R$ factors were computed explicitly; in particular, for those
codes which form an implicit representation of $Q$, the conversion to
an explicit representation was included in the runtime measurement.
The results show that TSQR scales better than CGS or MGS, and
significantly outperforms ScaLAPACK's QR.  Also, using the recursive
local QR in TSQR, rather than LAPACK's QR, more than doubles
performance.  CholeskyQR gets the best performance of all the
algorithms, but at the expense of significant loss of orthogonality.

Table \ref{tbl:TSQR:cluster:weak:m} shows the results of the second
set of experiments on the Pentium III cluster.  In these experiments,
the ratio $\lceil m/P \rceil = 100000$ and the total number of columns
$n = 50$ were kept constant as the number of processors $P$ varied
from 1 to 64.  This was meant to illustrate weak scaling with respect
to the total number of rows $m$ in the matrix.  If an algorithm scales
perfectly, then all the numbers in that algorithm's column should be
constant.  Unlike in the previous set of experiments, for those
algorithms which compute an implicit representation of the $Q$ factor,
that representation was left implicit.  The results show that TSQR
scales well.  In particular, when using TSQR with the recursive local
QR factorization, there is almost no performance penalty for moving
from one processor to two, unlike with CGS, MGS, and ScaLAPACK's QR.
Again, the recursive local QR significantly improves TSQR performance;
here it is the main factor in making TSQR perform better than
ScaLAPACK's QR.

Table \ref{tbl:TSQR:BGL:strong} shows the results of the third set of
experiments, which was performed on the BlueGene/L cluster ``Frost.''
These data show performance per processor (Mflop / s / (number of
processors)) on a matrix of constant dimensions $10^6 \times 50$, as
the number of processors was increased.  This illustrates strong
scaling.  If an algorithm scales perfectly, than all the numbers in
that algorithm's column should be constant.  Two different versions of
ScaLAPACK's QR factorization were used: \lstinline!PDGEQR2!  is the
textbook Householder QR panel factorization, and \lstinline!PDGEQRF!
is the blocked version which tries to coalesce multiple trailing
matrix updates into one.  The results again show that TSQR scales at
least as well as ScaLAPACK's QR factorization, which unlike TSQR is
presumably highly tuned on this platform.  Furthermore, using the
recursive local QR factorization with TSQR makes its performance
competitive with that of ScaLAPACK.

\subsubsection{Conclusions}

Both the Pentium III and BlueGene/L platforms have relatively slow
processors with a relatively low-latency interconnect.  TSQR was
optimized for the opposite case of fast processors and expensive
communication.  Nevertheless, TSQR outperforms ScaLAPACK's QR by
$6.8\times$ on 16 processors (and $3.5\times$ on 64 processors) on the
Pentium III cluster, and successfully competes with ScaLAPACK's QR on
the BlueGene/L machine.


\section{Parallel 2-D QR factorization}\label{S:CAQR}

The parallel CAQR (``Communication-Avoiding QR'') algorithm uses
parallel TSQR to perform a right-looking QR factorization of a dense
matrix $A$ on a two-dimensional grid of processors $P = P_r \times
P_c$.  The $m \times n$ matrix (with $m \geq n$) is distributed using
a 2-D block cyclic layout over the processor grid, with blocks of
dimension $b \times b$.  We assume that all the blocks have the same
size; we can always pad the input matrix with zero rows and columns to
ensure this is possible.  For a detailed description of the 2-D block
cyclic layout of a dense matrix, please refer to
\cite{scalapackusersguide}, in particular to the section entitled
``Details of Example Program \#1.''  There is also an analogous
sequential version of CAQR, which we summarize in Section
\ref{S:CAQR-seq} and describe in detail in Appendix
\ref{S:CAQR-seq-detailed}.  In summary, Table \ref{tbl:CAQR:par:model}
says that the number of arithmetic operations and words transferred is
roughly the same between parallel CAQR and ScaLAPACK's parallel QR
factorization, but the number of messages is a factor $b$ times lower
for CAQR.  For related work on parallel CAQR, see the second paragraph
of Section \ref{S:CAQR-seq}.

CAQR is based on TSQR in order to minimize communication.  At each
step of the factorization, TSQR is used to factor a panel of columns,
and the resulting Householder vectors are applied to the rest of the
matrix.  As we will show, the block column QR factorization as
performed in \lstinline!PDGEQRF! is the latency bottleneck of the
current ScaLAPACK QR algorithm.  Replacing this block column
factorization with TSQR, and adapting the rest of the algorithm to
work with TSQR's representation of the panel $Q$ factors, removes the
bottleneck.  We use the reduction-to-one-processor variant of TSQR, as
the panel's $R$ factor need only be stored on one processor (the pivot
block's processor).

CAQR is defined inductively.  We assume that the first $j-1$
iterations of the CAQR algorithm have been performed.  That is, $j-1$
panels of width $b$ have been factored and the trailing matrix has
been updated.  The active matrix at step $j$ (that is, the part of the
matrix which needs to be worked on) is of dimension
\[
(m - (j-1) b) \times (n - (j-1) b) = m_j \times n_j.
\]

\begin{figure}[htbp]
  \begin{center}        
    \leavevmode \includegraphics[scale=0.6]{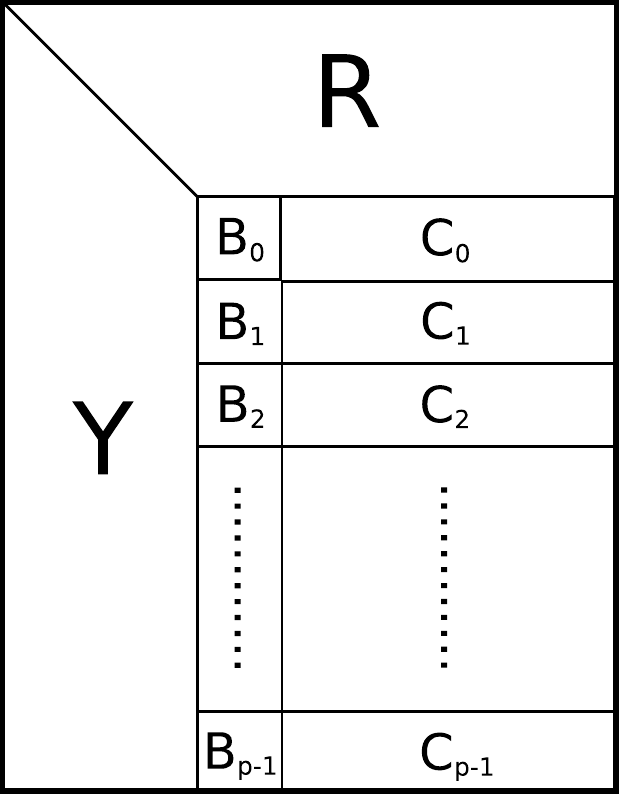}
    \caption{Step $j$ of the QR factorization algorithm.  First, the
      current panel of width $b$, consisting of the blocks $B_0$,
      $B_1$, $\dots$, $B_{p-1}$, is factorized using TSQR. Here, $p$
      is the number of blocks in the current panel.  Second, the
      trailing matrix, consisting of the blocks $C_0$, $C_1$, $\dots$,
      $C_{p-1}$, is updated.  The matrix elements above the current
      panel and the trailing matrix belong to the $R$ factor and will
      not be modified further by the QR factorization.}
    \label{fig:qr2d}
  \end{center}
\end{figure}

Figure \ref{fig:qr2d} shows the execution of the QR factorization.
For the sake of simplicity, we suppose that processors $0$, $\dots$,
$P_r - 1$ lie in the column of processes that hold the current panel $j$.
The $m_j \times b$ matrix $B$ represents the current panel $j$.  The
$m_j \times (n_j - b)$ matrix $C$ is the trailing matrix that needs to
be updated after the TSQR factorization of $B$.  For each processor
$p$, we refer to the first $b$ rows of its first block row of $B$ and
$C$ as $B_p$ and $C_p$ respectively.  

We first introduce some notation to help us refer to different parts
of a binary TSQR reduction tree.
\begin{itemize}
\item $level(i,k) = \left\lfloor \frac{i}{2^k} \right\rfloor$ denotes
  the node at level $k$ of the reduction tree which is assigned to a
  set of processors that includes processor $i$.  The initial stage of
  the reduction, with no communication, is $k = 0$.
\item $first\_proc(i,k) = 2^k level(i,k)$ is the index of the
  ``first'' processor associated with the node $level(i,k)$ at stage
  $k$ of the reduction tree.  In a reduction (not an all-reduction),
  it receives the messages from its neighbors and performs the local
  computation.
\item $target(i,k) = first\_proc(i,k) + (i + 2^{k-1}) \mod 2^k$ is the
  index of the processor with which processor $i$ exchanges data at
  level $k$ of the butterfly all-reduction algorithm.
\item $target\_first\_proc(i,k) = target(first\_proc(i,k)) =
  first\_proc(i,k) + 2^{k-1}$ is the index of the processor with which
  $first\_proc(i,k)$ exchanges data in an all-reduction at level $k$,
  or the index of the processor which sends its data to
  $first\_proc(i,k)$ in a reduction at level $k$.
\end{itemize}

Algorithm \ref{Alg:CAQR:j} outlines the right-looking parallel QR
decomposition.  At iteration $j$, first, the block column $j$ is
factored using TSQR.  We assume for ease of exposition that TSQR is
performed using a binary tree.  After the block column factorization
is complete, the matrices $C_p$ are updated as follows.  The update
corresponding to the QR factorization at the leaves of the TSQR tree
is performed locally on every processor.  The updates corresponding to
the upper levels of the TSQR tree are performed between groups of
neighboring trailing matrix processors as described in Section
\ref{SS:TSQR:localQR:trailing}.  Note that only one of the trailing
matrix processors in each neighbor group continues to be involved in
successive trailing matrix updates.  This allows overlap of
computation and communication, as the uninvolved processors can finish
their computations in parallel with successive reduction stages.

\begin{algorithm}[h!]
\caption{Right-looking parallel CAQR factorization}
\label{Alg:CAQR:j}
\begin{algorithmic}[1]
\For{$j = 1$ to $n/b$}
  \State{The column of processors that holds panel $j$ computes a TSQR
    factorization of this panel.  The Householder vectors are stored
    in a tree-like structure as described in Section 
    \ref{S:TSQR:impl}.}\label{Alg:CAQR:j:local-factor}
  \State{Each processor $p$ that belongs to the column of processes
    holding panel $j$ broadcasts along its row of processors the $m_j
    / P_r \times b$ rectangular matrix that holds the two sets of
    Householder vectors.  Processor $p$ also broadcasts two arrays of
    size $b$ each, containing the Householder multipliers $\tau_p$.}
  \State{Each processor in the same process row as processor $p$, $0
    \leq p < P_r$, forms $T_{p0}$ and updates its local trailing
    matrix $C$ using $T_{p0}$ and $Y_{p0}$. (This computation involves
    all processors.)}\label{Alg:CAQR:j:local-update}
  \For{$k = 1$ to $\log P_r$, the processors that lie in the same row 
    as processor $p$, where $0 \leq p < P_r$ equals 
    $first\_proc(p,k)$ or $target\_first\_proc(p,k)$,
    respectively.}
      \State{Processors in the same process row as
        $target\_first\_proc(p,k)$ form $T_{level(p,k),k}$ locally.  
	They also compute local pieces of 
	$W = Y_{level(p,k),k}^T C_{target\_first\_proc(p,k)}$,
        leaving the results distributed.  This computation is 
        overlapped with the communication in Line \ref{step_comm1}.}\label{Alg:CAQR:j:overlap1}

      \State{Each processor in the same process row as 
        $first\_proc(p,k)$ sends to the processor in the 
        same column and belonging to the row of processors of 
        $target\_first\_proc(p,k)$ the local pieces of
        $C_{first\_proc(p,k)}$.}\label{step_comm1}

      \State{Processors in the same process row as 
        $target\_first\_proc(p,k)$ compute local pieces of
        \[
        W = T_{level(p,k),k}^T \left( C_{first\_proc(p,k)} + W \right).
        \]}

      \State{Each processor in the same process row as 
        $target\_first\_proc(p,k)$ sends to the processor 
        in the same column and belonging to the process row
        of $first\_proc(p,k)$ the local pieces of $W$.}\label{step_comm2}

      \State{Processors in the same process row as
        $first\_proc(p,k)$ and
        $target\_first\_proc(p,k)$ each complete the rank-$b$
        updates $C_{first\_proc(p,k)} := C_{first\_proc(p,k)} - W$ and 
        $C_{target\_first\_proc(p,k)} := C_{target\_first\_proc(p,k)} -
        Y_{level(p,k),k} \cdot W$ locally.  The latter computation
        is overlapped with the communication in Line \ref{step_comm2}.}\label{Alg:CAQR:j:overlap2}
  \EndFor
\EndFor
\end{algorithmic}
\end{algorithm}

We see that CAQR consists of $\frac{n}{b}$ TSQR factorizations
involving $P_r$ processors each, and $n/b - 1$ applications of the
resulting Householder vectors.  Table \ref{tbl:CAQR:par:model}
expresses the performance model over a rectangular $P_r \times P_c$
grid of processors.  A detailed derivation of the model is given in
Appendix \ref{S:CAQR-par-detailed}.  According to the table, the
number of arithmetic operations and words transferred is roughly the
same between parallel CAQR and ScaLAPACK's parallel QR factorization,
but the number of messages is a factor $b$ times lower for CAQR.

\begin{table}[h]
\small
\centering
\begin{tabular}{l | l}
            & Parallel CAQR \\ \hline
\# messages & $\frac{3n}{b} \log P_r + \frac{2n}{b} \log P_c$ \\ \hline
\# words    & $\left( 
                   \frac{n^2}{P_c} 
                   + \frac{bn}{2} 
               \right) \log P_r
               + \left( 
                   \frac{mn - n^2/2}{P_r} + 2n 
               \right) \log P_c$ \\ \hline
\# flops    & $\frac{2n^2(3m-n)}{3P} 
               + \frac{bn^2}{2P_c} 
               + \frac{3bn(2m - n)}{2P_r} 
               + \left( \frac{4 b^2 n}{3} 
                   + \frac{n^2 (3b+5)}{2 P_c} 
               \right) \log P_r
               - b^2 n$ \\ \hline
\# divisions & $\frac{mn - n^2/2}{P_r} 
                + \frac{bn}{2} \left( \log P_r - 1 \right)$
             \\ \hline \hline
             & ScaLAPACK's \texttt{PDGEQRF} \\ \hline
\# messages & $3n \log P_r + \frac{2n}{b} \log P_c$ \\ \hline
\# words    & $\left( 
                   \frac{n^2}{P_c} 
                   + bn 
               \right) \log P_r
               + \left( 
                   \frac{mn - n^2/2}{P_r} 
                   + \frac{bn}{2} 
               \right) \log P_c$ \\ \hline
\# flops    & $\frac{2n^2(3m-n)}{3P} 
               + \frac{bn^2}{2P_c} 
               + \frac{3bn(2m - n)}{2P_r}
               - \frac{b^2 n}{3 P_r}$
            \\ \hline
\# divisions & $\frac{mn - n^2/2}{P_r}$ \\
\end{tabular}
\caption{Performance models of parallel CAQR and ScaLAPACK's
  \lstinline!PDGEQRF! when factoring an $m \times n$ matrix,
  distributed in a 2-D block cyclic layout on a $P_r \times P_c$ 
  grid of processors with square $b \times b$ blocks.  All terms are
  counted along the critical path.  In this table exclusively, ``flops'' 
  only includes floating-point additions and multiplications, not 
  floating-point divisions.  Some lower-order terms are omitted.  We
  generally assume $m \geq n$.  Note that the number of flops,
  divisions, and words transferred all roughly match between the two
  algorithms, but the number of messages is about $b$ times lower 
  for CAQR.}
\label{tbl:CAQR:par:model}
\end{table}

The parallelization of the computation is represented by the number of
multiplies and adds and by the number of divides, in Table
\ref{tbl:CAQR:par:model}.  We discuss first the parallelization of
multiplies and adds.  The first term for CAQR represents mainly the
parallelization of the local Householder update corresponding to the
leaves of the TSQR tree (the matrix-matrix multiplication in line
\ref{Alg:CAQR:j:local-update} of Algorithm \ref{Alg:CAQR:j}), and
matches the first term for \lstinline!PDGEQRF!.  The second term for
CAQR corresponds to forming the $T_{p0}$ matrices for the local
Householder update in line \ref{Alg:CAQR:j:local-update} of the
algorithm, and also has a matching term for \lstinline!PDGEQRF!.  The
third term for CAQR represents the QR factorization of a panel of
width $b$ that corresponds to the leaves of the TSQR tree (part of
line \ref{Alg:CAQR:j:local-factor}) and part of the local rank-$b$
update (triangular matrix-matrix multiplication) in line
\ref{Alg:CAQR:j:local-update} of the algorithm, and also has a
matching term for \lstinline!PDGEQRF!.

The fourth term in the number of multiplies and adds for CAQR
represents the redundant computation introduced by the TSQR
formulation. In this term, the number of flops performed for computing
the QR factorization of two upper triangular matrices at each node of
the TSQR tree is $(2/3) nb^2 \log(P_r)$.  The number of flops
performed during the Householder updates issued by each QR
factorization of two upper triangular matrices is $n^2 (3b+5)/(2 P_c)
\log(P_r)$.

The runtime estimation in Table \ref{tbl:CAQR:par:model} does not take
into account the overlap of computation and communication in lines
\ref{Alg:CAQR:j:overlap1} and \ref{step_comm1} of Algorithm
\ref{Alg:CAQR:j} or the overlap in steps \ref{step_comm2} and
\ref{Alg:CAQR:j:overlap2} of the algorithm.  Suppose that at each step
of the QR factorization, the condition
\[
     \alpha + \beta \frac{b (n_j - b)}{P_c} 
     > 
     \gamma b (b + 1) \frac{n_j - b}{P_c}
\]
is fulfilled.  This is the case for example when $\beta / \gamma >
b+1$.  Then the fourth non-division flops term that accounts for the
redundant computation is decreased by $n^2 (b+1) \log(P_r) / P_c$,
about a factor of $3$.

\begin{table}[h!]
\centering
\begin{tabular}{l | l}
             & Parallel CAQR w/ optimal $b$, $P_r$, $P_c$ \\ \hline
\# flops     & $\frac{2mn^2}{P} - \frac{2n^3}{3P}$ \\
\# messages  & $\frac{1}{4C}
                \sqrt{\frac{n P}{m}}
                \log^2\left( 
                    \frac{m P}{n} 
                \right) 
                \cdot \log\left( 
                    P \sqrt{\frac{m P}{n}}
                \right)$ \\
\# words & $\sqrt{\frac{m n^3}{P}} \log P
            - \frac{1}{4} \sqrt{\frac{n^5}{m P}} 
              \log\left( \frac{n P}{m} \right)$ \\
Optimal $b$    & $C \sqrt{\frac{m n}{P}} 
                  \log^{-2} \left( \frac{m P}{n} \right)$ \\
Optimal $P_r$  & $\sqrt{\frac{m P}{n}}$ \\
Optimal $P_c$  & $\sqrt{\frac{n P}{m}}$ 
               \\ \hline \hline
             & \texttt{PDGEQRF} w/ optimal $b$, $P_r$, $P_c$ \\ \hline
\# flops     & $\frac{2mn^2}{P} - \frac{2n^3}{3P}$ \\
\# messages & $\frac{n}{4C} 
                 \log\left( \frac{m P^5}{n} \right)
                 \log\left( \frac{m P}{n} \right)
               + \frac{3n}{2} \log\left( \frac{m P}{n} \right)$ \\
\# words & $\sqrt{\frac{m n^3}{P}} \log P
            - \frac{1}{4} \sqrt{\frac{n^5}{m P}} 
              \log\left( \frac{n P}{m} \right)$ \\
Optimal $b$    & $C \sqrt{\frac{mn}{P}} 
                  \log^{-1}\left( \frac{m P}{n} \right)$ \\
Optimal $P_r$  & $\sqrt{\frac{m P}{n}}$ \\
Optimal $P_c$  & $\sqrt{\frac{n P}{m}}$
               \\ \hline \hline
             & Theoretical lower bound \\ \hline
\# messages  & $\sqrt{\frac{n P}{2^{11} m}}$ \\ 
\# words     & $\sqrt{\frac{mn^3}{2^{11} P}}$ \\
\end{tabular}
\caption{Highest-order terms in the performance models of parallel
  CAQR, ScaLAPACK's \lstinline!PDGEQRF!, and theoretical lower bounds
  for each, when factoring an $m \times n$ matrix, distributed in a
  2-D block cyclic layout on a $P_r \times P_c$ grid of processors
  with square $b \times b$ blocks.  All terms are counted along the
  critical path.  The theoretical lower bounds assume that $n \geq
  2^{11} m / P$, i.e., that the matrix is not too tall and skinny.
  The parameter $C$ in both algorithms is a $\Theta(1)$ tuning 
  parameter.  In summary, if we choose $b$, $P_r$, and $P_c$ independently and
optimally for both algorithms, the two algorithms match in the number of flops and words
transferred, but CAQR sends a factor of $\Theta(\sqrt{mn/P})$ messages
fewer than ScaLAPACK QR.
This factor is the local memory requirement on each processor, up to a
small constant.}
\label{tbl:CAQR:par:model:opt}
\end{table}


The execution time for a square matrix ($m=n$), on a square grid of
processors ($P_r = P_c = \sqrt{P}$) and with more lower order terms
ignored, simplifies to:
\begin{equation}\label{Eq:CAQR:time_sq}
\begin{split}
T_{Par.\ CAQR}(n, n, \sqrt{P}, \sqrt{P}) = 
\gamma \left( \frac{4n^3}{3P} + \frac{3 n^2 b}{4 \sqrt{P}} \log P \right)  \\
 + \beta \frac{3n^2}{4\sqrt{P}} \log P
 + \alpha \frac{5n}{2b} \log P.
\end{split}
\end{equation}

\subsection{Choosing $b$, $P_r$, and $P_c$ to minimize
  runtime}\label{SS:CAQR:par:opt}

A simple 2-D block layout may not be optimal for all possible $m$,
$n$, and $P$.  In order to minimize the runtime of parallel CAQR, we
could use a general nonlinear optimization algorithm to minimize the
runtime model (see Table \ref{tbl:CAQR:par:model}) with respect to
$b$, $P_r$, and $P_c$.  However, some general heuristics can improve
insight into the roles of these parameters.  For example, we would
like to choose them such that the flop count is $2mn^2/P - (2/3) n^3/P$
plus lower-order terms, just as in ScaLAPACK's parallel QR
factorization.  We follow this heuristic in the following paragraphs.
In summary, if we choose $b$, $P_r$, and $P_c$ independently and
optimally for both parallel CAQR and ScaLAPACK's parallel Householder
QR, the two algorithms match in the number of flops and words
transferred, but CAQR sends a factor of $\Theta(\sqrt{mn/P})$ messages
fewer than ScaLAPACK QR (see Table \ref{tbl:CAQR:par:model:opt}).
This factor is the local memory requirement on each processor, up to a
small constant.

\subsubsection{Ansatz}\label{SSS:CAQR:par:opt:ansatz}

The parameters $b$, $P_r$ and $P_c$ are integers which must satisfy
the following conditions:
\begin{equation}\label{eq:CAQR:par:opt:ansatz:constraints}
\begin{aligned}
1 \leq P_r, P_c \leq P \\
P_r \cdot P_c = P      \\
1 \leq b \leq \frac{m}{P_r} \\
1 \leq b \leq \frac{n}{P_c} \\
\end{aligned}
\end{equation}
We assume in the above that $P_r$ evenly divides $m$ and that $P_c$
evenly divides $n$.  From now on in this section (Section
\ref{SSS:CAQR:par:opt:ansatz}), we implicitly allow $b$, $P_r$, and
$P_c$ to range over the reals.  The runtime models are sufficiently
smooth that for sufficiently large $m$ and $n$, we can round these
parameters back to integer values without moving far from the minimum
runtime.  Example values of $b$, $P_r$, and $P_c$ which satisfy the
constraints in Equation \eqref{eq:CAQR:par:opt:ansatz:constraints} are
\[
\begin{aligned}
P_r &= \sqrt{\frac{m P}{n}} \\
P_c &= \sqrt{\frac{n P}{m}} \\
b   &= \sqrt{\frac{m n}{P}} \\
\end{aligned}
\]
These values are chosen simultaneously to minimize the approximate
number of words sent, $n^2/P_c + mn/P_r - n^2/(2 P_r)$, and the approximate number
of messages, $5n/b$, where for simplicity we temporarily ignore
logarithmic factors and lower-order terms in Table
\ref{tbl:CAQR:par:model}.  This suggests using the following ansatz:
\begin{equation}\label{eq:CAQR:par:opt:ansatz}
\begin{aligned}
P_r &= K \cdot \sqrt{\frac{m P}{n}}, \\
P_c &= \frac{1}{K} \cdot \sqrt{\frac{n P}{m}}, \text{and} \\
b   &= B \cdot \sqrt{\frac{m n}{P}}, \\
\end{aligned}
\end{equation}
for general values of $K$ and $B \leq \min\{ K, 1/K \}$, since we can
thereby explore all possible values of $b$, $P_r$ and $P_c$ satisfying
\eqref{eq:CAQR:par:opt:ansatz:constraints}.

\subsubsection{Flops}\label{SSS:CAQR:par:opt:flops}

Using the substitutions in Equation \eqref{eq:CAQR:par:opt:ansatz},
the flop count (neglecting lower-order terms, including the division
counts) becomes
\begin{multline}\label{eq:CAQR:par:opt:ansatz:flops}
\frac{mn^2}{P} \left( 
    2 
    - B^2 
    + \frac{3B}{K} 
    + \frac{B K}{2} 
\right)
- \\ 
\frac{n^3}{P} \left( 
    \frac{2}{3} 
    + \frac{3B}{2K} 
\right)
+ \\
\frac{mn^2 \log\left( K \cdot \sqrt{\frac{mP}{n}} \right)}{P} \left( 
    \frac{4B^2}{3} 
    + \frac{3 B K}{2}
\right).
\end{multline}
We wish to choose $B$ and $K$ so as to minimize the flop count.  We
know at least that we need to eliminate the dominant $mn^2
\log(\dots)$ term, so that parallel CAQR has the same asymptotic flop
count as ScaLAPACK's \lstinline!PDGEQRF!.  This is because we know
that CAQR performs at least as many floating-point operations
(asymptotically) as \lstinline!PDGEQRF!, so matching the highest-order
terms will help minimize CAQR's flop count.

To make the high-order terms of \eqref{eq:CAQR:par:opt:ansatz:flops}
match the $2mn^2/P - 2n^3/(3P)$ flop count of ScaLAPACK's parallel QR
routine, while minimizing communication as well, we can pick $K=1$ and 
\[
B = o\left(
    \log^{-1}\left(
        \sqrt{\frac{ m P }{ n }}
    \right)
\right);
\]
we will use 
\begin{equation}\label{eq:CAQR:par:opt:flops:B}
B = C \log^{-2} \left(
    \sqrt{\frac{ m P }{ n }}
\right)
\end{equation}
for some positive constant $C$, for simplicity.

The above choices of $B$ and $K$ make the flop count as follows, with
some lower-order terms omitted:
\begin{equation}\label{eq:CAQR:par:opt:flops}
\frac{2mn^2}{P}
- \frac{2n^3}{3P}
+ \frac{3 C m n^2}{P \log\left( \frac{m P}{n} \right)}
\end{equation}
Thus, we can choose the block size $b$ so as to match the higher-order
terms of the flop count of ScaLAPACK's parallel QR factorization
\lstinline!PDGEQRF!.

\subsubsection{Number of messages}

Using the substitutions in Equations \eqref{eq:CAQR:par:opt:ansatz}
and \eqref{eq:CAQR:par:opt:flops:B} with $K = 1$, the number of
messages becomes
\begin{equation}
\label{eq:CAQR:par:opt:lat}
\frac{1}{C}
\sqrt{\frac{n P}{m}}
\cdot \log^2\left( \sqrt{\frac{m P}{n}} \right)
\cdot \log\left( P \sqrt{\frac{m P}{n}} \right).
\end{equation}
The best we can do with the latency is to make $C$ as large as
possible, which makes the block size $b$ as large as possible.  The
value $C$ must be a constant, however; specifically, the flop counts
require
\[
\begin{aligned}
C &= \Omega\left( \log^{-2} \left( K \sqrt{\frac{m P}{n}} \right)
\right)\,\text{and} \\
C &= O(1).
\end{aligned}
\]
We leave $C$ as a tuning parameter in the number of messages Equation
\eqref{eq:CAQR:par:opt:lat}.

\subsubsection{Communication volume}

Using the substitutions in Equation \eqref{eq:CAQR:par:opt:ansatz}
and \eqref{eq:CAQR:par:opt:flops:B}, the number of words transferred
between processors on the critical path, neglecting lower-order terms,
becomes
\begin{multline}\label{eq:CAQR:par:opt:bw}
\sqrt{\frac{m n^3}{P}} \log P 
- \frac{1}{4} \sqrt{\frac{n^5}{m P}} \log \left( \frac{n P}{m} \right)
+ \frac{C}{4} \sqrt{\frac{m n}{P}} \log^3\left( \frac{m P}{n} \right)
\approx \\
\sqrt{\frac{m n^3}{P}} \log P 
- \frac{1}{4} \sqrt{\frac{n^5}{m P}} \log \left( \frac{n P}{m} \right).
\end{multline}
In the second step above, we eliminated the $C$ term, as it is a
lower-order term (since $m \geq n$).  Thus, $C$ only has a significant
effect on the number of messages and not the number of words
transferred.

\subsubsection{Table of results}

Table \ref{tbl:CAQR:par:model:opt} shows the number of messages and
number of words used by parallel CAQR and ScaLAPACK when $P_r$, $P_c$,
and $b$ are independently chosen so as to minimize the runtime models,
as well as the optimal choices of these parameters.  In summary, if we
choose $b$, $P_r$, and $P_c$ independently and optimally for both
algorithms, the two algorithms match in the number of flops and words
transferred, but CAQR sends a factor of $\Theta(\sqrt{mn/P})$ messages
fewer than ScaLAPACK QR.  This factor is the local memory requirement
on each processor, up to a small constant.


\subsection{Look-ahead approach}

Our models assume that the QR factorization does not use a look-ahead
technique during the right-looking factorization.  With the look-ahead
right-looking approach, the communications are pipelined from left to
right.  At each step of factorization, we would model the latency cost
of the broadcast within rows of processors as $2$ instead of
$\log P_c$.

In the next section, we will describe the sequential CAQR algorithm.

\section{Sequential 2-D QR factorization}\label{S:CAQR-seq}

The sequential CAQR algorithm uses sequential TSQR to perform a QR
factorization of a dense matrix $A$.  It mostly follows the
out-of-DRAM QR factorization of Gunter et al.\
\cite{gunter2005parallel}.  The matrix is stored in a $P_r \times P_c$
2-D block layout, and the blocks are paged between slow and fast
memory and operated on in fast memory.  Sequential CAQR is based on
TSQR in order to minimize the number of messages as well as the number
of words transferred between slow and fast memory.  At each step of
the factorization, sequential TSQR is used to factor a panel of
columns.  Then, the resulting Householder vectors are applied to the
rest of the matrix.  In summary, Table \ref{tbl:CAQR:seq:model:opt}
shows that the number of arithmetic operations is about the same
between sequential CAQR and blocked sequential left-looking
Householder QR (see Appendix \ref{S:PFDGEQRF}, but CAQR transfers a
factor of $\Theta(m/\sqrt{W})$ fewer words and $\sqrt{W}$ fewer
messages between slow and fast memory than blocked left-looking
Householder QR, in which $W$ is the fast memory capacity.  

We note that the references
\cite{buttari2007class,buttari2007parallel,gunter2005parallel,kurzak2008qr,quintana-orti2008scheduling}
propose an algorithm called ``tiled QR,'' which is the same as our
sequential CAQR with square blocks.  However, they use it in parallel
on shared-memory platforms, especially single-socket multicore.  They
do this by exploiting the parallelism implicit in the directed acyclic
graph of tasks.  Often they use dynamic task scheduling, which we
could use but do not discuss in this paper.  Since the cost of
communication in the single-socket multicore regime is low, these
authors are less concerned than we are about minimizing latency.  We
also model and analyze communication costs in more detail than
previous authors did.

\begin{table}[h]
\centering
\begin{tabular}{l | l}
             & Sequential CAQR \\ \hline
\# messages  & $12 \frac{mn^2}{W^{3/2}}$ \\
\# words     & $3 \frac{mn^2}{\sqrt{W}}$ \\
Opt.\ $P$    & $4mn/W$ \\
Opt.\ $P_r$  & $2m / \sqrt{W}$ \\
Opt.\ $P_c$  & $2n / \sqrt{W}$
               \\ \hline \hline
             & ScaLAPACK's \texttt{PFDGEQRF} \\ \hline
\# messages  & $\frac{mn^2}{2W} + \frac{2mn}{W}$ \\
\# words     & $\frac{m^2 n^2}{2W} - \frac{m n^3}{6W}
                + \frac{3mn}{2} - \frac{3n^2}{4}$ \\
Opt.\ $b$    & $1$ \\
Opt.\ $c$    & $\approx \frac{W}{m}$
               \\ \hline \hline
             & Lower bound \\ \hline
\# messages  & $\frac{\frac{mn^2}{4} - \frac{n^2}{8} \left( 
                   \frac{n}{2} + 1 \right)}{\sqrt{8W^3}}$ \\
\# words     & $\frac{\frac{mn^2}{4} - \frac{n^2}{8} \left( 
                   \frac{n}{2} + 1 \right)}{\sqrt{8W}}$ \\

\end{tabular}
\caption{Highest-order terms in the performance models of sequential
  CAQR, ScaLAPACK's out-of-DRAM QR factorization \texttt{PFDGEQRF} 
  running on one processor, and theoretical lower bounds for each, when
  factoring an $m \times n$ matrix with a fast memory capacity of $W$
  words.  In the case of sequential CAQR, the matrix is arranged in a
  2-D block layout on a $P_r \times P_c$ grid of $P$ blocks.  (See 
  Appendix \ref{SS:CAQR-seq-detailed:opt}).  The optimal choices of
  these parameters result in square blocks (i.e., $m/P_r = n/P_c$).  
  In the case of \texttt{PFDGEQRF}, the $b$ parameter is the left 
  panel width and the $c$ parameter is the current panel width.  (See
  Appendix \ref{S:PFDGEQRF}.)}
\label{tbl:CAQR:seq:model:opt}
\end{table}

Let $P = P_r \cdot P_c$ be the number of blocks.  We assume without
loss of generality that $P_r$ evenly divides $m$ and $P_c$ evenly
divides $n$.  The dimensions of a single block of the matrix are
$m/P_r \times n/P_c$.  We assume that $m \geq n$ and that $m/P_r \geq
n/P_c$.  We also assume that fast memory has a capacity of $W$
floating-point words for direct use by the algorithms in this section,
neglecting lower-order amounts of additional work space.

Algorithm \ref{Alg:CAQR:seq:RL} is the sequential CAQR factorization
which we analyze in this work.  It is a right-looking algorithm.
Appendix \ref{S:CAQR-seq-detailed} explains that the left-looking and
right-looking variants perform essentially the same number of
floating-point operations, and send essentially the same number of
words in the same number of messages, so we need only analyze the
right-looking version.

\begin{algorithm}[h]
\caption{Right-looking sequential CAQR factorization}
\label{Alg:CAQR:seq:RL}
\begin{algorithmic}[1]
\State{Assume: $m \geq n$ and $\frac{m}{P_r} \geq \frac{n}{P_c}$}
\For{$J = 1$ to $P_c$}
  \State{Factor panel $J$ (in rows $(J-1)\frac{n}{P_c} + 1$ to $m$
         and columns $(J-1)\frac{n}{P_c}+1$ to $J\frac{n}{P_c}$)}
  \State{Update trailing panels to right, 
         (in rows $(J-1)\frac{n}{P_c} +1$ to $m$
         and columns $J\frac{n}{P_c}+1$ to $n$)
         using the current panel}
\EndFor
\end{algorithmic}
\end{algorithm}

In Appendix \ref{S:CAQR-seq-detailed}, we derive the following model
(Equation \eqref{eq:CAQR:seq:modeltime:P} in the Appendix) of the
runtime of sequential CAQR:
\begin{multline*}
T_{\text{Seq.\ CAQR}} (m,n,P_c,P_r) \approx
\alpha \left[ \frac{3}{2} P(P_c-1) \right] + \\
\beta \left[ \frac{3}{2} mn \left( P_c + \frac{4}{3} \right) 
             - \frac{1}{2} n^2 P_c + O\left(n^2 + nP\right) \right] + \\
\gamma \left[ 2n^2m - \frac{2}{3}n^3 \right].
\end{multline*}
In Appendix \ref{SS:CAQR-seq-detailed:opt}, we show that the choices
of $P$, $P_r$, and $P_c$ that minimize the runtime are given by $P =
4mn/W$, $P_r = 2m / \sqrt{W}$, and $P_c = 2n / \sqrt{W}$.  These
values yield a runtime of (Equation \eqref{eq:CAQR:seq:modeltime:P:opt} in Appendix
\ref{S:CAQR-seq-detailed})
\begin{multline*} 
T_{\text{Seq.\ CAQR}} (m,n,W) \approx
    \alpha \left[ 12 \frac{mn^2}{W^{3/2}} \right] +
    \beta \left[ 3 \frac{mn^2}{\sqrt{W}} + 
    O\left( \frac{mn^2}{W} \right) \right] + \\
    \gamma \left[ 2mn^2 - \frac{2}{3}n^3 \right].
\end{multline*}
We note that the bandwidth term is proportional to
$\frac{mn^2}{\sqrt{W}}$, and the latency term is $W$ times smaller,
both of which match (to within constant factors), the lower bounds on
bandwidth and latency described in Corollary~\ref{corollary:SeqCAQR}
in Section~\ref{sec:SeqCAQR}.  Furthermore, the flop count $2mn^2 -
2n^3/3$ is identical to the flop count of LAPACK's sequential QR
factorization routine \lstinline!DGEQRF!.

\subsection{Other Bandwidth Minimizing Sequential QR Algorithms}
\label{sec:seq_qr_other}

In this section we describe special cases in which previous
sequential QR algorithms also minimize bandwidth, although
they do not minimize latency.
In particular, we discuss 
two variants of Elmroth's and Gustavson's recursive
QR (RGEQR3 and RGEQRF \cite{elmroth2000applying}),
as well as LAPACK's DGEQRF.

The fully recursive routine RGEQR3 is analogous to Toledo's
fully recursive LU routine \cite{toledo1997locality}: Both
routines factor the left half of the matrix (recursively), 
use the resulting factorization of the left half to update 
the right half, and then factor the right half (recursively again).
The base case consists of a single column. The output of
RGEQR3 applied to an $m$-by-$n$ matrix returns the $Q$ 
factor in the form $I-YTY^T$, where $Y$ is the $m$-by-$n$
lower triangular matrix of Householder vectors,
and $T$ is an $n$-by-$n$ upper triangular matrix.
A simple recurrence for the number of memory references
of either RGEQR3 or Toledo's algorithm is
\begin{eqnarray}
\label{eqn:RGEQR3}
B(m,n) & = & \left\{ \begin{array}{ll} 
             B(m,\frac{n}{2}) + B(m-\frac{n}{2},\frac{n}{2}) + 
             O(\frac{mn^2}{\sqrt{W}})
             & {\rm if} \; mn > W \; {\rm and} \; n>1 \\
             mn & {\rm if} \; mn \leq W  \\
             m  & {\rm if} \; m > W \; {\rm and} \; n=1 
             \end{array} \right. \nonumber \\
    & \leq & \left\{ \begin{array}{ll} 
             2B(m,\frac{n}{2}) + 
             O(\frac{mn^2}{\sqrt{W}})
             & {\rm if} \; mn > W \; {\rm and} \; n>1 \\
             mn & {\rm if} \; mn \leq W  \\
             m  & {\rm if} \; m > W \; {\rm and} \; n=1 
             \end{array} \right. \nonumber \\
       & = & O(\frac{mn^2}{\sqrt{W}}) + mn
\end{eqnarray}
So RGEQR3 attains our bandwidth lower bound.
(The $mn$ term must be included to account for the case
when $n<\sqrt{W}$, since each of the $mn$ matrix entries
must be accessed at least once.)
However, RGEQR3 does 
a factor greater than one
times as many floating point operations
as sequential Householder QR.

Now we consider RGEQRF and DGEQRF, which are both
right-looking algorithms and differ only in how
they perform the panel factorization (by RGEQR3
and DGEQR2, resp.). Let $b$ be the width of
the panel in either algorithm. It is easy to
see that a reasonable estimate of the number of 
memory references just for the updates by all the panels
is the number of panels $\frac{n}{b}$ times the minimum
number of memory references for the average
size update $\Theta(\max(mn,\frac{mnb}{\sqrt{W}}))$,
or $\Theta(\max(\frac{mn^2}{b},\frac{mn^2}{\sqrt{W}}))$.
Thus we need to pick $b$ at least about as large
as $\sqrt{W}$ to attain the desired lower bound
$O(\frac{mn^2}{\sqrt{W}})$.

Concentrating now on RGEQRF, we get from 
inequality~(\ref{eqn:RGEQR3})
that the $\frac{n}{b}$ panel factorizations using RGEQR3
cost at most an additional \linebreak
$O(\frac{n}{b} \cdot [\frac{mb^2}{\sqrt{W}} + mb] )
= O( \frac{mnb}{\sqrt{W}} + mn)$
memory references, or $O(mn)$ if we pick $b=\sqrt{W}$.
Thus the total number of memory references for RGEQRF
with $b= \sqrt{W}$ is $O(\frac{mn^2}{\sqrt{W}} + mn)$
which attains the desired lower bound.

Next we consider LAPACK's DGEQRF.
In the worst case, a panel factorization by DGEQR2 will incur
one slow memory access per arithmetic operation,
and so $O(\frac{n}{b} \cdot mb^2 ) =  O(mnb)$ for all panel factorizations.
For the overall algorithm to be guaranteed to attain
minimal bandwidth, we need $mnb = O(\frac{mn^2}{\sqrt{W}})$,
or $b = O(\frac{n}{\sqrt{W}})$. Since $b$ must also be at 
least about $\sqrt{W}$, this means $W = O(n)$,
or that fast memory size may be at most large enough
to hold a few rows of the matrix, or may be much smaller.

RGEQR3 does not alway minimize latency. For example,
considering applying RGEQR3 to a single panel
with $n=\sqrt{W}$ columns and $m>W$ rows, stored
in a block-column layout with $\sqrt{W}$-by-$\sqrt{W}$
blocks stored columnwise, as above. Then a recurrence
for the number of messages RGEQR3 requires is
\begin{eqnarray*}
\label{eqn:RGEQR3_latency}
L(m,n) & = & \left\{ \begin{array}{ll} 
             L(m,\frac{n}{2}) + L(m-\frac{n}{2},\frac{n}{2}) + 
             O(\frac{m}{\sqrt{W}})
             & {\rm if} \; n>1 \\
             O(\frac{m}{\sqrt{W}}) & {\rm if} \; n = 1  
             \end{array} \right. \nonumber \\
       & = & O(\frac{mn}{\sqrt{W}}) = O(m) \; {\rm when} \; n = \sqrt{W} 
\end{eqnarray*}
which is larger than the minimum $O(\frac{mn}{W}) = O(\frac{m}{\sqrt{W}})$
attained by sequential TSQR when $n = \sqrt{W}$.

In contrast to DGEQRF, RGEQRF, and RGEQR3, 
CAQR minimizes flops, bandwidth and latency
for all values of $W$.


\section{Comparison of ScaLAPACK's parallel QR and CAQR}
\label{S:CAQR:counts}
\label{S:CAQR-counts}

Here, we compare ScaLAPACK's QR factorization routine
\lstinline!PDGEQRF! with parallel CAQR.  Table
\ref{tbl:CAQR:par:model:opt} summarizes the results of this
comparison: if we choose the $b$, $P_r$, and $P_c$ parameters
independently and optimally for both algorithms, the two algorithms
match in the number of flops and words transferred, but CAQR sends a
factor of $\Theta(\sqrt{mn/P})$ messages fewer than ScaLAPACK QR.
This factor is the local memory requirement on each processor, up to a
small constant.

\subsection{\lstinline!PDGEQRF! performance model}

We suppose that we decompose a $m \times n$ matrix with $m
\geq n$ which is distributed block cyclically over a $P_r$ by $P_c$
grid of processors, where $P_r \times P_c = P$.  The two-dimensional
block cyclic distribution uses square blocks of dimension $b \times
b$.  Equation \eqref{Eq:ScaLAPACK:time} represents the runtime
estimation of ScaLAPACK's QR, in which we assume that there is no
attempt to pipeline communications from left to right and some lower
order terms are omitted.
\begin{multline}
\label{Eq:ScaLAPACK:time}
T_{SC}( m, n, P_r, P_c ) = \\
\left[ 
    \frac{2 n^2}{3 P} \left( 3m - n \right)  
    + \frac{ 3(b+1) n \left(m - \frac{n}{2}\right) }{P_r}
    + \frac{b n^2}{2 P_c}  
    - b n \left( \frac{b}{3} + \frac{3}{2} \right)
\right] \gamma + \\
\left[ 
    \frac{ m n - \frac{n^2}{2} }{P_r} 
\right] \gamma_d + \\
\left[
    3 n \left( 
        1 + \frac{1}{b} 
    \right) \log P_r 
    + \frac{2n}{b} \log P_c 
\right] \alpha  + \\
\left[
    \left( 
        \frac{n^2}{P_c} + n (b+2) 
    \right) \log P_r
    + \left(
        \frac{1}{P_r} \left( mn - \frac{n^2}{2} \right) +
        \frac{nb}{2} 
    \right) \log P_c
\right] \beta \\
\end{multline}
Compare with a less detailed but similar performance estimation in
\cite{scalapackusersguide}, in particular Tables 5.1 and 5.8 (routine
\lstinline!PxGELS!, whose main cost is invoking \lstinline!PDGEQRF!)
and Equation (5.1).

When $P_r = P_c = \sqrt{P}$ and $m=n$, and ignoring more lower-order
terms, Equation \eqref{Eq:ScaLAPACK:time} simplifies to
\begin{equation}
\label{Eq:ScaLAPACK:time:square}
T_{SC}(n,n,\sqrt{P}, \sqrt{P}) = \gamma \frac{4}{3} \frac{n^3}{P} 
           + \beta \frac{3}{4} \log{P} \frac{n^2}{\sqrt{P}}
           + \alpha \left( \frac{3}{2} + \frac{5}{2b} \right) n \log P
\end{equation}

\subsection{Choosing $b$, $P_r$, and $P_c$ to minimize
  runtime}\label{SS:CAQR-counts:par:opt}
\label{SS:PDGEQRF:par:opt}

This paper, and this section in particular, aim to show that parallel
CAQR performs better than ScaLAPACK's parallel QR factorization
\lstinline!PDGEQRF!.  In order to make a fair comparison between the
two routines, we need to choose the parameters $b$, $P_r$, and $P_c$
so as to minimize the runtime of ScaLAPACK QR.  Even though
\lstinline!PDGEQRF! accepts input matrices in general 2-D block cyclic
layouts, users may prefer a 2-D block layout for simplicity.  However,
this may not be optimal for all possible $m$, $n$, and $P$.  In order
to minimize the runtime of ScaLAPACK's parallel QR factorization, we
could use a general nonlinear optimization algorithm to minimize the
runtime model (Equation \eqref{Eq:ScaLAPACK:time} in Section
\ref{S:CAQR-counts}) with respect to $b$, $P_r$, and $P_c$.  However,
some general heuristics can improve insight into the roles of these
parameters.  For example, we would like to choose them such that the
flop count is $(2mn^2 - 2n^3/3) / P$ plus lower-order terms, which in
terms of floating-point operations would offer the best possible
speedup for parallel Householder QR.  We follow this heuristic in the
following paragraphs.

\subsubsection{Ansatz}\label{SSS:CAQR-counts:par:opt:ansatz}
\label{SSS:PDGEQRF:par:opt:ansatz}

Just as with parallel CAQR (see
\eqref{eq:CAQR:par:opt:ansatz:constraints} in Section
\ref{SSS:CAQR:par:opt:ansatz}), the parameters $b$, $P_r$, and $P_c$
must satisfy the following conditions:
\begin{equation}\label{eq:PDGEQRF:par:opt:ansatz:constraints}
\begin{aligned}
1 \leq P_r, P_c \leq P \\
P_r \cdot P_c = P      \\
1 \leq b \leq \frac{m}{P_r} \\
1 \leq b \leq \frac{n}{P_c} \\
\end{aligned}
\end{equation}
As in Section \ref{SSS:CAQR:par:opt:ansatz}, we assume in the above
that $P_r$ evenly divides $m$ and that $P_c$ evenly divides $n$.  From
now on in this section (Section \ref{SSS:PDGEQRF:par:opt:ansatz}), we
implicitly allow $b$, $P_r$, and $P_c$ to range over the reals.  The
runtime models are sufficiently smooth that for sufficiently large $m$
and $n$, we can round these parameters back to integer values without
moving far from the minimum runtime.  Example values of $b$, $P_r$,
and $P_c$ which satisfy the constraints in Equation
\eqref{eq:PDGEQRF:par:opt:ansatz:constraints} are
\[
\begin{aligned}
P_r &= \sqrt{\frac{m P}{n}} \\
P_c &= \sqrt{\frac{n P}{m}} \\
b   &= \sqrt{\frac{m n}{P}} \\
\end{aligned}
\]
These values are chosen simultaneously to minimize the approximate
number of words sent, $n^2/P_c + mn/P_r$, and the approximate number
of messages, $3n + 5n/b$, where for simplicity we temporarily ignore
logarithmic factors and lower-order terms in Table
\ref{tbl:CAQR:par:model}.  This suggests using the following ansatz:
\begin{equation}\label{eq:PDGEQRF:par:opt:ansatz}
\begin{aligned}
P_r &= K \cdot \sqrt{\frac{m P}{n}}, \\
P_c &= \frac{1}{K} \cdot \sqrt{\frac{n P}{m}}, \text{and} \\
b   &= B \cdot \sqrt{\frac{m n}{P}}, \\
\end{aligned}
\end{equation}
for general values of $K$ and $B \leq \min\{ K, 1/K \}$, since we can
thereby explore all possible values of $b$, $P_c$ and $P_c$ satisfying
\eqref{eq:PDGEQRF:par:opt:ansatz:constraints}.  For simplicity, this
is the same ansatz as in Section \ref{SSS:CAQR:par:opt:ansatz}.

\subsubsection{Flops}\label{SSS:CAQR-counts:par:opt:flops}

Using the substitutions in Equation \eqref{eq:PDGEQRF:par:opt:ansatz},
the flop count (neglecting lower-order terms, including the division
counts) becomes
\begin{equation}\label{eq:PDGEQRF:par:opt:ansatz:flops}
\frac{mn^2}{P} \left(
    2 
    - \frac{B^2}{3} 
    + \frac{3B}{K} 
    + \frac{BK}{2}
\right)
-
\frac{n^3}{P} \left(
    \frac{2}{3} 
    + \frac{3B}{2K}
\right).
\end{equation}
As in Section \ref{SSS:CAQR:par:opt:flops}, we wish to choose $B$ and
$K$ so as to minimize the flop count.  We can do this by making the
flop count $(2mn^2 - 2n^3/3) / P$, because that is the best possible
parallel flop count for a parallelization of standard Householder QR.
To make the high-order terms of Equation
\eqref{eq:PDGEQRF:par:opt:ansatz:flops} $(2mn^2 - 2n^3/3) / P$, while
minimizing communication as well, we can pick $K=1$ and
\begin{equation}\label{eq:PDGEQRF:par:opt:flops:B}
B = 
\tilde{C} \log^{-c}\left( \sqrt{\frac{m P}{n}} \right) = 
(\tilde{C} \cdot 2^c) \log^{-c}\left( \frac{m P}{n} \right) = 
C \log^{-c}\left( \frac{m P}{n} \right)
\end{equation}
for some positive constant $C$ ($C = 2^c \tilde{C}$ for some positive
constant $\tilde{C}$) and positive integer $c \geq 1$.  Unlike in
Section \ref{SSS:CAQR:par:opt:flops}, a $log^{-1}(\dots)$ term
suffices to make the flop count $(2mn^2 - 2n^3/3) / P$ plus
lower-order terms.  This is because the parallel CAQR flop count
(Equation \eqref{eq:CAQR:par:opt:ansatz:flops} in Section
\ref{SSS:CAQR:par:opt:flops}) involves an additional $(4B^2/3 + 3BK/2)
mn^2 \log\left( \dots \right)$ term which must be made $O(mn^2
\log^{-1}(\dots))$.  We will choose $c$ below in order to minimize
communication.

\subsubsection{Number of messages}\label{SSS:CAQR-counts:par:opt:lat}

If we use the substitutions in Equations
\eqref{eq:PDGEQRF:par:opt:ansatz} and
\eqref{eq:PDGEQRF:par:opt:flops:B}, the number of messages becomes
\begin{multline}\label{eq:PDGEQRF:par:opt:lat:cC}
\text{Messages}_{\text{\texttt{PDGEQRF}}}(m, n, P, c, C) = \\
\frac{3n}{2} \log\left( \frac{m P}{n} \right)
+ 
\frac{n}{C} 
\left(
    2 \log P
    + 
    \frac{1}{2} \log\left( \frac{m P}{n} \right)
\right)
\left(
    \frac{ \log\left( \frac{m P}{n} \right) }{ 2 }
\right)^c
\end{multline}
In Section \ref{SSS:CAQR-counts:par:opt:flops}, we argued that the
parameter $c$ must satisfy $c \geq 1$.  As long as
\[
\frac{\log\left( \frac{m P}{n} \right)}{2} > 1
\]
is satisfied, it is clear from Equation
\eqref{eq:PDGEQRF:par:opt:lat:cC} that choosing $c = 1$ minimizes the
number of messages.  This results in a number of messages of
\begin{multline}\label{eq:PDGEQRF:par:opt:lat}
\text{Messages}_{\text{\texttt{PDGEQRF}}}(m, n, P, C) = \\
\frac{3n}{2} \log\left( \frac{m P}{n} \right)
+ 
\frac{n}{2C} 
\left(
    2 \log P
    + 
    \frac{1}{2} \log\left( \frac{m P}{n} \right)
\right)
\left(
    \log\left( \frac{m P}{n} \right)
\right) = \\
\frac{3n}{2} \log\left( \frac{m P}{n} \right)
+ 
\frac{n}{C} \log P \log\left( \frac{m P}{n} \right)
+
\frac{n}{4C} \left( \log\left( \frac{m P}{n} \right) \right)^2.
\end{multline}
The parameter $C$ must be $o(\log( m P / n ))$ in order to minimize
the number of messages (see Equation
\eqref{eq:PDGEQRF:par:opt:lat:cC}).  This means that the third term in
the last line of the above equation is dominant, as we assume $m \geq
n$.  Making $C$ larger thus reduces the number of messages.  However,
in practice, a sufficiently large $C$ may make the first term ($1.5 n
\log(m P / n)$) significant.  Thus, we leave $C$ as a tuning
parameter.

\subsubsection{Communication volume}\label{SSS:CAQR-counts:par:opt:bw}

Using the substitutions in Equation \eqref{eq:PDGEQRF:par:opt:ansatz}
and \eqref{eq:PDGEQRF:par:opt:flops:B}, the number of words transferred
between processors on the critical path, neglecting lower-order terms,
becomes
\begin{multline}\label{eq:PDGEQRF:par:opt:bw:cC}
\text{Words}_{\text{\texttt{PDGEQRF}}}(m, n, P, c, C) = \\
\sqrt{\frac{m n^3}{P}} \log P
- \frac{1}{4} \sqrt{\frac{n^5}{m P}} \log\left( \frac{n P}{m} \right)
+ \frac{C \log^{-c}\left( \frac{m P}{n} \right)}{4} 
  \sqrt{\frac{m n^3}{P}} 
  \log\left( \frac{n P}{m} \right).
\end{multline}
In Section \ref{SSS:CAQR-counts:par:opt:lat}, we argued for choosing
$c = 1$ in order to minimize the number of messages.  In that case,
the number of words transferred is
\begin{multline}\label{eq:PDGEQRF:par:opt:bw}
\text{Words}_{\text{\texttt{PDGEQRF}}}(m, n, P, C) = \\
\sqrt{\frac{m n^3}{P}} \log P
- \frac{1}{4} \sqrt{\frac{n^5}{m P}} \log\left( \frac{n P}{m} \right)
- \frac{C}{4} \sqrt{\frac{m n^3}{P}} \approx \\
\sqrt{\frac{m n^3}{P}} \log P
- \frac{1}{4} \sqrt{\frac{n^5}{m P}} \log\left( \frac{n P}{m} \right).
\end{multline}
The third term in the second line is a lower-order term (it is
subsumed by the first term), since $C = O(1)$.

\subsubsection{Table of results}

Table \ref{tbl:CAQR:par:model:opt} shows the number of messages and
number of words used by parallel CAQR and ScaLAPACK when $P_r$, $P_c$,
and $b$ are chosen so as to minimize the runtime model, as well as the
optimal choices of these parameters.  If we choose $b$, $P_r$, and
$P_c$ independently and optimally for both algorithms, the two
algorithms match in the number of flops and words transferred, but
CAQR sends a factor of $\Theta(\sqrt{mn/P})$ messages fewer than
ScaLAPACK QR.  This factor is the local memory requirement on each
processor, up to a small constant.

\section{Parallel CAQR performance estimation}\label{S:CAQR:perfest}

We use the performance model developed in the previous section to
estimate the performance of parallel CAQR on three computational
systems, IBM POWER5, Peta, and Grid, and compare it to ScaLAPACK's
parallel QR factorization routine \lstinline!PDGEQRF!.  Peta is a
model of a petascale machine with $8100$ processors, and Grid is a
model of $128$ machines connected over the Internet.  Each processor in
Peta and Grid can be itself a parallel machine, but our models
consider the parallelism only between these parallel machines.

We expect CAQR to outperform ScaLAPACK, in part because it uses a
faster algorithm for performing most of the computation of each panel
factorization (\lstinline!DGEQR3! vs.\ \lstinline!DGEQRF!), and in
part because it reduces the latency cost.  Our performance model uses
the same time per floating-point operation for both CAQR and
\lstinline!PDGEQRF!.  Hence our model evaluates the improvement due
only to reducing the latency cost.

We evaluate the performance using matrices of size $n \times n$,
distributed over a $P_r \times P_c$ grid of $P$ processors using a 2D
block cyclic distribution, with square blocks of size $b \times b$.
For each machine we estimate the best performance of CAQR and
\lstinline!PDGEQRF! for a given problem size $n$ and a given number of
processors $P$, by finding the optimal values for the block size $b$
and the shape of the grid $P_r \times P_c$ in the allowed ranges.  The
matrix size $n$ is varied in the range $10^3$, $10^{3.5}$, $10^4$,
$\dots$, $10^{7.5}$.  The block size $b$ is varied in the range $1$,
$5$, $10$, $\dots$, $50$, $60$, $\dots$, $\min(200, m/P_r, n/P_c)$.
The number of processors is varied from $1$ to the largest power of
$2$ smaller than $p_{max}$, in which $p_{max}$ is the maximum number
of processors available in the system.  The values for $P_r$ and $P_c$
are also chosen to be powers of two.

We describe now the parameters used for the three parallel machines.
The available memory on each processor is given in units of 8-byte
(IEEE 754 double-precision floating-point) words.  When we evaluate
the model, we set the $\gamma$ value in the model so that the modeled
floating-point rate is 80\% of the machine's peak rate, so as to
capture realistic performance on the local QR factorizations.  This
estimate favors ScaLAPACK rather than CAQR, as ScaLAPACK requires more
communication and CAQR more floating-point operations.  The inverse
network bandwidth $\beta$ has units of seconds per word.  The
bandwidth for Grid is estimated to be the Teragrid backbone bandwidth
of $40$ GB/sec divided by $p_{max}$.

\begin{itemize}
\item \textbf{IBM POWER5:} $p_{max} = 888$, peak flop rate is $7.6$
  Gflop/s, $mem = 5 \cdot 10^8$ words, $\alpha = 5 \cdot 10^{-6}$ s,
  $\beta = 2.5 \cdot 10^{-9}$ s/word ($1 / \beta = 400$ Mword/s $ =
  3.2$ GB/s).

\item \textbf{Peta:} $p_{max} = 8192$, peak flop rate is $500$
  Gflop/s, $mem = 62.5 \cdot 10^9$ words, $\alpha = 10^{-5}$ s, $\beta
  = 2 \cdot 10^{-9}$ s/word ($1 / \beta = 500$ Mword/s $= 4$ GB/s).

\item \textbf{Grid:} $p_{max} = 128$, peak flop rate is 10 Tflop/s,
  $mem = 10^{14}$ words, $\alpha = 10^{-1}$ s, $\beta = 25 \cdot
  10^{-9}$ s/word ($1 / \beta = 40$ Mword/s $= .32$ GB/s).
\end{itemize}

There are $13$ plots shown for each parallel machine.  The first three
plots display for specific $n$ and $P$ values our models of
\begin{itemize}
  \item the best speedup obtained by CAQR, with respect to the runtime
    using the fewest number of processors with enough memory to hold
    the matrix (which may be more than one processor),
  \item the best speedup obtained by \lstinline!PDGEQRF!, computed
    similarly, and
  \item the ratio of \lstinline!PDGEQRF! runtime to CAQR runtime.
\end{itemize}
The next ten plots are divided in two groups of five.  The first group
presents performance results for CAQR and the second group presents
performance results for \lstinline!PDGEQRF!.  The first two plots of
each group of five display the corresponding optimal values of $b$ and
$P_r$ obtained for each combination of $n$ and $P$.  (Since $P_c = P /
P_r$, we need not specify $P_c$ explicitly.)  The last $3$ plots of
each group of $5$ give the computation time to total time ratio, the
latency time to total time ratio, and the bandwidth time to total time
ratio.

The white regions in the plots signify that the problem needed too
much memory with respect to the memory available on the machine.  Note
that in our performance models, the block size $b$ has no meaning on
one processor, because there is no communication, and the term $4 n^3
/ (3 P)$ dominates the computation.  Thus, for one processor, we set
the optimal value of $b$ to 1 as a default.

CAQR leads to significant improvements with respect to
\lstinline!PDGEQRF! when the latency represents an important fraction
of the total time, as for example when a small matrix is computed on a
large number of processors.  On IBM POWER5, the best improvement is
predicted for the smallest matrix in our test set ($n = 10^3$), when
CAQR will outperform \lstinline!PDGEQRF! by a factor of $9.7$ on $512$
processors.  On Peta, the best improvement is a factor of $22.9$,
obtained for $n = 10^4$ and $P = 8192$.  On Grid, the best improvement
is obtained for one of the largest matrix in our test set
$m=n=10^{6.5}$, where CAQR outperforms \lstinline!PDGEQRF! by a factor
of $5.3$ on $128$ processors.

\subsection{Performance prediction on IBM POWER5}

\begin{figure}
  \begin{center}
    \mbox{
      \subfigure[Speedup CAQR]{\includegraphics[scale=0.35]{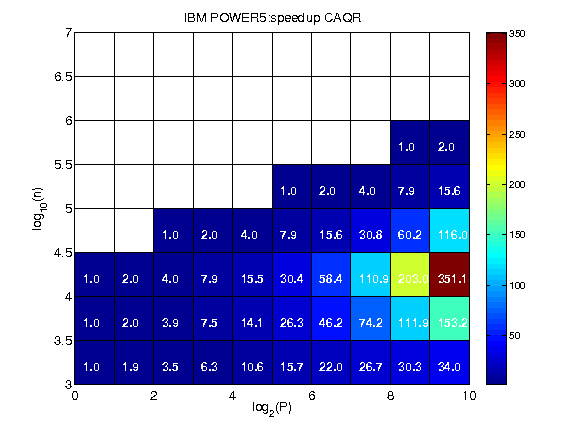}\label{ibmp5_spdCAQR}} 
      \subfigure[Speedup PDGEQRF]{\includegraphics[scale=0.35]{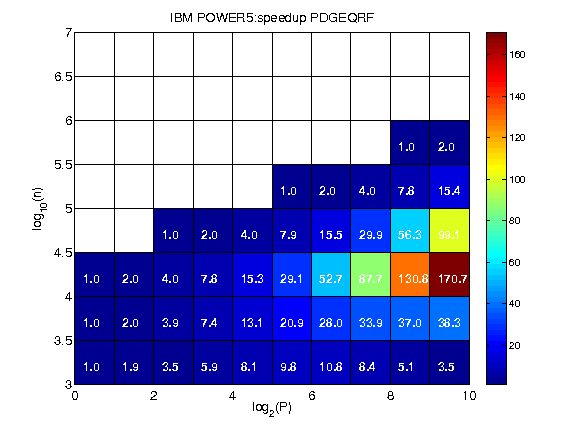}} 
    }
    \subfigure[Comparison]
	      {\includegraphics[scale=0.35]{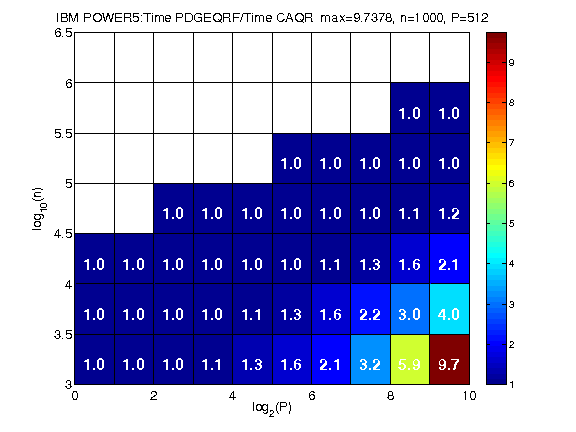}\label{ibmp5_CMP}}
  \end{center}
  \caption{\label{fig:PerfComp_ibmp5}Performance prediction comparing
    CAQR and \lstinline!PDGEQRF! on IBM POWER5.}
\end{figure}

\begin{figure}
  \begin{center}
    \mbox{
      \subfigure[Optimal $b$]{\includegraphics[scale=0.35]{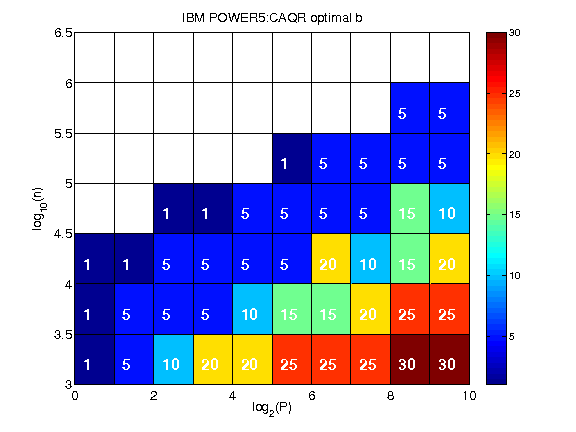}\label{ibmp5_CAQRb}}
      \subfigure[Optimal $P_r$]{\includegraphics[scale=0.35]{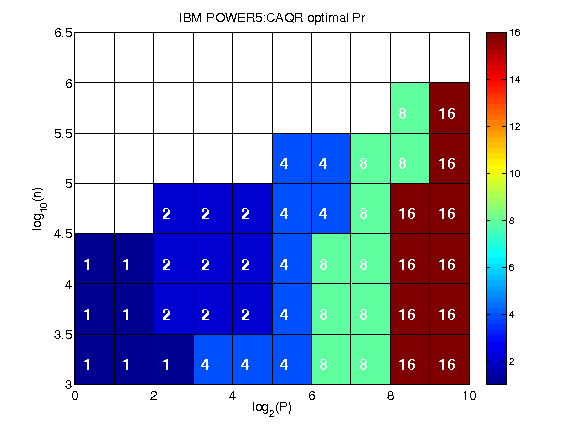}\label{ibmp5_CAQRPr}} 
    }
    \mbox{
      \subfigure[Fraction of time in computation]{\includegraphics[scale=0.35]{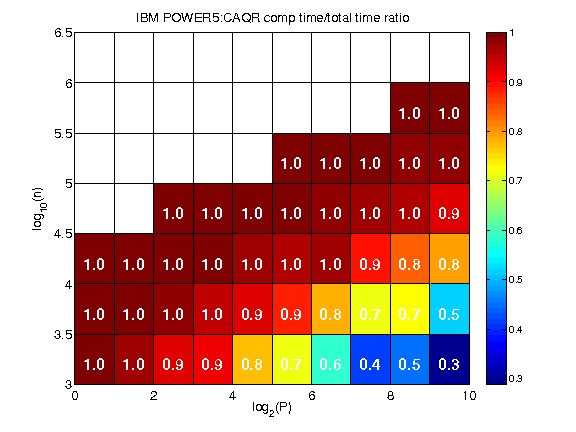}\label{ibmp5_CAQRcompR}}
      \subfigure[Fraction of time in latency]
		{\includegraphics[scale=0.35]{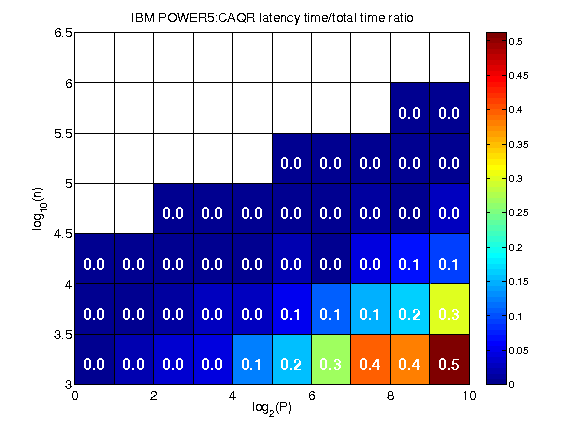}\label{ibmp5_CAQRlatR}}
    }
    \mbox{
      \subfigure[Fraction of time in bandwidth]{\includegraphics[scale=0.35]{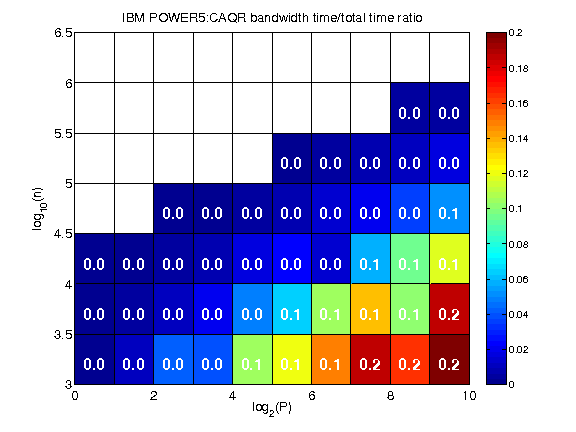}\label{ibmp5_CAQRbwR}}
    }
  \end{center}
  \caption{\label{fig:PerfCAQR_ibmp5}Performance prediction for
    CAQR on IBM POWER5.}
\end{figure}

\begin{figure}
  \begin{center}
    \mbox{
      \subfigure[Optimal $b$]{\includegraphics[scale=0.35]{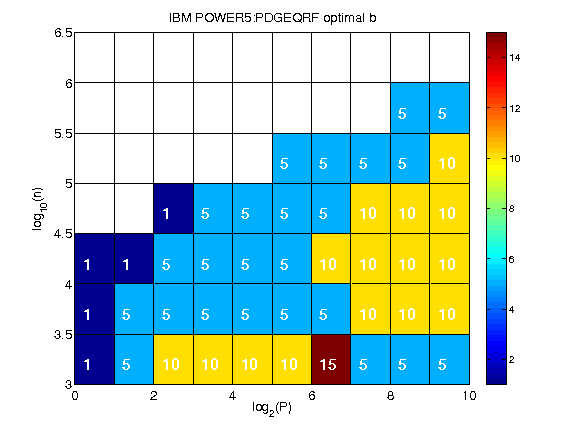}\label{ibmp5_QRFb}} 
      \subfigure[Optimal
	$P_r$]{\includegraphics[scale=0.35]{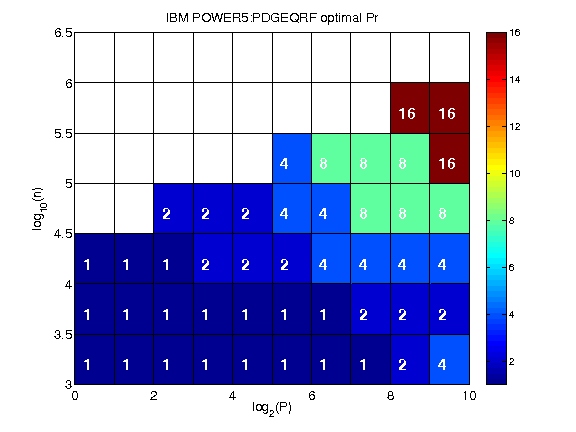}\label{ibmp5_QRFPr}} 
    }
    \mbox{
      \subfigure[Fraction of time in computation]{\includegraphics[scale=0.35]{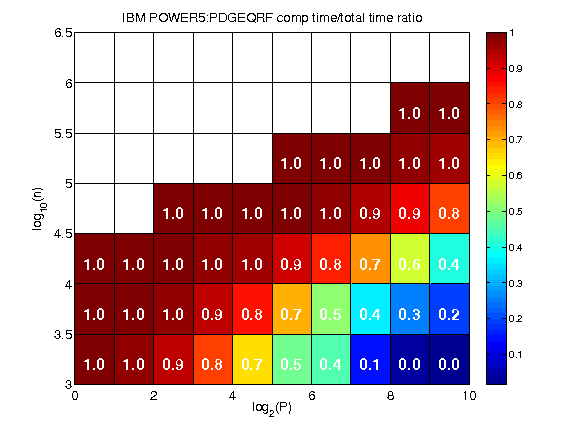}\label{ibmp5_QRFcompR}} 
      \subfigure[Fraction of time in latency]
		{\includegraphics[scale=0.35]{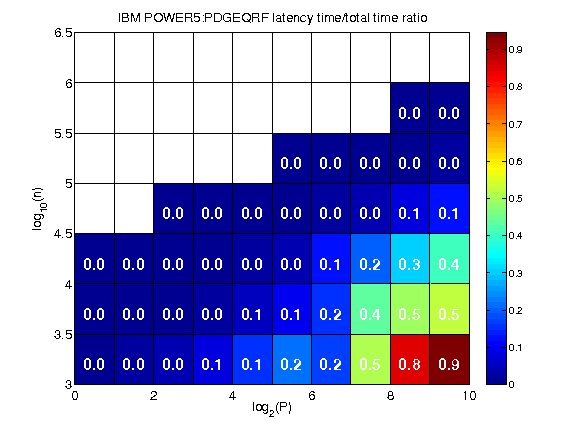}\label{ibmp5_QRFlatR}}
    }
    \mbox{
      \subfigure[Fraction of time in bandwidth]{\includegraphics[scale=0.35]{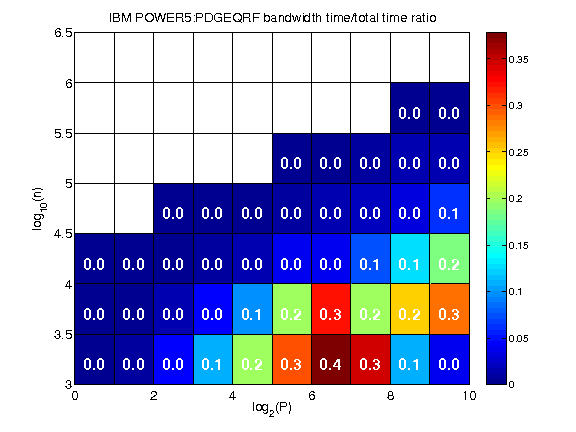}\label{ibmp5_QRFbwR}} 
    }
  \end{center}
  \caption{\label{fig:PerfPDGEQRF_ibmp5}Performance prediction for
    \lstinline!PDGEQRF! on IBM POWER5.}
\end{figure}

Figures \ref{fig:PerfComp_ibmp5}, \ref{fig:PerfCAQR_ibmp5}, and
\ref{fig:PerfPDGEQRF_ibmp5} depict modeled performance on the IBM
POWER 5 system.  CAQR has the same estimated performance as
\lstinline!PDGEQRF! when the computation dominates the total time.
But it outperforms \lstinline!PDGEQRF! when the fraction of time spent
in communication due to latency becomes significant.  The best
improvements are obtained for smaller $n$ and larger $P$, as displayed
in Figure \ref{ibmp5_CMP}, the bottom right corner.  For the smallest
matrix in our test set ($n = 10^3$), we predict that CAQR will
outperform \lstinline!PDGEQRF! by a factor of $9.7$ on $512$
processors.  As shown in Figure \ref{ibmp5_QRFlatR}, for this matrix,
the communication dominates the runtime of \lstinline!PDGEQRF!, with a
fraction of $0.9$ spent in latency.  For CAQR, the time spent in
latency is reduced to a fraction of $0.5$ of the total time from $0.9$
for PDGEQRF, and the time spent in computation is a fraction of $0.3$
of the total time.  This is illustrated in Figures
\ref{ibmp5_CAQRcompR} and \ref{ibmp5_CAQRlatR}.  

Another performance comparison consists in determining the improvement
obtained by taking the best performance independently for CAQR and
\lstinline!PDGEQRF!, when varying the number of processors from $1$ to
$512$.  For $n=10^3$, the best performance for CAQR is obtained when
using $P=512$ and the best performance for \lstinline!PDGEQRF! is
obtained for $P = 64$.  This leads to a speedup of more than $3$ for
CAQR compared to \lstinline!PDGEQRF!.  For any fixed $n$, we can take
the number of processors $P$ for which \lstinline!PDGEQRF! would
perform the best, and measure the speedup of CAQR over
\lstinline!PDGEQRF! using that number of processors.  We do this in
Table \ref{tbl:CAQR:par:POWER5:best}, which shows that CAQR always is
at least as fast as \lstinline!PDGEQRF!, and often significantly
faster (up to $3 \times$ faster in some cases).

Figure \ref{fig:PerfComp_ibmp5} shows that CAQR should scale well,
with a speedup of $351$ on $512$ processors when $m = n = 10^4$.  A
speedup of $116$ with respect to the parallel time on $4$ processors
(the fewest number of processors with enough memory to hold the
matrix) is predicted for $m=n=10^{4.5}$ on $512$ processors.  In these
cases, CAQR is estimated to outperform \lstinline!PDGEQRF! by factors
of $2.1$ and $1.2$, respectively.

Figures \ref{ibmp5_CAQRPr} and \ref{ibmp5_QRFPr} show that
\lstinline!PDGEQRF! has a smaller value for optimal $P_r$ than CAQR.
This trend is more significant in the bottom left corner of Figure
\ref{ibmp5_QRFPr}, where the optimal value of $P_r$ for
\lstinline!PDGEQRF! is $1$.  This corresponds to a 1D block column
cyclic layout.  In other words, \lstinline!PDGEQRF! runs faster by
reducing the $3 n \log{P_r}$ term of the latency cost of Equation
\eqref{Eq:ScaLAPACK:time} by choosing a small $P_r$.
\lstinline!PDGEQRF! also tends to have a better performance for a
smaller block size than CAQR, as displayed in Figures
\ref{ibmp5_CAQRb} and \ref{ibmp5_QRFb}.  The optimal block size $b$
varies from $1$ to $15$ for \lstinline!PDGEQRF!, and from $1$ to $30$
for CAQR.

\begin{table}
\begin{center}        
\begin{tabular}{r|c|c}
$\log_{10} n$ & Best $\log_2 P$ for \lstinline!PDGEQRF! & CAQR speedup
\\ \hline
3.0 & 6     & 2.1 \\
3.5 & 8     & 3.0 \\
4.0 & 9     & 2.1 \\
4.5 & 9     & 1.2 \\
5.0 & 9     & 1.0 \\
5.5 & 9     & 1.0 \\
\end{tabular}
\end{center}
\caption{Estimated runtime of \lstinline!PDGEQRF! divided by estimated
  runtime of CAQR on a square $n \times n$ matrix, on the IBM POWER5
  platform, for those values of $P$ (number of processors) for which
  \lstinline!PDGEQRF! performs the best for that problem size.}
\label{tbl:CAQR:par:POWER5:best}
\end{table}

\subsection{Performance prediction on Peta}

\begin{figure}
  \begin{center}
    \mbox{ \subfigure[Speedup
      CAQR]{\includegraphics[scale=0.35]{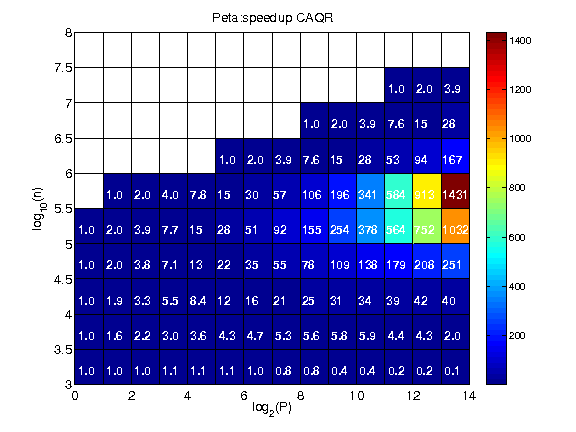}\label{peta_spdCAQR}}
      \subfigure[Speedup
      PDGEQRF]{\includegraphics[scale=0.35]{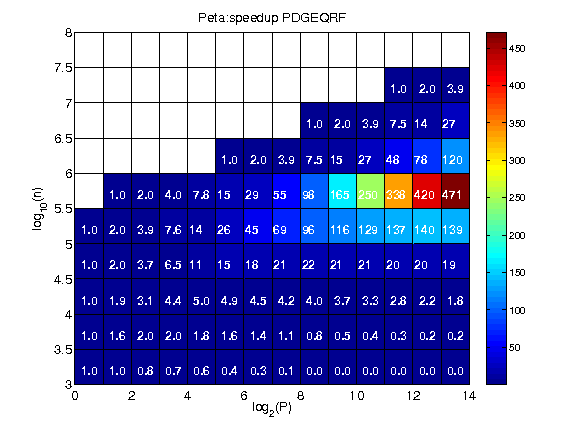}\label{peta_spdQRF}}
      } \subfigure[Comparison]
      {\includegraphics[scale=0.35]{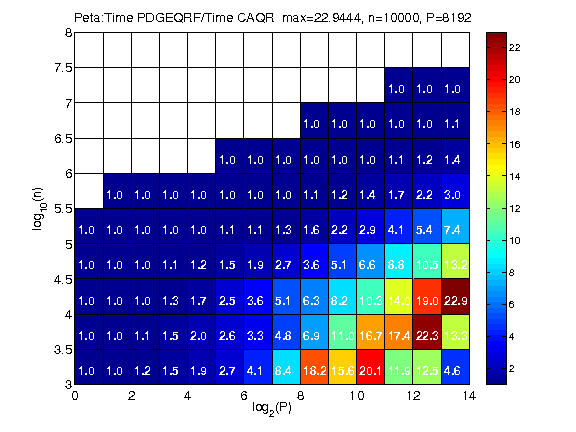}\label{peta_cmp}}
  \end{center}
  \caption{\label{fig:PerfComp_Peta}Performance prediction comparing
    CAQR and PDGEQRF on Peta.}
\end{figure}

\begin{figure}
  \begin{center}
    \mbox{ \subfigure[Optimal $b$]{\includegraphics[scale=0.35]{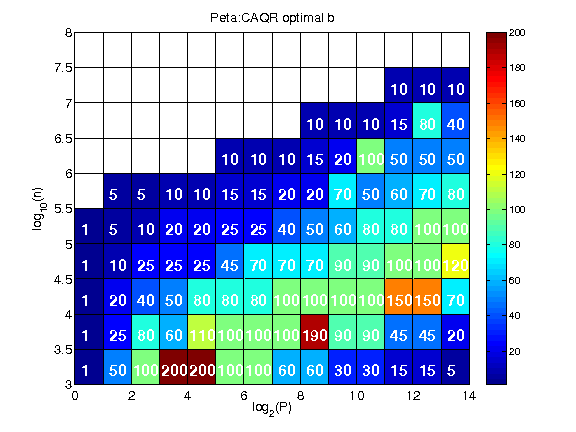}\label{peta_CAQRoptb}}
      \subfigure[Optimal $P_r$]{\includegraphics[scale=0.35]{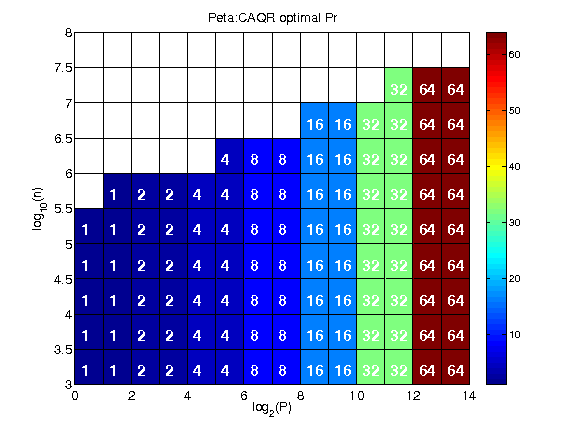}\label{peta_CAQRoptPr}}
    } 
    \mbox{ \subfigure[Fraction of time in
	computation]{\includegraphics[scale=0.35]{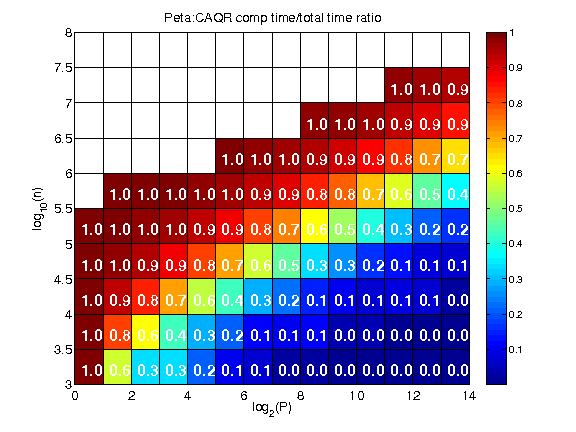}\label{peta_CAQRcompR}}
      \subfigure[Fraction of time in latency]
      {\includegraphics[scale=0.35]{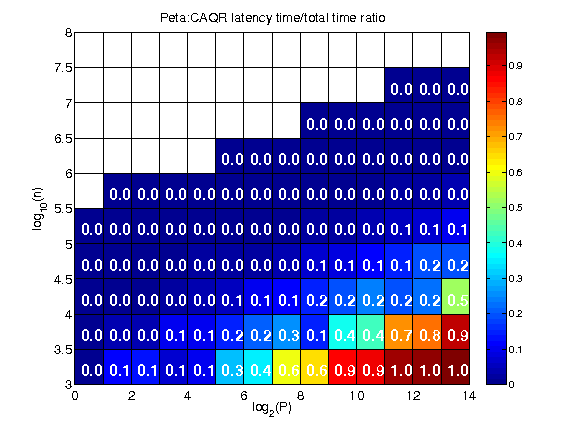}\label{peta_CAQRlatR}}
      }
    \mbox{ 
      \subfigure[Fraction of time in bandwidth]{\includegraphics[scale=0.35]{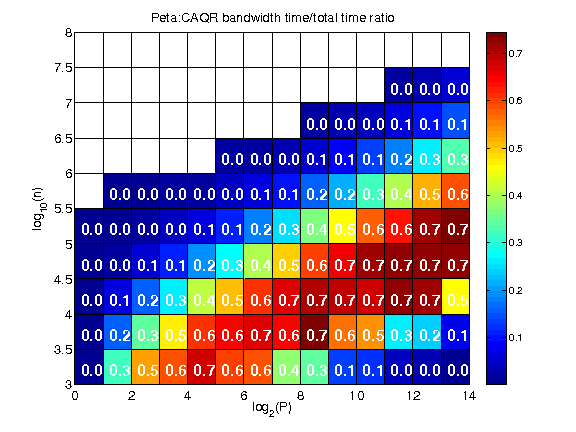}\label{peta_CAQRbwR}}
    }
  \end{center}
  \caption{\label{fig:PerfCAQR_Peta}Performance prediction for CAQR on
    Peta.}
\end{figure}

\begin{figure}
  \begin{center}
    \mbox{ \subfigure[Optimal $b$]{\includegraphics[scale=0.35]{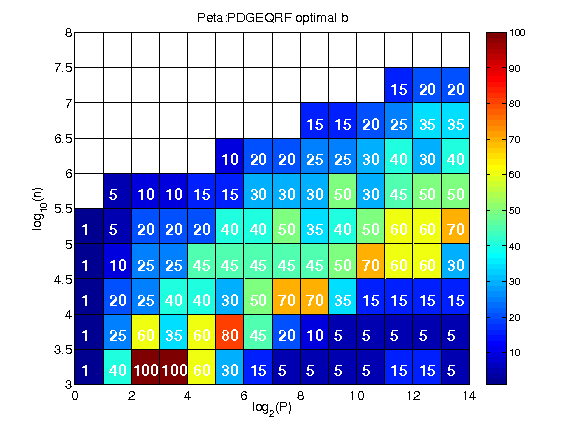}\label{peta_QRFoptb}}
      \subfigure[Optimal $P_r$]{\includegraphics[scale=0.35]{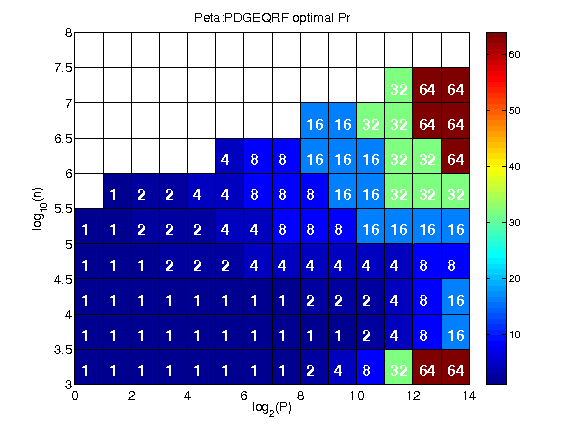}\label{peta_QRFoptPr}}
    } 
    \mbox{ 
      \subfigure[Fraction of time in computation]{\includegraphics[scale=0.35]{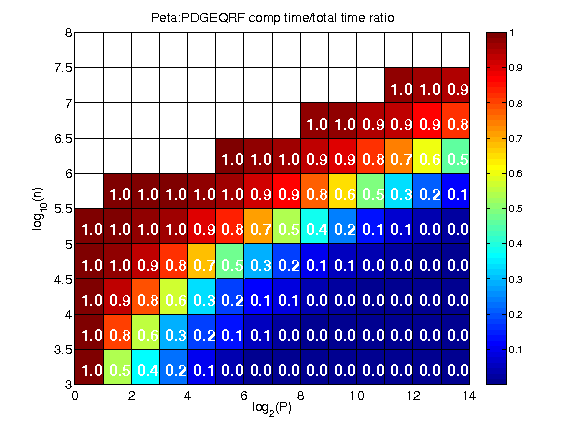}\label{peta_QRFcompR}}
      \subfigure[Fraction of time in latency]{\includegraphics[scale=0.35]{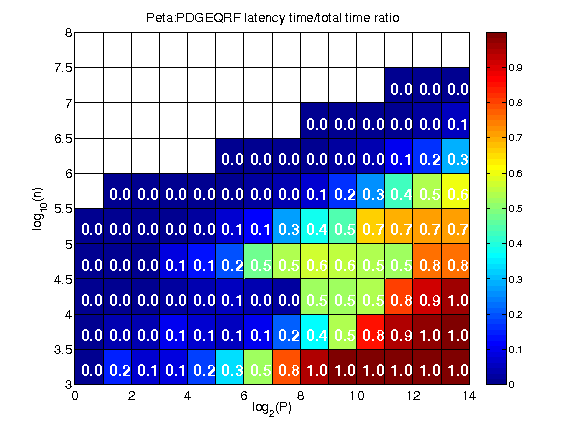}\label{peta_QRFlatR}}
    }
    \mbox{ 
      \subfigure[Fraction of time in bandwidth]{\includegraphics[scale=0.35]{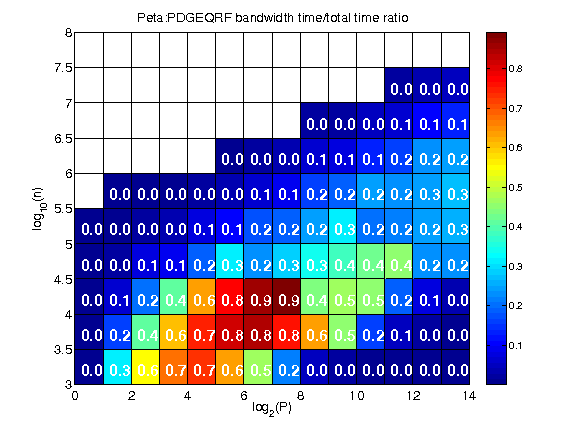}\label{peta_QRFbwR}}
    }
  \end{center}
  \caption{\label{fig:PerfPDGEQRF_Peta}Performance prediction for
    PDGEQRF on Peta.}
\end{figure}

Figures \ref{fig:PerfComp_Peta}, \ref{fig:PerfCAQR_Peta}, and
\ref{fig:PerfPDGEQRF_Peta} show our performance estimates of CAQR and
\lstinline!PDGEQRF! on the Petascale machine.  The estimated division
of time between computation, latency, and bandwidth for
\lstinline!PDGEQRF! is illustrated in Figures \ref{peta_QRFcompR},
\ref{peta_QRFlatR}, and \ref{peta_QRFbwR}.  In the upper left corner
of these figures, the computation dominates the total time, while in
the right bottom corner the latency dominates the total time.  In the
narrow band between these two regions, which goes from the left bottom
corner to the right upper corner, the bandwidth dominates the time.
CAQR decreases the latency cost, as can be seen in Figures
\ref{peta_CAQRcompR}, \ref{peta_CAQRlatR}, and \ref{peta_CAQRbwR}.
There are fewer test cases for which the latency dominates the time
(the right bottom corner of Figure \ref{peta_CAQRlatR}).  This shows
that CAQR is expected to be effective in decreasing the latency cost.
The left upper region where the computation dominates the time is
about the same for both algorithms.  Hence for CAQR there are more
test cases for which the bandwidth term is an important fraction of
the total time.

Note also in Figures \ref{peta_QRFoptPr} and \ref{peta_CAQRoptPr} that
optimal $P_r$ has smaller values for \lstinline!PDGEQRF! than for
CAQR.  There is an interesting regularity in the value of optimal
$P_r$ for CAQR.  CAQR is expected to have its best performance for
(almost) square grids.

As can be seen in Figure \ref{peta_spdCAQR}, CAQR is expected to show
good scalability for large matrices.  For example, for $n = 10^{5.5}$,
a speedup of $1431$, measured with respect to the time on $2$
processors, is obtained on $8192$ processors. For $n=10^{6.4}$ a
speedup of $166$, measured with respect to the time on $32$
processors, is obtained on $8192$ processors.

CAQR leads to more significant improvements when the latency
represents an important fraction of the total time.  This corresponds
to the right bottom corner of Figure \ref{peta_cmp}.  The best
improvement is a factor of $22.9$, obtained for $n = 10^4$ and $P =
8192$.  The speedup of the best CAQR compared to the best
\lstinline!PDGEQRF! for $n=10^4$ when using at most $P=8192$
processors is larger than $8$, which is still an important
improvement.  The best performance of CAQR is obtained for $P=4096$
processors and the best performance of \lstinline!PDGEQRF! is obtained
for $P=16$ processors.

Useful improvements are also obtained for larger matrices.  For $n =
10^6$, CAQR outperforms \lstinline!PDGEQRF! by a factor of $1.4$.
When the computation dominates the parallel time, there is no benefit
from using CAQR.  However, CAQR is never slower.  For any fixed $n$,
we can take the number of processors $P$ for which \lstinline!PDGEQRF!
would perform the best, and measure the speedup of CAQR over
\lstinline!PDGEQRF! using that number of processors.  We do this in
Table \ref{tbl:CAQR:par:Peta:best}, which shows that CAQR always is at
least as fast as \lstinline!PDGEQRF!, and often significantly faster
(up to $7.4 \times$ faster in some cases).

\begin{table}
\begin{tabular}{r|c|c}
$\log_{10} n$ & Best $\log_2 P$ for \lstinline!PDGEQRF! & CAQR speedup
\\ \hline
3.0 & 1     & 1 \\
3.5 & 2--3  & 1.1--1.5 \\
4.0 & 4--5  & 1.7--2.5 \\
4.5 & 7--10 & 2.7--6.6 \\
5.0 & 11--13 & 4.1--7.4 \\
5.5 & 13     & 3.0       \\
6.0 & 13     & 1.4       \\
\end{tabular}
\caption{Estimated runtime of \lstinline!PDGEQRF! divided by estimated
  runtime of CAQR on a square $n \times n$ matrix, on the Peta platform,
  for those values of $P$ (number of processors) for which
  \lstinline!PDGEQRF! performs the best for that problem size.}
\label{tbl:CAQR:par:Peta:best}
\end{table}

\subsection{Performance prediction on Grid}

\begin{figure}
\begin{center}
  \mbox{
    \subfigure[Speedup CAQR]{\includegraphics[scale=0.35]{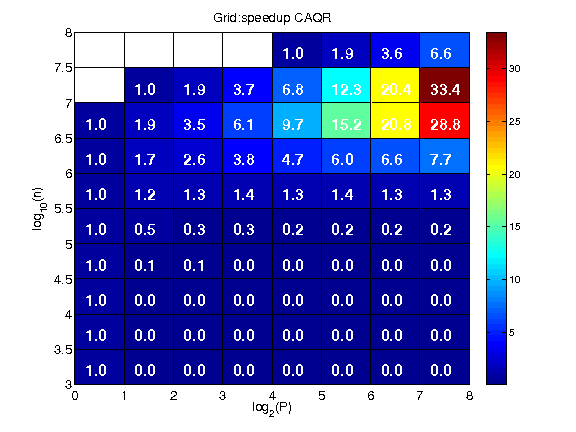}\label{grid_CAQRspdup}} 
    \subfigure[Speedup PDGEQRF]{\includegraphics[scale=0.35]{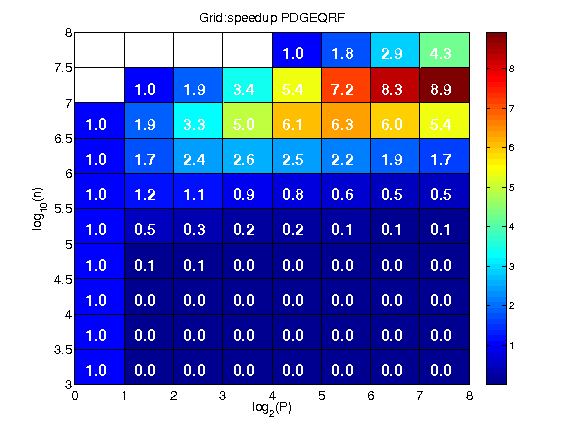}} 
  }
  \subfigure[Comparison]{\includegraphics[scale=0.35]{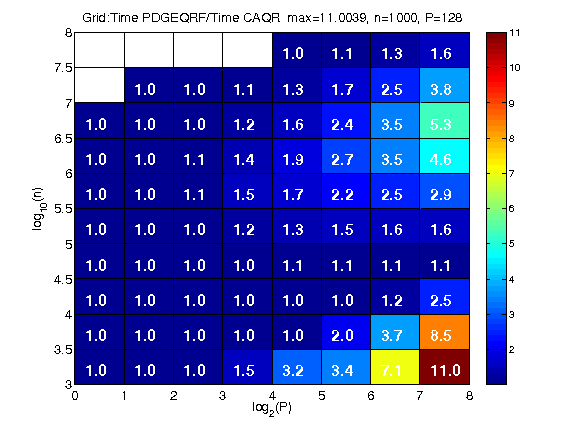}}
\end{center}
\caption{\label{fig:PerfComp_Grid}Performance prediction comparing
  CAQR and PDGEQRF on Grid.}
\end{figure}

\begin{figure}
  \begin{center}%
    \mbox{
      \subfigure[Optimal $b$]{\includegraphics[scale=0.35]{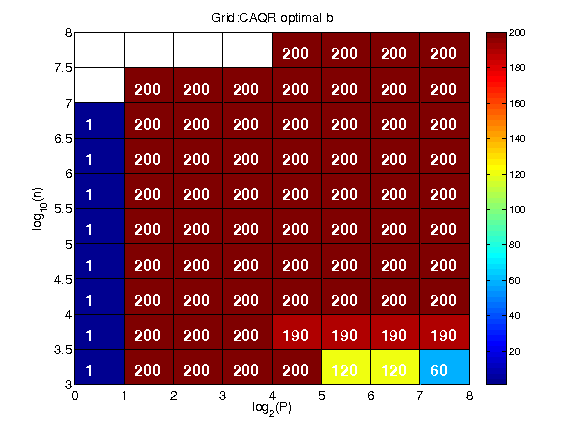}\label{grid_CAQRoptb}}
      \subfigure[Optimal $P_r$]{\includegraphics[scale=0.35]{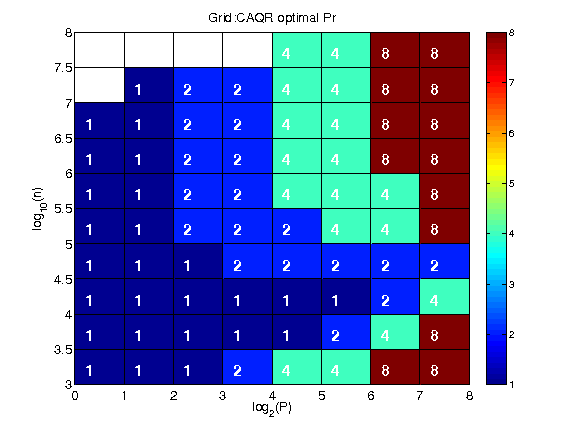}\label{grid_CAQRoptPr}} 
    }
    \mbox{
      \subfigure[Fraction of time in computation]{\includegraphics[scale=0.35]{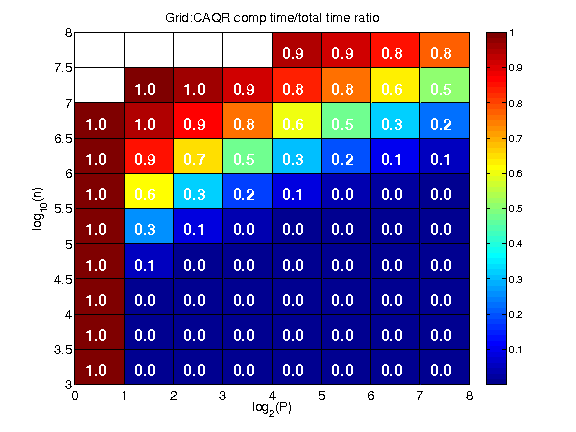}\label{grid_CAQRcompR}}
      \subfigure[Fraction of time in latency]{\includegraphics[scale=0.35]{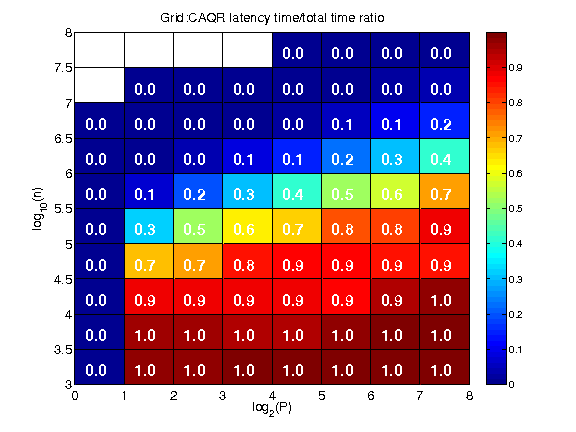}\label{grid_CAQRlatR}}
    }
    \mbox{
      \subfigure[Fraction of time in bandwidth]{\includegraphics[scale=0.35]{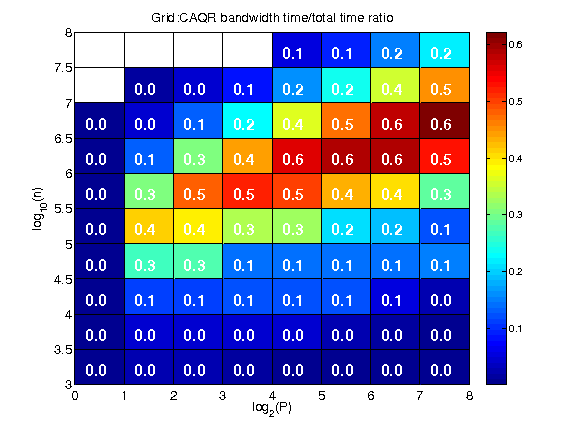}\label{grid_CAQRbwR}}
    }
  \end{center}
  \caption{\label{fig:PerfCAQR_Grid}Performance prediction for CAQR on
    Grid.}
\end{figure}

\begin{figure}
  \begin{center}
    \mbox{ 
      \subfigure[Optimal $b$]{\includegraphics[scale=0.35]{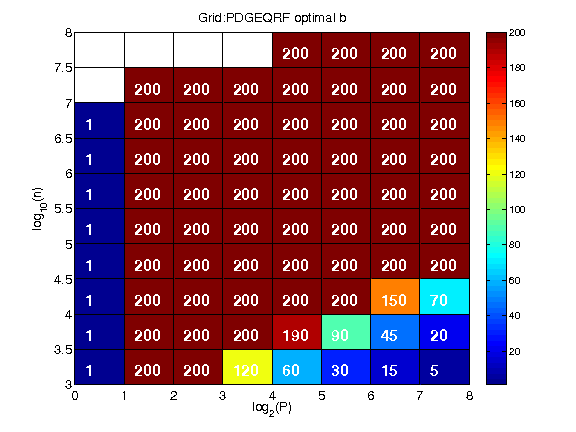}\label{grid_QRFoptb}}
      \subfigure[Optimal $P_r$]{\includegraphics[scale=0.35]{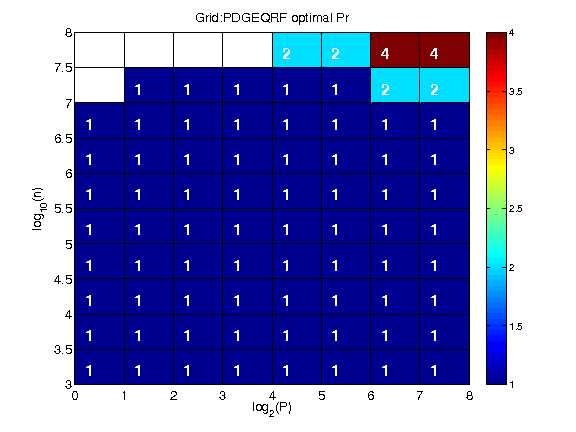}\label{grid_QRFoptPr}}
    } 
    \mbox{ 
      \subfigure[Fraction of time in computation]{\includegraphics[scale=0.35]{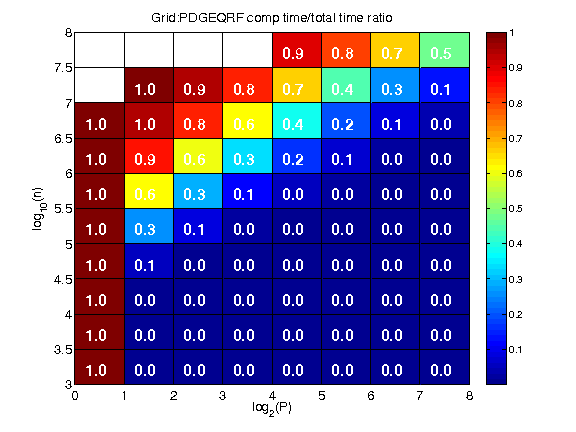}\label{grid_QRFcompR}}
      \subfigure[Fraction of time in latency]{\includegraphics[scale=0.35]{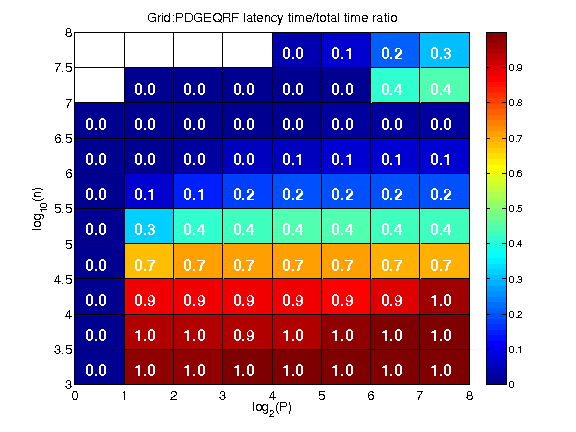}\label{grid_QRFlatR}}
    }
    \mbox{ 
      \subfigure[Fraction of time in bandwidth]{\includegraphics[scale=0.35]{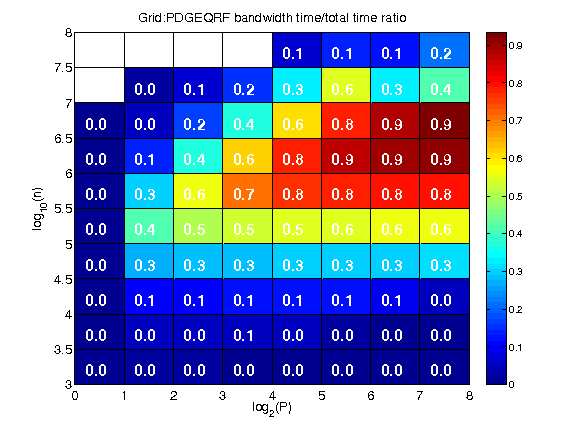}\label{grid_QRFbwR}}
    }
  \end{center}
  \caption{\label{fig:PerfPDGEQRF_Grid}Performance prediction for
    PDGEQRF on Grid.}
\end{figure}

The performance estimation obtained by CAQR and \lstinline!PDGEQRF! on
the Grid is displayed in Figures \ref{fig:PerfComp_Grid},
\ref{fig:PerfCAQR_Grid}, and \ref{fig:PerfPDGEQRF_Grid}.  For small
values of $n$ both algorithms do not obtain any speedup, even on small
number of processors.  Hence we discuss performance results for values
of $n$ bigger than $10^5$.

As displayed in Figures \ref{grid_CAQRoptb} and \ref{grid_QRFoptb},
the optimal block size for both algorithms is very often $200$, the
largest value in the allowed range.  The optimal value of $P_r$ for
\lstinline!PDGEQRF! is equal to $1$ for most of the test cases (Figure
\ref{grid_QRFoptPr}), while CAQR tends to prefer a square grid (Figure
\ref{grid_CAQRoptPr}).  This suggests that CAQR can successfully
exploit parallelism within block columns, unlike \lstinline!PDGEQRF!.

As can be seen in Figures \ref{grid_QRFcompR}, \ref{grid_QRFlatR}, and
\ref{grid_QRFbwR}, for small matrices, communication latency dominates
the total runtime of \lstinline!PDGEQRF!.  For large matrices and
smaller numbers of processors, computation dominates the runtime.  For
the test cases situated in the band going from the bottom left corner
to the upper right corner, bandwidth costs dominate the runtime.  The
model of \lstinline!PDGEQRF! suggests that the best way to decrease
the latency cost with this algorithm is to use, in most test cases, a
block column cyclic distribution (the layout obtained when $P_r = 1$).
In this case the bandwidth cost becomes significant.

The division of time between computation, latency, and bandwidth has a
similar pattern for CAQR, as shown in Figures \ref{grid_CAQRcompR},
\ref{grid_CAQRlatR}, and \ref{grid_CAQRbwR}.  However, unlike
\lstinline!PDGEQRF!, CAQR has as optimal grid shape a square or almost
square grid of processors, which suggests that CAQR is more scalable.

The best improvement is obtained for one of the largest matrix in our
test set $m=n=10^{6.5}$, where CAQR outperforms \lstinline!PDGEQRF! by
a factor of $5.3$ on $128$ processors.  The speedup obtained by the
best CAQR compared to the best \lstinline!PDGEQRF! is larger than $4$,
and the best performance is obtained by CAQR on $128$ processors,
while the best performance of \lstinline!PDGEQRF! is obtained on $32$
processors.

CAQR is predicted to obtain reasonable speedups for large problems on
the Grid, as displayed in Figure \ref{grid_CAQRspdup}.  For example,
for $n = 10^7$ we note a speedup of $33.4$ on $128$ processors
measured with respect to $2$ processors.  This represents an
improvement of $1.6$ over \lstinline!PDGEQRF!.  For the largest matrix
in the test set, $n=10^{7.5}$, we note a speedup of $6.6$ on $128$
processors, measured with respect to $16$ processors.  This is an
improvement of $3.8$ with respect to \lstinline!PDGEQRF!.

As with the last model, for any fixed $n$, we can take the number of
processors $P$ for which \lstinline!PDGEQRF! would perform the best,
and measure the speedup of CAQR over \lstinline!PDGEQRF! using that
number of processors.  We do this in Table
\ref{tbl:CAQR:par:Grid:best}, which shows that CAQR always is at
least as fast as \lstinline!PDGEQRF!, and often significantly faster
(up to $3.8 \times$ faster in some cases).

\begin{table}
\begin{tabular}{r|c|c}
$\log_{10} n$ & Best $\log_2 P$ for \lstinline!PDGEQRF! & CAQR speedup
\\ \hline
6.0 & 3     & 1.4 \\
6.5 & 5     & 2.4 \\
7.0 & 7     & 3.8 \\
7.5 & 7     & 1.6 \\
\end{tabular}
\caption{Estimated runtime of \lstinline!PDGEQRF! divided by estimated
  runtime of CAQR on a square $n \times n$ matrix, on the Grid
  platform, for those values of $P$ (number of processors) for which
  \lstinline!PDGEQRF! performs the best for that problem size.}
\label{tbl:CAQR:par:Grid:best}
\end{table}

\section{Lower bounds on communication for QR}\label{S:lowerbounds}

In this section, we review known lower bounds on communication 
bandwidth for parallel and sequential matrix-matrix multiplication 
of matrices stored in 1-D and 2-D block cyclic layouts, 
extend some of them to the rectangular case, and then extend them to QR, 
showing that our sequential and parallel TSQR and CAQR 
algorithms have optimal communication complexity with respect
to both bandwidth (in a Big-Oh sense, and sometimes modulo
polylogarithmic factors).

We will also use the simple fact that if $B$ is a lower bound on the number
of words that must be communicated to implement an algorithm,
and if $W$ is the size of the local memory (in the parallel case) 
or fast memory (in the sequential case),
so that $W$ is the largest possible size of a message,
then $B/W$ is a lower bound on the latency, i.e. the number
of messages needed to move $B$ words into or out of the memory.
We use this to derive lower bounds on latency, which are
also attained by our algorithms (again in a Big-Oh sense,
and sometimes modulo polylogarithmic factors).

We begin in section~\ref{SS:MMlowerbounds} by reviewing
known communication complexity bounds for matrix multiplication,
due first to Hong and Kung \cite{hong1981io} in the sequential case,
and later proved more simply and extended to the parallel case
by Irony, Toledo and Tiskin \cite{irony2004communication}.  
It is easy to extend lower bounds for matrix multiplication
to lower bounds for LU decomposition via the following 
reduction of matrix multiplication to LU decomposition:
\begin{equation}\label{eq:GEMM-to-LU}
\begin{pmatrix}
I & 0 & -B \\
A & I & 0 \\
0 & 0 & I \\
\end{pmatrix}
=
\begin{pmatrix}
I &   &   \\
A & I &   \\
0 & 0 & I \\
\end{pmatrix}
\begin{pmatrix}
I & 0 & -B \\
  & I & A \cdot B \\
  &   & I \\
\end{pmatrix}.
\end{equation}
See \cite{grigori2008calu} for an implementation of 
parallel LU that attains these bounds.
It is reasonable to expect that lower bounds for matrix multiplication
will also apply (at least in a Big-Oh sense) to other 
one-sided factorizations, such as QR.
Though as we will see, QR is not quite so simple.

All this work assumes commutative and associative
reorderings of the conventional $O(n^3)$ 
matrix multiplication
algorithm, and so excludes faster algorithms using distributivity
or special constants, such as those of
Strassen \cite{strassen1969gaussian} or Coppersmith and Winograd \cite{coppersmith1982asymptotic},
and their use in asymptotically fast versions of
LU and QR \cite{FastLinearAlgebraIsStable}.
Extending communication lower bounds to these asymptotically
faster algorithms is an open problem.

\subsection{Matrix Multiplication Lower Bounds}
\label{SS:MMlowerbounds}

Hong and Kung \cite{hong1981io},
and later Irony, Toledo and Tiskin
\cite{irony2004communication} considered the multiplication
of two $n$-by-$n$ matrices $C = A \cdot B$
using commutative and associative (but not distributive)
reorderings of the usual $O(n^3)$ algorithm.

In the sequential case, they assume that $A$ and $B$ initially reside
in slow memory, that there is a fast memory of size $W < n^2$,
and that the product $C = A \cdot B$ must be computed and eventually 
reside in slow memory.
They bound from below the number of
words that need to be moved between slow memory and fast memory
to perform this task:
\begin{equation}\label{eqn_MatMul_seq_bw_lowerbound}
{\rm \#\ words\ moved} \geq \frac{n^3}{2 \sqrt{2} W^{1/2}} - W \approx
\frac{n^3}{2 \sqrt{2} W^{1/2}}  \; \; .
\end{equation}
Since only $W$ words can be moved in one message, this also provides
a lower bound on the number of messages:
\begin{equation}\label{eqn_MatMul_seq_lat_lowerbound}
{\rm \#\ messages} \geq \frac{n^3}{2 \sqrt{2} W^{3/2}} - 1 \approx
\frac{n^3}{2 \sqrt{2} W^{3/2}} \; \; .
\end{equation}
In the rectangular case, where $A$ is $n$-by-$r$, $B$ is $r$-by-$m$,
and $C$ is $n$-by-$m$, so that the number of arithmetic operations in
the standard algorithm is $2mnr$, the above two results still apply,
but with $n^3$ replaced by $mnr$.

Irony, Toledo and Tiskin also consider the parallel case. There is
actually a spectrum of algorithms, from the so-called 2D case,
that use little extra memory beyond that needed to store equal fractions
of the matrices $A$, $B$ and $C$
(and so about $3n^2/P$ words for each of $P$ processors, in the square case), 
to the 3D case, where each input matrix is replicated up to $P^{1/3}$ times,
so with each processor needing memory of size $n^2/P^{2/3}$ in the square case.
We only consider the 2D case, which is the conventional, memory 
scalable approach.
In the 2D case, with square matrices, Irony et al show that
if each processor has $\mu n^2 /P$ words of local
memory, and $P \geq 32 \mu^3$, then at least one of the processors
must send or receive at least the following number of words:
\begin{equation}\label{eqn_MatMul_par_bw_lowerbound}
{\rm \#\ words\ sent\ or\ received} \geq \frac{n^2}{4 \sqrt{2} (\mu P)^{1/2}} 
\end{equation}
and so using at least the following number of messages
(assuming a maximum message size of $n^2/P$):
\begin{equation}\label{eqn_MatMul_par_lat_lowerbound}
{\rm \#\ messages} \geq \frac{P^{1/2}}{4 \sqrt{2} (\mu)^{3/2}} \; \; .
\end{equation}

We wish to extend this to the case of rectangular matrices.
We do this in preparation for analyzing CAQR in the rectangular case.
The proof is a simple extension of Thm.~4.1 in 
\cite{irony2004communication}.  

\theorem{
Consider the conventional matrix multiplication algorithm applied to
$C = A \cdot B$ where
$A$ is $n$-by-$r$, $B$ is $r$-by-$m$, and $C$ is $n$-by-$m$.
implemented on a $P$ processor distributed memory parallel computer.
Let $\bar{n}$, $\bar{m}$ and $\bar{r}$ be the sorted values of
$n$, $m$, and $r$, i.e. $\bar{n} \geq \bar{m} \geq \bar{r}$.
Suppose each processor has $3\bar{n}\bar{m}/P$ words of local memory,
so that it can fit 3 times as much as $1/P$-th of the largest of the
three matrices. Then as long as
\begin{equation}\label{eqn:RectMM}
\bar{r} \geq \sqrt{\frac{864 \bar{n}\bar{m}}{P}}
\end{equation}
(i.e. none of the matrices is ``too rectangular'')
then the number of words at least one processor must send or
receive is 
\begin{equation}\label{eqn:RectMM_bw}
{\rm \#\ words\ moved} \geq 
\frac{\sqrt{\bar{n}\bar{m}} \cdot \bar{r}}
{\sqrt{96 P}}
\end{equation}
and the number of messages is 
\begin{equation}\label{eqn:RectMM_lat}
{\rm \#\ messages} \geq 
\frac{ \sqrt{P} \cdot \bar{r}}
{ \sqrt{864 \bar{n}\bar{m}} } 
\end{equation}
}

\rm

\begin{proof}
We use (\ref{eqn_MatMul_seq_bw_lowerbound}) 
with $\bar{m} \bar{n} \bar{r}/P$ substituted for $n^3$,
since at least one processor does this much arithmetic,
and $W = 3\bar{n}\bar{m}/P$ words of local memory.
The constants in inequality (\ref{eqn:RectMM})
are chosen so that the first term in
(\ref{eqn_MatMul_seq_bw_lowerbound}) is at least $2W$,
and half the first term is a lower bound.
\end{proof}

\subsection{Lower bounds for TSQR}\label{SS:lowerbounds:1d}

TSQR with data stored in a 1D layout 
is simpler than the general CAQR case, and does
not depend on the above bounds for matrix multiplication.

\subsubsection{Sequential TSQR}

In the sequential case, the $m \times n$ matrix $A$ must be read from
slow memory into fast memory at least once, if we assume that fast
memory is empty at the start of the computation, and the answer written
out to slow memory.
Thus, the number of
words transferred (the bandwidth lower bound) is at least $2mn$.
As described in Table~\ref{tbl:QR:perfcomp:seq} 
in Section~\ref{S:TSQR:perfcomp}, and in more detail in Appendix~B,
our sequential TSQR moves 
\[
\frac{mn^2}{W - n(n+1)/2}
+ 2mn 
- \frac{n(n+1)}{2} 
\]
words.  Since we assume $W \geq \frac{3}{2} n^2$, this is little more
than the lower bound $2mn$.
In contrast,
blocked left-looking Householder QR moves
\[
+ \frac{m^2 n^2}{2 W} 
- \frac{m n^3}{6W}
+ \frac{3mn}{2} 
- \frac{3n^2}{4} 
\]
words, where the first and second terms combined can be $O(\frac{mn}{W})$ 
times larger than the lower bound (see Table~\ref{tbl:QR:perfcomp:seq});
note that $\frac{mn}{W}$ is how many times larger the matrix is than
the fast memory.

The number of slow memory reads and writes (the latency lower bound)
is at least $2mn/W$.
As described in the same sections as before, our sequential TSQR sends
\[
\frac{2mn}{W - n(n+1)/2}
\]
words, which is close to the lower bound.  
In contrast, blocked left-looking Householder QR sends 
\[
\frac{2mn}{W}
+ 
\frac{mn^2}{2W}
\] 
messages, which can be $O(n)$ times larger than the
lower bound (see Table~\ref{tbl:QR:perfcomp:seq}).

\subsubsection{Parallel TSQR}
\label{sec:par_TSQR_comm_bounds}

In the parallel case, we prefer for 1-D layouts to distinguish between
the minimum number of messages per processor, and the number of
messages along the critical path.  For example, one can perform a
reduction linearly, so that each processor only sends one message to
the next processor.  This requires $P - 1$ messages along the critical
path, but only one message per processor.  A lower bound on the
minimum number of sends or receives performed by any processor is also
a lower bound on the number of messages along the critical path.  The
latter is more difficult to analyze for 2-D layouts, so we only look
at the critical path for 1-D layouts.  By the usual argument that any
nontrivial function of data distributed across $P$ processors requires
at least $\log_2 P$ messages to compute, the critical path length
$C_{\text{1-D}}(m,n,P)$ satisfies
\begin{equation}\label{eq:lowerbound:1d:par:lat}
C_{\text{1-D}}(m,n,P) \geq \log_2 P.
\end{equation}
This is also the number of messages required per processor along the
critical path. This lower bound is obviously attained by parallel TSQR
based on a binary tree.

Appendix~\ref{S:CommLowerBoundsFromCalculus}
shows formally that for any reduction tree computing the QR
decomposition of $\frac{m}{p} \times n$ matrices at its leaves, each
path from the leaves to the root must send at least $n(n+1)/2$ words
of information along each edge. This means the bandwidth cost is at
least $n(n+1)/2$ times the length of critical path, or at least
$\log(P)n(n+1)/2$.  This is clearly attained by TSQR (see Table
\ref{tbl:QR:perfcomp:par}).

\subsection{Lower Bounds for CAQR}\label{SS:lowerbounds:2d}

Now we need to extend our analysis of matrix multiplication.
We assume all variables are real; extensions to the complex
case are straightforward.
Suppose $A = QR$ is $m$-by-$n$, $n$ even, 
so that 
\[
\bar{Q}^T \cdot \bar{A} 
\equiv 
\left(
    Q( 1:m, 1:\frac{n}{2} )
\right)^T 
\cdot 
A( 1:m, \frac{n}{2}+1:n )
= 
R( 1:\frac{n}{2}, \frac{n}{2}+1:n )
\equiv 
\bar{R} \; \; .
\]
It is easy to see that $\bar{Q}$
depends only on the first $\frac{n}{2}$ columns of $A$, and
so is independent of $\bar{A}$. The obstacle to directly
applying existing lower bounds for matrix multiplication 
of course is that $\bar{Q}$ is not represented as an explicit 
matrix, and $\bar{Q}^T \cdot \bar{A}$ is not implemented by 
straightforward matrix multiplication. 
Nevertheless, we argue that the same data dependencies
as in matrix multiplication can be found inside many
implementations
of $\bar{Q}^T \cdot \bar{A}$,
and that therefore
the geometric ideas underlying the analysis in
\cite{irony2004communication} still apply.
Namely, there are two data structures $\tilde{Q}$
and $\tilde{A}$ indexed with pairs of subscripts
$(j,i)$ and $(j,k)$ respectively with the following properties.
\begin{itemize}
\item $\tilde{A}$ stores $\bar{A}$ as well as all intermediate
results which may overwrite $\bar{A}$. 
\item $\tilde{Q}$ represents $\bar{Q}$, i.e., an $m$-by-$\frac{n}{2}$
orthogonal matrix. Such a matrix is a member of the Stiefel manifold
of orthogonal matrices, and is known to require 
$\frac{mn}{2} - \frac{n}{4}(\frac{n}{2}+1)$ independent parameters
to represent, with column $i$ requiring $m-i$ parameters,
although a particular algorithm may represent
$\bar{Q}$ using more data.
\item The algorithm operates mathematically independently on each column
of $\bar{A}$, i.e., methods like that of Strassen are excluded.
This means that the algorithm performs at least 
$\frac{mn}{2} - \frac{n}{4}(\frac{n}{2}+1)$ multiplications 
on each $m$-dimensional column vector of $\bar{A}$
(see subsection~\ref{SS:lowerbounds:2d:flops} for a proof),
and does the same operations on each column of $\bar{A}$. 
\item For each $(i,k)$ indexing $\bar{R}_{i,k}$, which is the 
component of the $k$-th column $\bar{A}_{:,k}$ of $\bar{A}$ in
the direction of the $i$-th column $\bar{Q}_{:,i}$ of $\bar{Q}$,
it is possible to identify at least $m-i$ common components of $\tilde{A}_{:,k}$
and of $\tilde{Q}_{:,i}$ such that a parameter associated with
$\tilde{Q}_{j,i}$ is multiplied by a value stored in $\tilde{A}_{j,k}$.
\end{itemize}
The last point, which says that $\bar{Q}^T \cdot \bar{A}$
has at least the same dependencies as matrix multiplication, 
requires illustration.
\begin{itemize}
\item Suppose $\bar{Q}$ is represented as a product of $\frac{n}{2}$
Householder reflections with a projection $\hat{Q}$ onto the
first $\frac{n}{2}$ coordinates,
$\bar{Q} = 
(I - \tau_1 u_1u_1^T)
	\cdots 
(I - \tau_{n/2} u_{n/2} u_{n/2}^T) 
\hat{Q}$,
normalized in the 
conventional way where the topmost nonzero entry of each $u_j$ is one,
and $\hat{Q}$ consists of the first $n/2$ columns of the $n$-by-$n$
identity matrix.
Then $\tilde{Q}_{j,i} = u_i(j)$ is multiplied by some intermediate value of
$\bar{A}_{j,k}$, i.e. $\tilde{A}_{j,k}$.
\item Suppose $\bar{Q}$ is represented as a product of block
Householder transformations $(I-Z_1U_1^T) \cdots (I-Z_f U_f^T) \hat{Q}$
where $U_g$ and $Z_g$ are $m$-by-$b_g$ matrices, $U_g$ consisting
of $b_g$ Householder vectors side-by-side. 
Again associate $\tilde{Q}_{j,i}$ with
the $j$-th entry of the $i$-th Householder vector $u_i(j)$.
\item Recursive versions of QR \cite{elmroth1998new} apply 
blocked Householder transformations organized so as to better
use BLAS3, but still let us use the approach of the last bullet.
\item Suppose $\bar{Q}$ is represented as a product of 
$\frac{mn}{2} - \frac{n}{4}(\frac{n}{2}+1)$ Givens rotations,
each one creating a unique subdiagonal zero entry in $A$ which is
never filled in.  There are
many orders in which these zeros can be created, and possibly
many choices of row that each Givens rotation may rotate with to zero
out its desired entry. 
If the desired zero entry in $A_{j,i}$ is created by
the rotation in rows $j'$ and $j$, $j'<j$, 
then associate $\tilde{Q}_{j,i}$ with the value
of the cosine in the Givens rotation, since this will be multiplied
by $\bar{A}_{j,k}$.
\item Suppose, finally, that we use CAQR to perform the
QR decomposition, so that $\bar{Q} = Q_1 \cdots Q_f \hat{Q}$, where
each $Q_g$ is the result of TSQR on $b_g$ columns.
Consider without loss of generality $Q_1$, which operates
on the first $b_1$ columns of $A$.
We argue that TSQR still produces $m-i$ parameters associated
with column $i$ as the above methods. Suppose there are $P$
row blocks, each of dimension $\frac{m}{P}$-by-$b_1$.
Parallel TSQR initially does QR independently on each block, using
any of the above methods; we associate multipliers as 
above with the subdiagonal entries in each block. Now
consider the reduction tree that combines $q$ different $b_1$-by-$b_1$
triangular blocks at any particular node.
As described in subsection~\ref{SS:TSQR:localQR:structured}, this generates
$(q-1)b_1(b_1+1)/2$ parameters that multiply the equal number of entries 
of the $q-1$ triangles being zeroed out, and so can be associated with
appropriate entries of $\tilde{Q}$. Following the reduction tree, we
see that parallel TSQR produces exactly 
as many parameters as Householder reduction,
and that these may be associated one-for-one with all subdiagonal
entries of $\tilde{Q}(:,1:b_1)$ and $\tilde{A}(:,1:b_1)$ as above.
Sequential TSQR reduction is analogous.
\end{itemize}
We see that we have only tried to capture the dependencies
of a fraction of the arithmetic operations performed by various
QR implementations; this is all we need for a lower bound.

Now we resort to the geometric approach of 
\cite{irony2004communication}: Consider a three dimensional 
block of lattice points, indexed by $(i,j,k)$.
Each point on the $(i,0,k)$ face is associated with $\bar{R}_{i,k}$,
for $1 \leq i,k \leq \frac{n}{2}$.
Each point on the $(0,j,k)$ face is associated with $\tilde{A}_{j,k}$,
for $1 \leq k \leq \frac{n}{2}$ and $1 \leq j \leq m$.
Each point on the $(i,j,0)$ face is associated with $\tilde{Q}_{j,i}$,
for $1 \leq i \leq \frac{n}{2}$ and $1 \leq j \leq m$.
Finally, each interior point $(i,j,k)$ for 
$1 \leq i,k \leq \frac{n}{2}$ and $1 \leq j \leq m$ represents the
multiplication $\tilde{Q}_{j,i} \cdot \tilde{A}_{j,k}$.
The point is that the multiplication at $(i,j,k)$ cannot occur
unless $\tilde{Q}_{j,i}$ and $\tilde{A}_{j,k}$ are together in memory.

Finally, we need the Loomis-Whitney inequality \cite{loomis1949inequality}:
Suppose $V$ is a set of lattice points in 3D, 
$V_i$ is projection of $V$ along $i$ onto the $(j,k)$ plane,
and similarly for $V_j$ and $V_k$. Let $|V|$ denote the cardinality of $V$,
i.e. counting lattice points. Then 
$|V|^2 \leq |V_i| \cdot |V_j| \cdot |V_k|$.

Finally, we can state

\lemma{
Suppose a processor with local (fast) memory of size $W$ is participating
in the QR decomposition of an $m$-by-$n$ matrix, $m \geq n$, using an
algorithm of the sort discussed above.
There may or may not be other processors participating (i.e. this lemma covers
the sequential and parallel cases). Suppose the processor performs $F$
multiplications. Then the processor must move
the following number of works into or out of its memory:
\begin{equation}\label{Thm:1_bw}
{\rm \#\ of\ words\ moved} \geq \frac{F}{(8W)^{1/2}} - W
\end{equation}
using at least the following number of messages:
\begin{equation}\label{Thm:1_lat}
{\rm \#\ of\ messages} \geq \frac{F}{(8W^3)^{1/2}} - 1
\end{equation}
}
\label{lemma:LB}

\rm
\begin{proof}
The proof closely follows that of Lemma~3.1 in \cite{irony2004communication}.
We decompose the computation into phases. Phase $l$ begins when the
total number of words moved into and out of memory is exactly $lW$.
Thus in each phase, except perhaps the last, the memory loads and stores
exactly $W$ words.

The number of words $n_A$ from different $\tilde{A}_{jk}$ that the processor
can access in its memory during a phase is $2W$, since each word was either
present at the beginning of the phase or read during the phase.
Similarly the number of coefficients $n_Q$ from different $\tilde{Q}_{ji}$
also satisfies $n_Q \leq 2W$. Similarly, the number $n_R$ of locations
into which intermediate results like $\tilde{Q}_{ji} \cdot \tilde{A}_{jk}$
can be accumulated or stored is at most $2W$. Note that these intermediate
results could conceivably be stored or accumulated in $\tilde{A}$ because 
of overwriting; this does not affect the upper bound on $n_R$.

By the Loomis-Whitney inequality, the maximum number of useful multiplications
that can be done during a phase (i.e. assuming intermediate results are
not just thrown away) is bounded by 
$\sqrt{n_A \cdot n_Q \cdot n_R} \leq \sqrt{8W^3}$. Since the processor does
$F$ multiplications, the number of full phases required is at least
\[
\left\lfloor \frac{F}{\sqrt{8W^3}} \right\rfloor \geq \frac{F}{\sqrt{8W^3}} -1
\]
so the total number of words moved is $W$ times larger, i.e. at least
\[
{\rm \#\ number\ of\ words\ moved} \geq
\frac{F}{\sqrt{8W}} -W \; \; .
\]
The number of messages follows by dividing by $W$, the maximum
message size.
\end{proof}
\rm

\subsubsection{Sequential CAQR}
\label{sec:SeqCAQR}

\corollary{
Consider a single processor computing the QR decomposition of
an $m$-by-$n$ matrix with $m \geq n$, using an algorithm of the
sort discussed above. Then the number of words moved between
fast and slow memory is at least
\begin{equation}\label{Thm:2_bw}
{\rm \#\ of\ words\ moved} \geq \frac
{\frac{mn^2}{4}-\frac{n^2}{8}(\frac{n}{2}+1)}
{(8W)^{1/2}} - W
\end{equation}
using at least the following number of messages:
\begin{equation}\label{Thm:2_lat}
{\rm \#\ of\ messages} \geq \frac
{\frac{mn^2}{4}-\frac{n^2}{8}(\frac{n}{2}+1)}
{(8W^3)^{1/2}} - 1
\end{equation}
}
\label{corollary:SeqCAQR}
\rm
\begin{proof}
The proof follows easily from Lemma~\ref{lemma:LB} by
using the lower bound
$F \geq {\frac{mn^2}{4}-\frac{n^2}{8}(\frac{n}{2}+1)}$ on the number
of multiplications by any algorithm in the class discussed above
(see Lemma~\ref{lemma:F_lowerbound} in
subsection~\ref{SS:lowerbounds:2d:flops} for a proof).
\end{proof}

\rm
The lower bound could be increased by a constant factor by the
using specific number of multiplications (say $mn^2 - \frac{1}{3}n^3$
using Householder reductions), instead of arguing more generally based on
the number of parameters needed to represent orthogonal matrices.

Comparing to Equation \eqref{eq:CAQR:seq:modeltime:P:opt} in Appendix
\ref{S:CAQR-seq-detailed} or the presentation in Section \ref{S:CAQR-seq},
we see that CAQR attains these bounds to within a constant factor.

\subsubsection{Parallel CAQR}

\corollary{
Consider a parallel computer with $P$ processors
and $W$ words of memory per processor
computing the QR decomposition of
an $m$-by-$n$ matrix with $m \geq n$, using an algorithm of the
sort discussed above. 
Then the number of words sent and received
by at least one processor
is at least
\begin{equation}\label{Thm:3a_bw}
{\rm \#\ of\ words\ moved} \geq \frac
{\frac{mn^2}{4}-\frac{n^2}{8}(\frac{n}{2}+1)}
{P(8W)^{1/2}} - W
\end{equation}
using at least the following number of messages:
\begin{equation}\label{Thm:3a_lat}
{\rm \#\ of\ messages} \geq \frac
{\frac{mn^2}{4}-\frac{n^2}{8}(\frac{n}{2}+1)}
{P(8W^3)^{1/2}} - 1
\end{equation}
In particular, when each processor has $W = mn/P$ words of memory
and the matrix is not too rectangular, $n \geq \frac{2^{11}m}{P}$,
then the number of words sent and received
by at least one processor is at least
\begin{equation}\label{Thm:3b_bw}
{\rm \#\ of\ words\ moved} \geq 
\sqrt{\frac{m n^3}{2^{11}P}}
\end{equation}
using at least the following number of messages:
\begin{equation}\label{Thm:3b_lat}
{\rm \#\ of\ messages} \geq 
\sqrt{\frac{nP}{2^{11}m}} \; \; .
\end{equation}
In particular, in the square case $m=n$, we get that
as long as $P \geq 2^{11}$, 
then the number of words sent and received
by at least one processor is at least
\begin{equation}\label{Thm:4b_bw}
{\rm \#\ of\ words\ moved} \geq 
{\frac{n^2}{2^{11/2}P^{1/2}}}
\end{equation}
using at least the following number of messages:
\begin{equation}\label{Thm:4b_lat}
{\rm \#\ of\ messages} \geq 
\sqrt{\frac{P}{2^{11}}} \; \; .
\end{equation}
}
\rm
\begin{proof}
The result follows from the previous Corollary, since
at least one processor has to do $1/P$-th of the work.
\end{proof}

Comparing to Equations \eqref{eq:CAQR:par:opt:lat} and
\eqref{eq:CAQR:par:opt:bw} in Section \ref{SS:CAQR:par:opt}, we see
that CAQR attains these bounds to within a polylog factor.

\subsection{Lower Bounds on Flop Counts for QR}
\label{SS:lowerbounds:2d:flops}

This section proves lower bounds on arithmetic for {\em any} ``columnwise'' 
implementation of QR, by which we mean one whose operations can be reordered
so as to be left looking, i.e. the operations that compute columns $i$
of $Q$ and $R$ depend on data only in columns 1 through $i$ of $A$.
The mathematical dependencies are such that columns $i$ of $Q$ and $R$
do only depend on columns 1 through $i$ of $A$, but saying that operations
only depend on these columns eliminates algorithms like Strassen.
(It is known that QR can be done asymptotically as fast as any fast
matrix multiplication algorithm like Strassen, and stable as well
\cite{FastLinearAlgebraIsStable}.)

This section says where the lower bound on
$F$ comes from that is used in the proof of 
Corollary~\ref{corollary:SeqCAQR} above.

The intuition is as follows. 
Suppose $A = QR$ is $m$-by-$(j+1)$, 
so that 
\[
\bar{Q}^T \cdot \bar{A} 
\equiv (Q(1:m,1:j))^T \cdot A(1:m,j+1) = 
R(1:j, j+1)
\equiv \bar{R} \; \; .
\]
where $\bar{Q}$ only depends on the first $j$ columns of $A$, and
is independent of $\bar{A}$. As an arbitrary $m$-by-$j$ orthogonal 
matrix, a member of the Stiefel manifold of dimension
$mj-j(j+1)/2$, $\bar{Q}$ requires $mj-j(j+1)/2$ independent
parameters to represent. We will argue that no matter how $\bar{Q}$
is represented, i.e. without appealing to the special structure of
Givens rotations or Householder transformations, that unless
$mj-j(j+1)/2$ multiplications are performed to compute $\bar{R}$
it cannot be computed correctly, because it cannot depend on
enough parameters. 

Assuming for a moment that this is true, we get a lower bound
on the number of multiplications needed for QR on an $m$-by-$n$ matrix 
by summing 
$\sum_{j=1}^{n-1} [mj-j(j+1)/2] = \frac{mn^2}{2} - \frac{n^3}{6} + O(mn)$.
The two leading terms are half the multiplication count for Householder
QR (and one fourth of the total operation count, including additions).
So the lower bound is rather tight.

Again assuming this is true, we get a lower bound on the
value $F$ in Corollary~\ref{corollary:SeqCAQR} by multiplying 
$\frac{n}{2} \cdot (m\frac{n}{2} - \frac{n}{2}(\frac{n}{2}+1)/2)
= \frac{mn^2}{4} - \frac{n^2}{8}(\frac{n}{2}+1) \leq F$.

Now we prove the main assertion, that $mj-j(j+1)/2$ multiplications are needed
to compute the single column $\bar{R} = \bar{Q}^T \cdot \bar{A}$, no matter how 
$\bar{Q}$ is represented. We model the computation as a DAG (directed
acyclic graph) of operations with the following properties, which we
justify as we state them.
\begin{enumerate}
\item There are $m$ input nodes labeled by the $m$ entries of $\bar{A}$,
$a_{1,j+1}$ through $a_{m,j+1}$. 
We call these $\bar{A}$-input nodes for short.
\item There are at least $mj-j(j+1)/2$ input nodes labeled by parameters
representing $\bar{Q}$, since this many parameters are needed to
represent a member of the Stiefel manifold.
We call these $\bar{Q}$-input nodes for short.
\item There are two types of computation nodes, addition and multiplication.
In other words, we assume that we do not do divisions, square roots, etc. 
Since we are only doing matrix multiplication, this is reasonable.
We note that any divisions or square roots in the overall algorithm
may be done in order to compute the parameters represented $\bar{Q}$.
Omitting these from consideration only lowers our lower bound
(though not by much).
\item There are no branches in the algorithm. In other words, the
way an entry of $\bar{R}$ is computed does not depend on the numerical
values. This assumption reflects current algorithms, but could in fact 
be eliminated as explained later.
\item Since the computation nodes only do multiplication and addition, 
we may view the output of each node as a polynomial in entries of $\bar{A}$
and parameters representing $\bar{Q}$.
\item We further restrict the operations performed so that the output of
any node must be a homogeneous linear polynomial in the entries of $\bar{A}$.
In other words, we never multiply two quantities depending on entries
of $\bar{A}$ to get a quadratic or higher order polynomial,
or add a constant or parameter depending on $\bar{Q}$ to an entry of
$\bar{A}$.  This is
natural, since the ultimate output is linear and homogeneous in $\bar{A}$, 
and any higher degree polynomial terms or constant terms would have to 
be canceled away.  No current or foreseeable algorithm (even Strassen based) 
would do this, and numerical stability would likely be lost.
\item There are $j$ output nodes labeled by the entries of $\bar{R}$,
$r_{1,j+1}$ through $r_{j,j+1}$.
\end{enumerate}

The final requirement means that multiplication nodes are only allowed
to multiply $\bar{Q}$-input nodes and homogeneous linear functions
of $\bar{A}$, including $\bar{A}$-input nodes.
Addition nodes may add homogeneous linear functions of $\bar{A}$ 
(again including $\bar{A}$-input nodes), but not add $\bar{Q}$-input nodes
to homogeneous linear functions of $\bar{A}$.
We exclude the possibility of adding or multiplying $\bar{Q}$-input nodes,
since the results of these could just be represented as additional
$\bar{Q}$-input nodes.

Thus we see that the algorithm represented by the DAG just described
outputs $j$ polynomials that are homogeneous and linear in $\bar{A}$.
Let $M$ be the total number of multiplication nodes in the DAG.
We now want to argue that unless $M \geq mj-j(j+1)/2$,
these output polynomials cannot possibly compute the right answer. 
We will do this by arguing that the dimension of
a certain algebraic variety they define is both bounded above by $M$,
and the dimension must be at least $mj-j(j+1)/2$ to get the right answer.

Number the output nodes from $1$ to $j$.
The output polynomial representing node $i$ can be written as
$\sum_{k=1}^m p_{k,i} (\bar{Q}) a_{k,j+1}$, where $p_{k,i}(\bar{Q})$ is
a polynomial in the values of the $\bar{Q}$-input nodes. According
to our rules for DAGs above, only multiplication nodes can introduce
a dependence on a previously unused $\bar{Q}$-input node, so
all the $p_{k,i}(\bar{Q})$ 
can only depend on $M$ independent parameters.

Finally, viewing each output node as a vector
of $m$ coefficient polynomials $(p_{1,i} (\bar{Q}),...,p_{m,i} (\bar{Q}))$,
we can view the entire output as a vector of $mj$ coefficient polynomials
$V(\bar{Q}) = (p_{1,1}(\bar{Q}),...,p_{m,j}(\bar{Q}))$, 
depending on $M$ independent parameters. 
This vector of length $mj$ needs to represent the set of
all $m$-by-$j$ orthogonal matrices. But the Stiefel manifold of
such orthogonal matrices has dimension $mj-j(j+1)/2$, so the surface
defined by $V$ has to have at least this dimension, i.e. $M \geq mj-j(j+1)/2$.

As an extension, we could add branches to our algorithm by noting that the
output of our algorithm would be piecewise polynomials, on regions
whose boundaries are themselves defined by varieties in the
same homogeneous linear polynomials. We can apply the above argument
on all the regions with nonempty interiors to argue that the same
number of multiplications is needed.

In summary, we have proven

\lemma{
Suppose we are doing the QR factorization of an $m$-by-$n$ 
matrix using any ``columnwise'' algorithm in the sense described 
above. Then at least $mn - j(j+1)/2$ multiplications are required
to compute column $j+1$ of $R$, and at least
$\frac{mn^2}{4} - \frac{n^2}{8}(\frac{n}{2} + 1)$
multiplications to compute columns $\frac{n}{2}+1$ through $n$
of $R$.
}
\label{lemma:F_lowerbound}
\rm

\section{Lower bounds on parallelism}\label{S:limits-to-par}\label{S:QR:opt:par}

We base this paper on the premise that communication costs matter more
than computation costs.  Many authors developed parallel algorithms
and bounds based on a PRAM model, which assumes that communication is
essentially free.  Their bounds are nevertheless of interest because
they provide fundamental limits to the amount of parallelism that can
be extracted, regardless of the cost of communication.

A number of authors have developed parallel QR factorization methods
based on Givens rotations (see e.g.,
\cite{sameh78:_stabl_solver,modi84:_given,cosnard83:_qr}).  Givens
rotations are a good model for a large class of QR factorizations,
including those based on Householder reflections.  This is because all
such algorithms have a similar dataflow graph (see e.g.,
\cite{leoncini1999parallel}), and are all based on orthogonal linear
transformations (so they are numerically stable).  Furthermore, these
bounds also apply to methods that perform block Givens rotations, if
we consider each block as an ``element'' of a block matrix.

\subsection{Minimum critical path length}\label{SS:limits-to-par:critpath}

Cosnard, Muller, and Robert proved lower bounds on the critical path
length $Opt(m,n)$ of any parallel QR algorithm of an $m \times n$
matrix based on Givens rotations \cite{cosnard86}.  They assume any
number of processors, any communication network, and any initial data
distribution; in the extreme case, there many be $mn$ processors, each
with one element of the matrix.  In their class of algorithms, a
single step consists of computing one Givens rotation and zeroing out
one matrix entry.  Their first result concerns matrices for which
$\lim_{m,n \to \infty} m/n^2 = 0$.  This includes the case when $m =
n$.  Then the minimum critical path length is
\begin{equation}\label{eq:optimal:qr:squarish}
C_{\text{Par.\ QR, squarish}} \geq 2n + o(n).
\end{equation}
A second complexity result is obtained for the case when $m \to
\infty$ and $n$ is fixed -- that is, ``tall skinny'' matrices.  Then,
the minimum critical path length is
\begin{equation}\label{eq:optimal:qr:tallskinny}
C_{\text{Par.\ QR, tall skinny}} \geq
\log_2 m + (n-1) \log_2\left( \log_2 m \right) + 
o\left( \log_2 \left( \log_2 m \right) \right).
\end{equation}

The above bounds apply to 1-D and 2-D block (cyclic) layouts if we
consider each ``row'' as a block row, and each ``column'' as a block
column.  One step in the computation involves computing one block
Givens rotation and applying it (i.e., either updating or eliminating
the current block).  Then, Equation \eqref{eq:optimal:qr:squarish}
shows in the case of a square matrix that the critical path length is
twice the number of block columns.  (This makes sense, because the
current panel must be factored, and the trailing matrix must be
updated using the current panel factorization; these are two dependent
steps.)  In the case of a tall skinny matrix in a 1-D block row
layout, Equation \eqref{eq:optimal:qr:tallskinny} shows that the
critical path length is $\log_2 (m/P)$, in which $P$ is the number of
processors.  (The $(n-1) \log_2\left( \log_2 m \right)$ term does not
contribute, because there is only one block column, so we can say that
$n = 1$.)

\subsection{Householder or Givens QR is P-complete}

Leoncini et al.\ show that any QR factorization based on Householder
reductions or Givens rotations is P-complete
\cite{leoncini1999parallel}.  This means that if there exists an
algorithm that can solve this problem using a number of processors
polynomial in the number of matrix entries, in a number of steps
polynomial in the logarithm of the number of matrix entries
(``polylogarithmic''), then \emph{all} tractable problems for a
sequential computer (the set P) can be solved in parallel in
polylogarithmic time, given a polynomial number of processors (the set
NC).  This ``P equals NC'' conclusion is considered unlikely, much as
``P equals NP'' is considered unlikely.  

Note that one could compute the QR factorization of a matrix $A$ by
multiplying $A^T \cdot A$, computing the Cholesky factorization $R
\cdot R^T$ of the result, and then performing $Q := A R^{-1}$.  We
describe this method (``CholeskyQR'') in detail in Section
\ref{S:TSQR:perfcomp}.  Csanky shows arithmetic NC algorithms for
inverting matrices and solving linear systems, and matrix-matrix
multiplication also has an arithmetic NC algorithm
\cite{csanky1976fast}.  Thus, we could construct a version of
CholeskyQR that is in arithmetic NC.  However, this method is highly
inaccurate in floating-point arithmetic.  Not only is CholeskyQR
itself inaccurate (see Section \ref{S:TSQR:stability}), Demmel
observes that Csanky's arithmetic NC linear solver is so unstable that
it loses all significant digits when inverting $3I_{n \times n}$ in
IEEE 754 double-precision arithmetic, for $n \geq 60$
\cite{demmel1992trading}.  As far as we know, there exists no stable,
practical QR factorization that is in arithmetic NC.

\section{Extending algorithms and optimality proofs to general architectures}
\label{S:hierarchies}

Our TSQR and CAQR algorithms have been described and analyzed 
in most detail for simple machine models: either
sequential with two levels of memory hierarchy (fast and slow),
or a homogeneous parallel machine, where each processor is
itself sequential. Real computers are more complicated, with
many levels of memory hierarchy and many levels of parallelism
(multicore, multisocket, multinode, multirack, \dots) all with
different bandwidths and latencies. So it is natural to ask
whether our algorithms and optimality proofs can be extended
to these more general situations. We have briefly described
how TSQR could be extended to general
reduction trees in Section~\ref{SS:TSQR:GeneralTrees}, which
could in turn be chosen depending on the architecture.
But we have not discussed CAQR, which we do here.

We again look at the simpler case of matrix multiplication
for inspiration. Consider the sequential case, with
$k$ levels of memory hierarchy instead of 2, where
level 1 is fastest and smallest with $W_1$ words of memory,
level 2 is slower and larger with $W_2$ words of memory,
and so on, with level $k$ being slowest and large enough
to hold all the data. By dividing this hierarchy into
two pieces, levels $k$ through $i+1$ ("slow")
and $i$ through 1 ("fast"), we can apply the theory in
Section~\ref{SS:MMlowerbounds} to get lower bounds
on bandwidth and latency for moving data between levels $i$
and $i+1$ of memory. So our goal expands to finding
a matrix multiplication algorithm that attains not just
1 set of lower bounds, but $k-1$ sets of lower bounds,
one for each level of the hierarchy. 

Fortunately, as is well known, the standard approach to tiling 
matrix multiplication achieves all these lower bounds simultaneously, 
by simply applying it recursively: level $i+1$ holds submatrices
of dimension $O(\sqrt{W_{i+1}})$, and multiplies them by tiling
them into submatrices of dimension $O(\sqrt{W_i})$, and so on.

The analogous observation is true of parallel matrix multiplication
on a hierarchical parallel processor where each node in the parallel 
processor is itself a parallel processor (multicore, multisocket, 
multirack, \dots).

We believe that this same recursive hierarchical approach applies
to CAQR (and indeed much of linear algebra) but there is a catch:
Simple recursion does not work, because the subtasks are not all
simply smaller QR decompositions. Rather they are a mixture of
tasks, including small QR decompositions and operations like 
matrix multiplication. Therefore we still expect that the same
hierarchical approach will work: if a subtask is matrix multiplication
then it will be broken into smaller matrix multiplications as
described above, and if it is QR decomposition, it will be broken into
smaller QR decompositions and matrix multiplications.

There are various obstacles to this simple approach.
First, the small QR decompositions generally have structure,
e.g., a pair of triangles. To exploit this structure fully
would complicate the recursive decomposition. (Or we could 
choose to ignore this structure, perhaps only on the smaller
subproblems, where the overhead would dominate.)

Second, it suggests that the data structure with which the
matrix is stored should be hierarchical as well, with
matrices stored as subblocks of subblocks \cite{elmroth2004recursive}.
This is certainly possible, but it differs significantly
from the usual data structures to which users are accustomed.
It also suggests that recent approaches based on decomposing
dense linear algebra operations into DAGs of subtasks \cite{buttari2007class,boboulin2008issues,kurzak2008qr,quintana-orti2008scheduling,quintana-orti2008design}
may need to be hierarchical, rather than have a single layer
of tasks. A single layer is a good match for the single socket 
multicore architectures that motivate these systems, but may
not scale well to, e.g., petascale architectures.

Third, it is not clear whether this approach best
accommodates machines that mix hierarchies of parallelism
and memory. For example, a multicore / multisocket / multirack 
computer will have also have disk, DRAM and various caches,
and it remains to be seen whether straightforward recursion
will minimize bandwidth and latency everywhere that
communication takes place within such an architecture.

Fourth and finally, all our analysis has assumed homogeneous
machines, with the same flop rate, bandwidth and latency
in all components. This assumption can be violated in
many ways, from having different bandwidth and latency
between racks, sockets, and cores on a single chip,
to having some specialized floating point units like GPUs.

It is most likely that an adaptive, ``autotuning'' approach
will be needed to deal with some of these issues, just
as it has been used for the simpler case of a matrix
multiplication.  Addressing all these issues is future work.


\appendix

\section*{Appendix}

\section{Structured local Householder QR flop counts}\label{S:localQR-flops}
\label{S:flopcounts} 

Here, we summarize floating-point operation counts for local
structured Householder QR factorizations of various matrices of
interest.  We count operations for both the factorization, and for
applying the resulting implicitly represented $Q$ or $Q^T$ factor to a
dense matrix.  Unless otherwise mentioned, we omit counts for BLAS 3
variants of structured Householder QR factorizations, as these
variants require more floating-point operations.  Presumably, the use
of a BLAS 3 variant indicates that small constant factors and
lower-order terms in the arithmetic operation count matter less to
performance than the BLAS 3 optimization.

\subsection{General formula}\label{SS:localQR-flops:general}

\subsubsection{Factorization}

Algorithm \ref{Alg:QR:qnxn} in Section
\ref{SS:TSQR:localQR:structured} shows a column-by-column Householder
QR factorization of the $qn \times n$ matrix of upper triangular $n
\times n$ blocks, using structured Householder reflectors.  We can
generalize this to an $m \times n$ matrix $A$ with a different nonzero
pattern, as long as the trailing matrix updates do not create nonzeros
below the diagonal in the trailing matrix.  This is true for all the
matrix structures encountered in the local QR factorizations in this
report.  A number of authors discuss how to predict fill in general
sparse QR factorizations; see, for example,
\cite{gilbert1992predicting}.  We do not need this level of
generality, since the structures we exploit do not cause fill.

The factorization proceeds column-by-column, starting from the left.
For each column, two operations are performed: computing the
Householder reflector for that column, and updating the trailing
matrix.  The cost of computing the Householder vector of a column
$A(j:m,j)$ is dominated by finding the norm of $A(j:m,j)$ and scaling
it.  If this part of the column contains $k_j$ nonzeros, this
comprises about $4k_j$ flops, not counting comparisons.  We assume
here that the factorization never creates nonzeros in the trailing
matrix; a necessary (but not sufficient) condition on $k_j$ is that it
is nondecreasing in $j$.

The trailing matrix update involves applying a length $m - j + 1$
Householder reflector, whose vector contains $k_j$ nonzeros, to the $m
- j + 1 \times c_j$ trailing matrix $C_j$.  The operation has the
following form:
\[
(I - \tau v_j v_j^T) C_j = C_j - v_j (\tau_j (v^T C_j)),
\]
in which $v_j$ is the vector associated with the Householder
reflector.  The first step $v_j^T C_j$ costs $2 c_j k_j$ flops, as we
do not need to compute with the zero elements of $v_j$.  The result is
a $1 \times c_j$ row vector and in general dense, so scaling it by
$\tau_j$ costs $c_j$ flops.  The outer product with $v_j$ then costs
$c_j k_j$ flops, and finally updating the matrix $C_j$ costs $c_j k_j$
flops (one for each nonzero in the outer product).  The total is $4
c_j k_j + c_j$.

When factoring an $m \times n$ matrix, $c_j = n - j$.  The total
number of arithmetic operations for the factorization is therefore
\begin{equation}\label{eq:localQR:gen:factor}
  \text{Flops}_{\text{Local, Factor}}(m,n) = \sum_{j=1}^n 4(n-j)k_j + 4k_j + (n-j)\,\text{flops.}
\end{equation}
We assume in this formula that $m \geq n$; otherwise, one would have
to sum from $j=1$ to $\min\{m,n\}$ in Equation
\eqref{eq:localQR:gen:factor}.

\subsubsection{Applying implicit $Q$ or $Q^T$ factor}

Applying an $m \times m$ $Q$ or $Q^T$ arising from QR on an $m \times
n$ matrix to an $m \times c$ matrix $C$ is like performing $n$
trailing matrix updates, except that the trailing matrix size $c$
stays constant.  This gives us an arithmetic operation count of
\begin{equation}\label{eq:localQR:gen:apply}
  \text{Flops}_{\text{Local, Apply}}(m,n,c) = \sum_{j=1}^n (4c k_j +
  4k_j + c)\, \text{flops.}
\end{equation}
The highest-order term in the above is
\[
4c \sum_{j=1}^n k_j.
\]
Note that the summation $\sum_{j=1}^n k_j$ is simply the number of
nonzeros in the collection of $n$ Householder vectors.  We assume in
this formula that $m \geq n$; otherwise, one would have to sum from
$j=1$ to $\min\{m,n\}$ in Equation \eqref{eq:localQR:gen:apply}.

\subsection{Special cases of interest}

\subsubsection{One block -- sequential TSQR}

\paragraph{Factorization}

The first step of sequential TSQR involves factoring a single $m/P
\times n$ input block.  This is the special case of a full matrix,
and thus the flop count is
\begin{equation}
\label{eq:localQR:seq:1block:factor}
\text{Flops}_{\text{Seq, 1 block, factor}}(m,n,P) = 
\frac{2mn^2}{P} - \frac{2n^3}{3} + O\left( \frac{mn}{P} \right),
\end{equation}
where we use $k_j = m/P - j + 1$ in Equation \eqref{eq:localQR:gen:factor}.
We assume in this formula that $m/P \geq n$.  This operation requires
keeping the following in fast memory:
\begin{itemize}
\item One $m/P \times n$ block of the input matrix $A$
\item Scratch space (at most about $n^2$ words)
\end{itemize}
The two-block factorization below consumes more fast memory, so it
governs the fast memory requirements in the sequential TSQR
factorization (see Section \ref{SSS:localQR-flops:seq:2blocks} for the
fast memory requirements).

\paragraph{Applying $Q$ or $Q^T$}

The cost in flops of applying an $m/P \times m/P$ $Q$ factor that
comes from the QR factorization of an $m/P \times n$ matrix to an $m/P
\times c$ matrix $C$ is
\begin{equation}
\label{eq:localQR:seq:1block:apply}
\text{Flops}_{\text{Seq, 1 block, apply}}(m,n,P) = 
\frac{4cmn}{P} - 2cn^2 + O\left( \frac{cm}{P} \right),
\end{equation}
where we use $k_j = m/P - j + 1$ in Equation
\eqref{eq:localQR:gen:apply}.  We assume in this formula that $m/P
\geq n$.

\subsubsection{Two blocks -- sequential TSQR}\label{SSS:localQR-flops:seq:2blocks}

\paragraph{Factorization}

For a $2m/P \times n$ local factorization with the top $m/P \times n$
block upper triangular and the lower $m/P \times n$ block full, we
have $k_j = 1 + m/P$ nonzeros in the $j^{\text{th}}$ Householder
reflector.  Thus, the flop count of the local QR factorization is
\begin{equation}
\label{eq:localQR:seq:2blocks:factor}
\text{Flops}_{\text{Seq, 2 blocks, factor}}(m,n,P) = 
\frac{2mn^2}{P} + \frac{2mn}{P} + O(n^2),
\end{equation}
using Equation \eqref{eq:localQR:gen:factor}.  We assume in this
formula that $m/P \geq n$.  For the case $m/P = n$ (two square $n
\times n$ blocks), this specializes to $k_j = 1+n$ and thus the flop
count is
\[
2 n^3 + O(n^2).
\]  
Without exploiting structure, the flop count would have been 
\[
\frac{10}{3} n^3 + O(n^2).
\]
Thus, the structured approach requires only about $3/5$ times as many
flops as standard Householder QR on the same $2n \times n$ matrix.

This operation requires keeping the following in fast memory:
\begin{itemize}
\item The previous block's $R$ factor ($n(n+1)/2$ words)
\item One $m/P \times n$ block of the input matrix $A$
\item Scratch space (at most about $n^2$ words)
\end{itemize}
Neglecting scratch space, the total fast memory requirement is 
\[
\frac{mn}{P} + \frac{n(n+1)}{2}\,\text{words.}
\]
Assume that fast memory can hold $W$ floating-point words.  We assume
that $m/P \geq n$, so clearly we must have
\[
W \geq \frac{n(n+1)}{2} + n^2 = \frac{3}{2}n^2 + \frac{n}{2}
\]
in order to solve the problem at all, no matter what value of $P$ we
use.  If this condition is satisfied, we can then pick $P$ so as to
maximize the block size (and therefore minimize the number of
transfers between slow and fast memory) in our algorithm.  Neglecting
the lower-order term of scratch space, the block size is maximized
when $P$ is minimized, namely when
\[
P = \frac{mn}{W - \frac{n(n+1)}{2}}.
\]

\paragraph{Applying $Q$ or $Q^T$}

Given the $2m/P \times 2m/P$ $Q$ factor arising from the QR
factorization in the previous paragraph, the cost in flops of applying
it to a $2m/P \times c$ matrix is given by Equation
\eqref{eq:localQR:gen:apply} as 
\begin{equation}
\label{eq:localQR:seq:2blocks:apply}
\text{Flops}_{\text{Seq, 2 blocks, apply}}(m,n,P) = 
\frac{4cmn}{P} + O\left( \frac{mn}{P} \right).
\end{equation}
We assume in this formula that $m/P \geq n$.

This operation requires keeping the following in fast memory:
\begin{itemize}
\item Two input blocks ($cm/P$ words each)
\item The local $Q$ factor ($n + mn/P$ words)
\item Scratch space of size $c \times n$ for $Q^T C$
\end{itemize}
The total fast memory requirements are therefore 
\[
\frac{(2c+n)m}{P} + (c+1)n\,\text{words},
\]
plus a lower-order term for scratch and stack space.  Assume that fast
memory can hold $W$ floating-point words.  We assume that $m/P \geq
n$, so clearly we must have
\[
W \geq (2c+n)n + (c+1)n
\]
in order to solve the problem at all, no matter what value of $P$ we
use.  If this condition is satisfied, we can then pick $P$ so as to
maximize the block size (and therefore minimize the number of
transfers between slow and fast memory) in our algorithm.  Neglecting
the lower-order term of scratch space, the block size is maximized
when $P$ is minimized, namely when
\[
P = \frac{(2c+n) m}{W - (c+1)n}.
\]
This formula makes sense because we assume that $W \geq 2(c+n)n +
(c+1)n$, so the denominator is always positive.

\subsubsection{Two or more blocks -- parallel TSQR}

\paragraph{Factorization}

For two $m/P \times n$ upper triangular blocks grouped to form a $2m/P
\times n$ matrix, we have $k_j = 1 + j$ and therefore the flop count
from Equation \eqref{eq:localQR:gen:factor}
is 
\begin{equation}
\label{eq:localQR:par:2blocks:factor}
\text{Flops}_{\text{Par, 2 blocks, factor}}(n,P) = 
\frac{2}{3}n^3 + O(n^2).
\end{equation}
We assume in this formula that $m/P \geq n$.  Without exploiting
structure, the flop count would have been
\[
\frac{4mn^2}{P} - \frac{2}{3} n^3 + O(n^2).
\]
Therefore, exploiting structure makes the flop count independent of
$m$.  We can generalize this case to some number $q \geq 2$ of the
$m/P \times n$ upper triangular blocks, which is useful for performing
TSQR with tree structures other than binary.  Here, $q$ is the
branching factor of a node in the tree.  In that case, we have $k_j =
1 + (q-1)j$ nonzeros in the $j^{th}$ Householder reflector, and
therefore the flop count from Equation \eqref{eq:localQR:gen:factor}
is
\begin{equation}
\label{eq:localQR:par:qblocks:factor}
\text{Flops}_{\text{Par, $q$ blocks, factor}}(n,q,P) = 
\frac{2}{3}(q-1) n^3 + O(qn^2).
\end{equation}
Again, we assume in this formula that $m/P \geq n$.  In the case $m/P
= n$, the optimization saves up to $2/3$ of the arithmetic operations
required by the standard approach.

\paragraph{Applying $Q$ or $Q^T$}

Given the $2m/P \times 2m/P$ $Q$ factor from the QR factorization in
the previous paragraph, the cost in flops of applying it to a $2m/P
\times c$ matrix $C$ is, from Equation \eqref{eq:localQR:gen:apply},
\begin{equation}
\label{eq:localQR:par:2blocks:apply}
\text{Flops}_{\text{Par, 2 blocks, apply}}(n,c) = 
2(c+1)n^2.
\end{equation}
We assume in this formula that $m/P \geq n$.  For the more general
case mentioned above of a $qm/P \times n$ $Q$ factor (with $q \geq
2$), the cost in flops of applying it to a $qm/P \times c$ matrix is
\begin{equation}
\label{eq:localQR:par:qblocks:apply}
\text{Flops}_{\text{Par, $q$ blocks, apply}}(n,q,c) = 
2(q-1) c n^2 + O(qn^2).
\end{equation}
Again, we assume in this formula that $m/P \geq n$.


\section{Sequential TSQR performance model}\label{S:TSQR-seq-detailed}

\subsection{Conventions and notation}

The sequential TSQR factorization operates on an $m \times n$ matrix,
divided into $P$ row blocks.  We assume that $m/P \geq n$, and we
assume without loss of generality that $P$ evenly divides $m$ (if not,
the block(s) may be padded with zeros).  We do not model a general 1-D
block cyclic layout, as it is only meaningful in the parallel case.
We assume that fast memory has a capacity of $W$ floating-point words
for direct use by the algorithms in this section.  We neglect the
lower-order amount of additional work space needed.  We additionally
assume read and write bandwidth are the same, and equal to $1/\beta$.
For simplicity of analysis, we assume no overlapping of computation
and communication; overlap could potentially provide another twofold
speedup.

The individual block operations (QR factorizations and updates) may be
performed using any stable QR algorithm.  In particular, the
optimizations in Section \ref{S:TSQR:localQR} apply.  When counting
floating-point operations and determining fast memory requirements, we
use the structured QR factorizations and updates described in Section
\ref{SS:TSQR:localQR:structured} and analyzed in Section
\ref{S:localQR-flops}.  In practice, one would generally also use the
BLAS 3 optimizations in Section \ref{SS:TSQR:localQR:BLAS3structured};
we omit them here because if their blocking factor is large enough,
they increase the flop counts by a small constant factor.  They also
increase the fast memory requirements by a small constant factor.  The
interaction between the BLAS 3 blocking factor and the block
dimensions in TSQR is complex, and perhaps best resolved by
benchmarking and search rather than by a performance model.

\subsection{Factorization}\label{SS:TSQR-seq-detailed:factor}

We now derive the performance model for sequential TSQR.  The
floating-point operation counts and fast memory requirements in this
section were derived in Appendix \ref{S:localQR-flops}.  Sequential
TSQR first performs one local QR factorization of the topmost $m/P
\times n$ block alone, at the cost of 
\begin{itemize}
  \item $\frac{2mn^2}{P} - \frac{2n^3}{3} + O(mn/P)$ flops (see
    Appendix \ref{S:localQR-flops}, Equation 
    \eqref{eq:localQR:seq:1block:factor}),
  \item one read from secondary memory of size $mn/P$, and
  \item one write to secondary memory, containing both the implicitly
    represented $Q$ factor (of size $mn/P - n(n+1)/2$), and the $\tau$
    array (of size $n$).
\end{itemize}
Then it does $P-1$ local QR factorizations of two $m/P \times n$
blocks grouped into a $2m/P \times n$ block.  In each of these local
QR factorizations, the upper $m/P \times n$ block is upper triangular,
and the lower block is a full matrix.  Each of these operations
requires
\begin{itemize}
  \item $\frac{2mn^2}{P} + O(mn/P)$ flops (see Appendix
    \ref{S:localQR-flops}, Equation
    \eqref{eq:localQR:seq:2blocks:factor}), 
  \item one read from slow memory of size $mn/P$, and
  \item one write to slow memory of size $mn/P + n$ (the
    Householder reflectors and the $\tau$ array).
\end{itemize}
The resulting modeled runtime is
\begin{multline}
\label{eq:TSQR:seq:modeltime:factor}
T_{\text{Seq.\ TSQR}}(m,n,P) = \\
  \alpha \left( 2P \right) + 
  \beta \left( 2mn + nP - \frac{n(n-1)}{2} \right) + 
  \gamma \left( 2mn^2 - \frac{2n^3}{3} \right).
\end{multline}

The above performance model leaves $P$ as a parameter.  Here, we pick
$P$ so as to maximize fast memory usage.  This minimizes the number of
memory transfers between slow and fast memory, and thus minimizes the
latency term in the model.  Suppose that fast memory can only hold $W$
words of data for sequential TSQR.  According to Section
\ref{SSS:localQR-flops:seq:2blocks}, the best choice of $P$ is the
minimum, namely
\[
P_{\text{opt}} = \frac{mn}{W - \frac{n(n+1)}{2}},
\]
which minimizes the number of transfers between slow and fast memory.  This
gives us a modeled runtime of
\begin{multline}
\label{eq:TSQR:seq:modeltimeW:factor}
T_{\text{Seq.\ TSQR}}(m,n,W) = 
  \alpha \left( \frac{2mn}{W - \frac{n(n+1)}{2}} \right) + \\
  \beta \left( 2mn - \frac{n(n+1)}{2} + \frac{mn^2}{W -
      \frac{n(n+1)}{2}} \right) + \\
  \gamma \left( 2mn^2 - \frac{2n^3}{3} \right).
\end{multline}

\subsection{Applying $Q$ or $Q^T$}\label{SS:TSQR-seq-detailed:apply}

The $Q$ and $Q^T$ cases are distinguished only by the order of
operations.  Suppose we are applying $Q$ or $Q^T$ to the dense $m
\times c$ matrix $C$.  The top block row of $C$ receives both a
one-block update (see Equation \eqref{eq:localQR:seq:1block:apply} in
Section \ref{S:localQR-flops}) and a two-block update (see Equation
\eqref{eq:localQR:seq:2blocks:apply} in Section \ref{S:localQR-flops}),
whereas the remaining block rows of $C$ each receive only a two-block
update.

The total number of arithmetic operations is
\[
4cmn - 2 c n^2 + O\left( \frac{cm}{P} \right) + O(mn) \,\text{flops.}
\]
Each block of the matrix $C$ is read from slow memory once and written
back to slow memory once.  Furthermore, each block of the $Q$ factor
is read from slow memory once.  Thus, the total number of transfers
between slow and fast memory is $3P$, and the total number of words
transferred between slow and fast memory is
\begin{multline*}
2cm + P \left( \frac{mn}{P} + n \right) - \frac{n(n+1)}{2} = 
(2c + n)m + n \cdot P - \frac{n(n+1)}{2}.
\end{multline*}
Altogether, we get the model
\begin{multline}
\label{eq:TSQR:seq:modeltimeP:apply}
T_{\text{Seq.\ TSQR apply}}(m,n,c,P) =
  \alpha \left( 
      3P
  \right) + \\
  \beta \left( 
     (2c + n)m + n \cdot P - \frac{n(n+1)}{2}
  \right) + 
  \gamma \left( 
      4cmn - 2 c n^2 
  \right).
\end{multline}

The above performance model leaves $P$ as a parameter.  Here, we
minimize $P$ so as to maximize fast memory usage.  This minimizes the
number of memory transfers between slow and fast memory, and thus
minimizes the latency term in the model.  Suppose that fast memory can
only hold $W$ words of data for sequential TSQR.  The two-block update
steps dominate the fast memory requirements.  Section
\ref{SSS:localQR-flops:seq:2blocks} describes the fast memory
requirements of the two-block update in detail.  In summary, choosing
\[
P_{\text{opt}} = \frac{(2c+n) m}{W - (c+1)n}.
\]
minimizes the number of transfers between slow and fast memory.  This
gives us a modeled runtime of
\begin{multline}
\label{eq:TSQR:seq:modeltimeW:apply}
T_{\text{Seq.\ TSQR apply}}(m,n,c,W) = 
  \alpha \left( 
      \frac{3 (2c+n) m}{W - (c+1)n}
  \right) + \\
  \beta \left( 
      (2c+n)m + \frac{2cmn + mn^2}{W - (c+1)n} - \frac{n(n+1)}{2}    
  \right) + \\
  \gamma \left( 
      4cmn - 2 c n^2 
  \right).
\end{multline}



\section{Sequential CAQR performance model}\label{S:CAQR-seq-detailed}

\subsection{Conventions and notation}

Sequential CAQR operates on an $m \times n$ matrix, stored in a $P_r
\times P_c$ 2-D block layout.  We do not model a fully general block
cyclic layout, as it is only helpful in the parallel case for load
balancing.  Let $P = P_r \cdot P_c$ be the number of blocks (not the
number of processors, as in the parallel case -- here we only use one
processor).  We assume without loss of generality that $P_r$ evenly
divides $m$ and $P_c$ evenly divides $n$.  The dimensions of a single
block of the matrix are $M \times N$, in which $M = m/P_r$ and $N =
n/P_c$.  Analogously to our assumption in Appendix
\ref{S:TSQR-seq-detailed}, we assume that $m \geq n$ and that $M \geq
N$.  Our convention is to use capital letters for quantities related
to blocks and the block layout, and lowercase letters for quantities
related to the whole matrix independent of a particular layout.

We assume that fast memory has a capacity of $W$ floating-point words
for direct use by the algorithms in this section,  neglecting 
lower-order amounts of additional work space.

The individual block operations (QR factorizations and updates) may be
performed using any stable QR algorithm.  In particular, the
optimizations in Section \ref{S:TSQR:localQR} apply.  When counting
floating-point operations and determining fast memory requirements, we
use the structured QR factorizations and updates described in Section
\ref{SS:TSQR:localQR:structured} and analyzed in Appendices
\ref{S:localQR-flops} and \ref{S:TSQR-seq-detailed}.  In practice, one
would generally also use the BLAS 3 optimizations in Section
\ref{SS:TSQR:localQR:BLAS3structured}; we omit them here because if
their blocking factor is large enough, they increase the flop counts
by a small constant factor.  They also increase the fast memory
requirements by a small constant factor.  The interaction between the
BLAS 3 blocking factor and the block dimensions in CAQR is complex,
and perhaps best resolved by benchmarking and search rather than by a
performance model.

\subsection{Factorization outline}\label{SS:CAQR-seq-detailed:factor}

Algorithms \ref{Alg:CAQR:seq:LL:outline} and
\ref{Alg:CAQR:seq:RL:outline} outline left-looking resp.\
right-looking variants of the sequential CAQR factorization.  Since it
turns out that the left-looking and right-looking algorithms perform
essentially the same number of floating-point operations, and send
essentially the same number of words in the same number of messages,
we will only analyze the right-looking algorithm.  

\begin{algorithm}[h]
\caption{Outline of left-looking sequential CAQR factorization}
\label{Alg:CAQR:seq:LL:outline}
\begin{algorithmic}[1]
\State{Assume: $m \geq n$ and $\frac{m}{P_r} \geq \frac{n}{P_c}$}
\For{$J = 1$ to $P_c$}
  \For{$K = 1$ to $J-1$}
  \State{Update panel $J$ (in rows $1$ to $m$
         and columns $(J-1)\frac{n}{P_c}+1$ to $J\frac{n}{P_c}$)
         using panel $K$}
  \EndFor
  \State{Factor panel $J$ (in rows $(J-1)\frac{n}{P_c} +1$ to $m$
         and columns $(J-1)\frac{n}{P_c}+1$ to $J\frac{n}{P_c}$)}
\EndFor
\end{algorithmic}
\end{algorithm}

\begin{algorithm}[h]
\caption{Outline of right-looking sequential CAQR factorization}
\label{Alg:CAQR:seq:RL:outline}
\begin{algorithmic}[1]
\State{Assume: $m \geq n$ and $\frac{m}{P_r} \geq \frac{n}{P_c}$}
\For{$J = 1$ to $P_c$}
  \State{Factor panel $J$ (in rows $(J-1)\frac{n}{P_c} +1$ to $m$
         and columns $(J-1)\frac{n}{P_c}+1$ to $J\frac{n}{P_c}$)}
  \State{Update trailing panels to right, 
         (in rows $(J-1)\frac{n}{P_c} +1$ to $m$
         and columns $J\frac{n}{P_c}+1$ to $n$)
         using the current panel}
\EndFor
\end{algorithmic}
\end{algorithm}

Indeed, we need only to replace the loop in 
Algorithm~\ref{Alg:CAQR:seq:RL:outline} by a summation,
and the calls to ``factor panel'' and ``update panel''
with uses of the formulas for $T_{\text{Seq.\ TSQR}}()$
from equation~(\ref{eq:TSQR:seq:modeltime:factor}) and
for $T_{\text{Seq.\ TSQR apply}}()$ 
from equation~(\ref{eq:TSQR:seq:modeltimeP:apply}):

\begin{eqnarray} 
\label{eq:CAQR:seq:modeltime:P}
T_{\text{Seq.\ CAQR}} (m,n,P_c,P_r) & \leq &
          \sum_{J=1}^{P_c} 
              T_{\text{Seq.\ TSQR}}(m-(J-1)\frac{n}{P_c},\frac{n}{P_c},P_r) 
          \nonumber \\
&   &       + (P_c-J)T_{\text{Seq.\ TSQR apply}}(m-(J-1)\frac{n}{P_c},
                     \frac{n}{P_c},\frac{n}{P_c},P_r)  \nonumber \\
& = & \alpha \left[ \frac{3}{2} P(P_c-1) \right] +  \nonumber \\
&   & \beta \left[ \frac{3}{2} mn \left( P_c + \frac{4}{3} \right) 
                   - \frac{1}{2} n^2 P_c + O\left(n^2 + nP\right) \right] + \\
&   & \gamma \left[ 2n^2m - \frac{2}{3}n^3 \right]  \nonumber
\end{eqnarray}
where we have ignored lower order terms, 
and used $P_r$ as an upper
bound on the number of blocks in each column 
(in the last argument of each function), since this only increases the run
time slightly, and is simpler to evaluate than for the true number of
blocks $P_r - \lfloor (J-1) \frac{n P_r}{m P_c} \rfloor$.

\subsection{Choosing $P$, $P_r$ and $P_c$ to optimize runtime}\label{SS:CAQR-seq-detailed:opt}

From the above formula for
$T_{\text{Seq.\ CAQR}} (m,n,P_c,P_r)$, we see that the runtime is
an increasing function of $P_r$ and $P_c$, so that we would like
to choose them as small as possible, within the limits imposed by
the fast memory size $W$, to minimize the runtime.

The CAQR step requiring the most fast memory is the two-block update
of a panel.
Each trailing matrix block has $m/P_r$ rows and $n/P_c$ columns, so
given the formula in Section~\ref{SSS:localQR-flops:seq:2blocks}, the
fast memory requirement is
\[
\frac{3mn}{P} + \frac{n^2}{P_c^2} + \frac{n}{P_c}
\leq \frac{4mn}{P} +\sqrt{\frac{mn}{P}}
\,\text{words,}
\]
plus a lower-order term for scratch and stack space.  
For simplicity, we approximate this by $\frac{4mn}{P}$.
To minimize runtime, we want to minimize $P$ subject
to $\frac{4mn}{P} \leq W$, i.e. we choose $P = \frac{4mn}{W}$.
But we still need to choose $P_r$ and $P_c$ subject to 
$P_r P_c = P$.
 
Examining $T_{\text{Seq.\ CAQR}} (m,n,P_c,P_r)$, again, we see
that if $P$ is fixed, the runtime is also an increasing function
of $P_c$, which we therefore want to minimize. But we are assuming
$\frac{m}{P_r} \geq \frac{n}{P_c}$, or $P_c \geq \frac{nP_r}{m}$.
The optimal choice is therefore  $P_c = \frac{nP_r}{m}$
or $P_c = \sqrt{\frac{nP}{m}}$, which also means
$\frac{m}{P_r} = \frac{n}{P_c}$, i.e., the blocks in the algorithm
are square. This choice of $P_r = \frac{2m}{\sqrt{W}}$ and 
$P_c = \frac{2n}{\sqrt{W}}$ therefore minimizes the runtime,
yielding
\begin{eqnarray} 
\label{eq:CAQR:seq:modeltime:P:opt}
T_{\text{Seq.\ CAQR}} (m,n,W) & \leq &
      \alpha \left[ 12 \frac{mn^2}{W^{3/2}} \right] +
      \beta \left[ 3 \frac{mn^2}{\sqrt{W}} + 
      O\left( \frac{mn^2}{W} \right) \right] + \nonumber \\
&   &      \gamma \left[ 2mn^2 - \frac{2}{3}n^3 \right].
\end{eqnarray}

We note that the bandwidth term is proportional to $\frac{mn^2}{\sqrt{W}}$,
and the latency term is $W$ times smaller,
both of which match (to within constant factors), the lower bounds
on bandwidth and latency described in Corollary~\ref{corollary:SeqCAQR}
in Section~\ref{sec:SeqCAQR}.

\section{Parallel TSQR performance model}\label{S:TSQR-par-detailed}
\label{S:TSQR:counts}

\subsection{Conventions and notation}

The parallel TSQR factorization operates on an $m \times n$ matrix in
a 1-D block layout on $P$ processors.  We assume that $m/P \geq n$,
and we assume without loss of generality that $P$ evenly divides $m$
(if not, the block(s) may be padded with zeros).  Furthermore, we
assume that the number of processors $P$ is a power of two: $P =
2^{L-1}$ with $L = \log_2 P$.  For simplicity, we do not model a
general 1-D block cyclic layout here, and we assume no overlapping of
computation and communication (overlap could potentially provide
another twofold speedup).

The individual block operations (QR factorizations and updates) may be
performed using any stable QR algorithm.  In particular, the
optimizations in Section \ref{S:TSQR:localQR} apply.  When counting
floating-point operations and determining fast memory requirements, we
use the structured QR factorizations and updates described in Section
\ref{SS:TSQR:localQR:structured}.  In practice, one would generally
also use the BLAS 3 optimizations in Section
\ref{SS:TSQR:localQR:BLAS3structured}; we omit them here because if
their blocking factor is large enough, they increase the flop counts
by a small constant factor.  They also increase the fast memory
requirements by a small constant factor.  The interaction between the
BLAS 3 blocking factor and the block dimensions in TSQR is complex,
and perhaps best resolved by benchmarking and search rather than by a
performance model.

\subsection{Factorization}

We now derive the performance model for parallel TSQR on a binary tree
of $P$ processors.  We restrict our performance analysis to the block
row, reduction based Algorithm \ref{Alg:TSQR:allred:blkrow}.  The
all-reduction-based version has the same number of flops on the critical
path (the root process of the reduction tree), but it requires $2q$
parallel messages per level of the tree on a $q$-ary tree, instead of
just $q-1$ parallel messages to the parent node at that level in the
case of a reduction.  When counting the number of floating-point
operations for each step of the factorization, we use the counts
derived in Appendix \ref{S:localQR-flops}.

A parallel TSQR factorization on a binary reduction tree performs the
following computations along the critical path: 
\begin{itemize}
\item One local QR factorization of a fully dense $m/P \times n$
  matrix ($2mn^2/P - \frac{n^3}{3} + O(mn/P)$ flops)
\item $\log_2 P$ factorizations, each of a $2n \times n$ matrix
  consisting of two $n \times n$ upper triangular matrices
  ($\frac{2}{3}n^3 + O(n^2)$ flops)
\end{itemize}
Thus, the total flop count is
\[
\text{Flops}_{\text{Par.\ TSQR}}(m,n,P) =
\frac{2mn^2}{P} 
- \frac{n^3}{3} 
+ \frac{2}{3} n^3 \log_2 P 
+ O\left( \frac{mn}{P} \right).
\]
The factorization requires $\log_2 P$ messages, each of size
$n(n+1)/2$ (the $R$ factors at each step of the tree).

\subsection{Applying $Q$ or $Q^T$}\label{SS:TSQR-par-detailed:apply}

Suppose we are applying $Q$ or $Q^T$ to the dense $m \times c$ matrix
$C$.  (The $Q$ and $Q^T$ cases are distinguished only by the order of
operations.)  We assume the matrix $C$ is distributed in the same 1-D
block layout on $P$ processors as was the original matrix $A$.  The
total number of arithmetic operations is
\[
\text{Flops}_{\text{Par.\ TSQR apply}}(m,n,P) =
\frac{4cmn}{P} 
+ 2cn^2 (\log_2(P) - 1)
+ O\left( \frac{(c+n)m}{P} \right).
\]
Suppose that a reduction (rather than an all-reduction) is used to
apply the $Q$ factor.  Then, at each inner node of the reduction tree,
one processor receives an $n \times c$ block of the matrix $C$ from
its neighbor, updates it, and sends the result back to its neighbor.
So there are two messages per inner node, each of size $cn$.  This
gives a total of $2 \log_2 P$ messages, and a total communication
volume of $2 cn \log_2 P$ words.  If an all-reduction is used, there
is only one message per inner node along the critical path, and that
message is of size $cn$.  This gives a total of $\log_2 P$ messages,
and a total communication volume of $cn \log_2$ words.


\section{Parallel CAQR performance model}\label{S:CAQR-par-detailed}

\subsection{Conventions and notation}

In this section, we model the performance of the parallel CAQR
algorithm described in Section \ref{S:CAQR}.  Parallel CAQR operates
on an $m \times n$ matrix $A$, stored in a 2-D block cyclic layout on
a $P_r \times P_c$ grid of $P$ processors.  We assume without loss of
generality that $P_r$ evenly divides $m$ and that $P_c$ evenly divides
$n$ (if not, the block(s) may be padded with zeros).  We assume no
overlapping of computation and communication (overlap could
potentially provide another twofold speedup).

The individual block operations (QR factorizations and updates) may be
performed using any stable QR algorithm.  In particular, the
optimizations in Section \ref{S:TSQR:localQR} apply.  When counting
floating-point operations and determining fast memory requirements, we
use the structured QR factorizations and updates described in Section
\ref{SS:TSQR:localQR:structured}.  In practice, one would generally
also use the BLAS 3 optimizations in Section
\ref{SS:TSQR:localQR:BLAS3structured}; we omit them here because if
their blocking factor is large enough, they increase the flop counts
by a small constant factor.  They also increase the fast memory
requirements by a small constant factor.  The interaction between the
BLAS 3 blocking factor and the block dimensions in CAQR is complex,
and perhaps best resolved by benchmarking and search rather than by a
performance model.


\subsection{Factorization}\label{SS:CAQR-par-detailed:factor}

First, we count the number of floating point arithmetic operations
that CAQR performs along the critical path.  We compute first the cost
of computing the QR factorization using Householder transformations of
a $m \times n$ matrix A (using \lstinline!DGEQR2!).  The cost of
computing the $j$th Householder vector is given by the cost of
computing its Euclidian norm and then by scaling the vector.  This
involves $3(m-j+1)$ flops and $(m-j+1)$ divides.  The cost of updating
the trailing $A(j:m, j+1:n)$ matrix by $I-\tau v_j v_j^T$ is
$4(n-j)(m-j+1)$.  The total number of flops is:
\begin{multline*}
  3 \sum_{j=1}^{n} (m-j+1) + 4 \sum_{j=1}^{n-1} (n-j)(m-j+1) = \\
  2 m n^2 
  - \frac{2 n^3}{3} 
  + mn 
  + \frac{n^2}{2} 
  + \frac{n}{3} = \\
2 mn^2 - \frac{2 n^3}{3} + O(mn).
\end{multline*}
The total number of divides is around $mn - n^2 / 2$.

The Householder update of a matrix $(I - Y T^T Y^T) C$, where $Y$ is a
$m \times n$ unit lower trapezoidal matrix of Householder vectors and
$C$ is a $m \times q$ matrix, can be expressed as:
\begin{eqnarray*}
C = 
\begin{pmatrix}
C_0 \\
C_1 \\
\end{pmatrix}
= 
\left(
I -
\begin{pmatrix}
Y_0 \\
Y_1 \\
\end{pmatrix}
\cdot 
\begin{array}{cc} 	
T^T \\
\\
\end{array}
\cdot
\begin{pmatrix}
Y_0 \\
Y_1 \\
\end{pmatrix}^T
\right)
\begin{pmatrix}
C_0 \\
C_1 \\
\end{pmatrix}
\end{eqnarray*}
in which $Y_0$ is a $n \times n$ unit lower triangular matrix and
$Y_1$ is a rectangular matrix. The total number of flops is around $qn
(4m - n - 1) \approx qn (4m-n)$.  We described in Section
\ref{SS:TSQR:localQR:trailing} how to perform the trailing matrix
update.  The breakdown of the number of flops in each step is:
\begin{itemize}
\item $W = Y_0^T C_0$ $\rightarrow$ $n (n-1) q$ flops.
\item $ W = W + Y_1^T C_1$ $\rightarrow$ $2n (m-n) q$ flops.
\item $ W = T^T W$ $\rightarrow$ $n^2 q $ flops.
\item $C_0 = C_0 - Y_0 W$ $\rightarrow$ $n^2 q$ flops.
\item $C_1 = C_1 - Y_1 W$ $\rightarrow$ $2n (m-n) q$ flops.
\end{itemize}

We consider now the computation of the upper triangular matrix $T$ used in the
$(I - Y T Y^T)$ representation of the Householder vectors (\lstinline!DLARFT!
routine).  This consists of $n$ transformations of the form $(I - \tau v_i
v_i^T)$.  Consider $Y$, a $m \times n$ unit lower trapezoidal matrix of
Householder vectors.  The matrix $T$ is an upper triangular matrix of dimension
$n \times n$ that is obtained in $n$ steps.  At step $j$, the first $j-1$
columns of $T$ are formed.  The $j$-th column is obtained as follows:
\begin{multline*}
T(1:j, 1:j) = \\
\begin{pmatrix}
      T(1:j-1, 1:j-1) & - \tau T(1:j-1, 1:j-1) Y^T (:, 1:j-1) v_j) \\
                      & \tau \\
\end{pmatrix}
\end{multline*}
in which $v_j$ is the $j$th Householder vector of length $m-j+1$.  This
is obtained by computing first $w = - \tau Y^T(:, 1:j-1) v_j$ (matrix
vector multiply of $2 ( j-1) (m - j +1)$ flops ) followed by the
computation $T(1:j-1, j) = T(1:j-1, 1:j-1) w$ (triangular matrix
vector multiplication of $(j-1)^2$ flops).  The total cost of forming
$T$ is:
\begin{eqnarray*}
  m n^2 - \frac{n^3}{3} - mn + \frac{n^2}{2} - \frac{n}{6} 
\approx 
  m n^2 - \frac{n^3}{3}
\end{eqnarray*}

The new QR factorization algorithm also performs Householder updates
of the form
\begin{eqnarray*}
C = 
\left(
\begin{array}{cc} 	
C_0 \\
C_1 \\
\end{array}
\right)
= 
\left(
I -
\left(
\begin{array}{cc} 	
I \\
Y_1 \\
\end{array}
\right)
\cdot 
\begin{array}{cc} 	
T^T \\
\\
\end{array}
\cdot
\left(
\begin{array}{cc} 	
I \\
Y_1 \\
\end{array}
\right)^T
\right)
 \left(
\begin{array}{cc} 	
C_0 \\
C_1 \\
\end{array}
 \right)
\end{eqnarray*}
in which $Y_1$ is a $n \times n$ upper triangular matrix and $C$ is a $2n
\times q$ matrix.  The total number of flops is $3 n^2 q + 6 n q$.
The following outlines the number of floating-point operations corresponding 
to each operation performed during this update:
\begin{itemize}
\item $W = Y_1^T C_1$ $\rightarrow$ $n(n+1) q$ flops.
\item $ W = W + C_0$ $\rightarrow$ $n q$ flops.
\item $ W = T^T W$ $\rightarrow$ $n (n+1) q $ flops.
\item $C_0 = C_0 - W$ $\rightarrow$ $n q$ flops.
\item $C_1 = C_1 - Y_1 W$ $\rightarrow$ $n (n+2) q$ flops.
\end{itemize}

Forming the upper triangular matrix $T$ used in the above Householder
update corresponds now to computing $- \tau T(1:j-1, 1:j-1)
Y_1^T(1:j-1,1:j-1) v_j(n+1:n+j)$.  $v_j$ is the $j$th Householder
vector composed of $1$ in position $j$ and nonzeros in positions $n+1,
\ldots n+j+1$.  First $w = - \tau Y_1^T(1:j-1,1:j-1) v_j(n+1:2n)$ is
computed (triangular matrix vector multiply of $j ( j-1) $ flops),
followed by $T(1:j-1, j) = T(1:j-1, 1:j-1) w$ (triangular matrix
vector multiplication of $(j-1)^2$ flops).  The total number of flops
is
\begin{equation}
\sum_{j=1}^{n} j (j-1) + \sum_{j=1}^{n} (j-1)^2 \approx 2 \frac{n^3}{3}
\end{equation}

We are now ready to estimate the time of CAQR.
\begin{enumerate}
\item The column of processes that holds panel $j$ computes a TSQR
  factorization of this panel.  The Householder vectors are stored in a
  tree as described in Section \ref{S:TSQR:impl}.

\begin{multline}
\label{step1}
\gamma \left[ 
    \frac{2b^2 m_j}{P_r} 
    + \frac{2 b^3}{3} \left( \log P_r - 1 \right) 
\right] + \\
\gamma_d \left[ 
    \frac{m_j b}{P_r} 
    + \frac{b^2}{2} \left( \log P_r - 1 \right) 
\right] + \\
\alpha \log P_r 
+ \beta \frac{b^2}{2} \log P_r 
\end{multline}

\item Each processor $p$ that belongs to the column of processes
  holding panel $j$ broadcasts along its row of processors the $m_j /
  P_r \times b$ rectangular matrix that holds the two sets of
  Householder vectors.  Processor $p$ also broadcasts two arrays of
  size $b$ each, containing the Householder factors $\tau_p$.

\begin{equation}
\label{step2}
     \alpha \left( 2 \log P_c \right) + \beta \left( \frac{m_j b}{P_r} + 2 b \right) \log P_c
\end{equation}

\item Each processor in the same row template as processor $p$, $0
  \leq i < P_r$, forms $T_{p0}$ (first two terms in the number of
  flops) and updates its local trailing matrix $C$ using $T_{p0}$ and
  $Y_{p0}$ (last term in the number of flops). (This computation
  involves all processors and there is no communication.)

\begin{equation}
\label{step3}
  \left[ b^2 \frac{m_j}{P_r} - \frac{b^3}{3} 
+ 
b \frac{n_j - b}{P_c} \left( 4 \frac{m_j}{P_r} -
  b \right) \right] \gamma
\end{equation}

\item {\bf for} $k = 1$ to $\log P_r$ {\bf do}

Processors that lie in the same row as processor $p$, where $0 \leq p
< P_r$ equals $first\_proc(p,k)$ or $target\_first\_proc(p,k)$ perform:

\begin{enumerate}
\item Processors in the same template row as $target\_first\_proc(p,k)$ form
locally $T_{level(p,k),k}$.  They also compute local pieces of $W =$ \linebreak 
$Y_{level(p,k),k}^T C_{target\_first\_proc(p,k)}$, leaving the results distributed.
This computation is overlapped with the communication
in (\ref{ap_step_comm1}).

\begin{equation}
\label{step4}
   \left[ \frac{2b^3}{3}  +
   b (b+1) \frac{n_j - b}{P_c} \right] \gamma 
\end{equation}

\item 
\label{ap_step_comm1}
Each processor in the same row of the grid as $first\_proc(p,k)$ sends to the
  processor in the same column and belonging to the row of
  $target\_first\_proc(p,k)$ the local pieces of $C_{first\_proc(p,k)}$.

\begin{equation}
\label{step5}
     \alpha +\frac{b (n_j - b)}{P_c} \beta
\end{equation}

\item Processors in the same template row as $target\_first\_proc(p,k)$ compute
  local pieces of \\ 
  $W = T_{level(p,k),k}^T \left( C_{first\_proc(p,k)}+ W  \right)$.

\begin{equation}
\label{step6}
  \left( b (b + 2) \frac{n_j - b}{P_c} \right) \gamma
\end{equation}

\item 
\label{ap_step_comm2}
Each processor in the same template row as $target\_first\_proc(p,k)$
sends to the processor in the same column and belonging to the row
template of $first\_proc(p,k)$ the local pieces of $W$.

\begin{equation}
\label{step7}
     \alpha + \beta \left( \frac{b (n_j - b)}{P_c} \right)
\end{equation}

\item Processors in the same template row as $first\_proc(p,k)$ complete
locally the rank-$b$ update $C_{first\_proc(p,k)} = C_{first\_proc(p,k)} - W$
(number of flops in Equation~\ref{step8}).  Processors in the same
template row as $target\_first\_proc(p,k)$ locally complete the
rank-$b$ update \linebreak
$C_{target\_first\_proc(p,k)} = C_{target\_first\_proc(p,k)} - Y_{level(p,k),k} W$ (number of
flops in Equation~\ref{step9}).  The latter computation is overlapped
with the communication in~(\ref{ap_step_comm2}).

\begin{eqnarray}
\label{step8}
  \gamma \left( b \frac{n_j - b}{P_c} \right)  \\
\label{step9}
  \gamma \left( b (b+2) \frac{n_j - b}{P_c} \right)
\end{eqnarray}

\end{enumerate}

 {\bf end for}

\end{enumerate}

We can express the total computation time over a rectangular grid of
processors $T_{\text{Par.\ CAQR}}(m, n, P_r, P_c)$ as a sum over the number of iterations
of the previously presented steps.  The number of messages is $n/b (3
\log P_r + 2 \log P_c)$.  The volume of communication is:
\begin{eqnarray*}
 \sum_{j=1}^{n/b} \left(  
\frac{b^2}{2} \log P_r + \frac{m_j b}{P_r} \log P_c + 2b \log P_c + \frac{2 b (n_j -
  b)}{P_c} \log P_r
\right) = \\
\left( \frac{nb}{2} + \frac{n^2}{P_c} \right) \log P_r + 
\left( 2n + \frac{mn - n^2/2 }{P_r} \right) \log P_c
\end{eqnarray*}

The total number of flops corresponding to each step is given by the
following, in which ``(Eq.\ $S$)'' (for some number $S$) is a
reference to Equation ($S$) in this section.
\[
\begin{aligned}
\text{(Eq.\ \ref{step1})} & \sum_{j=1}^{n/b} &\approx
    &\frac{2nmb - n^2 b + nb^2 }{P_r} + \frac{2 b^2 n}{3} ( \log{P_r} - 1) \\
\text{(Eq.\ \ref{step3})} & \sum_{j=1}^{n/b} &\approx 
    &\frac{1}{P} \left( 
         2mn^2  
         - \frac{2}{3} n^3 
     \right) 
     + \frac{1}{P_r} \left( 
         mnb 
         + \frac{nb^2}{2} 
         - \frac{n^2b}{2} 
     \right) 
     + \frac{n^2b}{2P_c} 
     - \frac{b^2n}{3} \\
\text{(Eq.\ \ref{step4})} & \sum_{j=1}^{n/b} &\approx 
    &\left( 
         \frac{2 b^2 n}{3} 
         + \frac{n^2 (b+1)}{2 P_c} 
     \right) \log P_r \\
\text{(Eq.\ \ref{step6})} & \sum_{j=1}^{n/b} &\approx 
    &\frac{n^2 (b+2)}{2 P_c} \log P_r \\
\text{(Eq.\ \ref{step8})} & \sum_{j=1}^{n/b} &\approx 
    &\frac{n^2}{2 P_c} \log P_r \\
\text{(Eq.\ \ref{step9})} & \sum_{j=1}^{n/b} &\approx 
    &\frac{n^2 (b+2)}{2 P_c} \log P_r \\
\end{aligned}
\]


The total computation time of parallel CAQR can be estimated as:
\begin{multline}
\label{Eq:ap_time_latavoidQR}
T_{\text{Par.\ CAQR}} (m, n, P_r, P_c) = \\
\gamma \left[ 
    \frac{2n^2(3m-n)}{3P}
    + \frac{b n^2}{2 P_c} 
    + \frac{3bn(2m-n)}{2P_r}
    + \right. \\ \left. 
    \left( 
        \frac{4 b^2 n}{3} 
        + \frac{n^2 (3b+5)}{2 P_c}
    \right) \log P_r 
    - b^2 n
\right] + \\
\gamma_d \left[ 
    \frac{mn - n^2/2}{P_r}
    + \frac{bn}{2} \left( \log(P_r) - 1 \right)
\right] + \\
\alpha \left[
    \frac{3n}{b} \log P_r
    + \frac{2n}{b} \log P_c
\right] + \\
\beta \left[
    \left( 
        + \frac{n^2}{P_c} 
        \frac{bn}{2} 
    \right) \log P_r
    + \left( 
        \frac{ mn - \frac{n^2}{2} }{P_r} 
        + 2n 
    \right) \log P_c
\right].
\end{multline}


\section{Sca\-LA\-PACK's out-of-DRAM QR fac\-tor\-i\-za\-tion \texttt{PFDGEQRF}}
\label{S:seqLL} 
\label{S:PFDGEQRF}

LAPACK Working Note \#118 describes an out-of-DRAM QR factorization
routine \lstinline!PFDGEQRF!, which is implemented as an extension of
ScaLAPACK \cite{dazevedo1997design}.  It uses ScaLAPACK's existing
parallel in-DRAM panel factorization (\lstinline!PDGEQRF!) and update
(\lstinline!PDORMQR!) routines.  Thus, it is able to exploit
parallelism within each of these steps, assuming that the connection
to the filesystem is shared among the processors.  It can also take
advantage of the features of parallel filesystems for block reads and
writes.

We use the algorithm and communication pattern underlying
\lstinline!PFDGEQRF! as a model for a reasonable sequential
out-of-DRAM implementation of Householder QR.  This means we assume
that all operations in fast memory run sequentially, and also that the
connection to slow memory is sequential.  These assumptions are fair,
because we can always model a parallel machine as a faster sequential
machine, and model multiple connections to slow memory as a single
higher-bandwidth connection.  From now on, when we say
\lstinline!PFDGEQRF! without further qualifications, we mean our
sequential out-of-DRAM model.  We also describe the algorithm for
applying the $Q$ factor computed and implicitly stored by
\lstinline!PFDGEQRF!, either as $Q$ or as $Q^T$, to an $m \times r$
matrix $B$ which need not fit entirely in fast memory.

We will show that the estimated runtime of \lstinline!PFDGEQRF!, as a
function of the fast memory size $W$, is
\begin{multline*}
T_{\text{\texttt{PFDGEQRF}}}(m, n, W) = 
\alpha \left[
    \frac{2mn}{W} + \frac{mn^2}{2W} - \frac{n}{2}
\right] + \\
\beta \left[
    \frac{3mn}{2} - \frac{3n^2}{4} 
    + \frac{m n}{W} \left(
        \frac{m n}{2} 
        - \frac{n^2}{6}.
    \right)
\right] + \\
\gamma \left[
    2mn^2 - \frac{2n^3}{3} + O(mn) 
\right].
\end{multline*}

\subsection{Conventions and notation}

Algorithm \ref{Alg:seqLL:outline} computes the Householder QR
factorization of an $m \times n$ matrix.  We assume that $m \geq n$
and that the matrix is sufficiently large to warrant not storing it
all in fast memory.  Algorithm \ref{Alg:seqLL:outline:apply} applies
the implicitly stored $Q$ factor from this factorization, either as
$Q$ or as $Q^T$, to an $m \times r$ matrix $B$.  We assume in both
cases that fast memory has a capacity of $W$ floating-point words.
When computing how much fast memory an algorithm uses, we neglect
lower-order terms, which may include scratch space.  We additionally
assume read and write bandwidth are the same, and equal to $1/\beta$.
For simplicity of analysis, we assume no overlapping of computation
and communication; overlap could potentially provide another twofold
speedup.

The individual panel operations (QR factorizations and updates) may be
performed using any stable QR algorithm.  In particular, the
optimizations in Section \ref{S:TSQR:localQR} apply.  When counting
floating-point operations and determining fast memory requirements, we
use the flop counts for standard (unstructured) QR factorizations and
updates analyzed in Appendix \ref{S:localQR-flops}.  In practice, one
would generally also use the BLAS 3 optimizations described in Section
\ref{SS:TSQR:localQR:BLAS3structured}; we omit them here because if
their blocking factor is large enough, they increase the flop counts
by a small constant factor.  They also increase the fast memory
requirements by a small constant factor.  The interaction between the
BLAS 3 blocking factor and the panel widths $b$ and $c$ is complex,
and perhaps best resolved by benchmarking and search rather than by a
performance model.

\subsection{Factorization}\label{SS:seqLL:factor}\label{SS:seqLL:fact}

\begin{algorithm}[h]
\caption{Outline of ScaLAPACK's out-of-DRAM QR factorization (\texttt{PFDGEQRF})}
\label{Alg:PFDGEQRF:outline}\label{Alg:seqLL:outline}
\begin{algorithmic}[1]
\For{$j = 1$ to $n-c$ step $c$}
  \State{Read current panel (columns $j : j + c - 1$) from slow memory}
  \For{$k = 1$ to $j - 1$ step $b$}
    \State{Read left panel (columns $k : k + b - 1$) from slow memory}
    \State{Apply left panel to current panel (in fast memory)}
  \EndFor
  \State{Factor current panel (in fast memory)}
  \State{Write current panel to slow memory}
\EndFor
\end{algorithmic}
\end{algorithm}

\lstinline!PFDGEQRF! is a left-looking QR factorization method.  The
code keeps two panels in fast memory: a left panel of fixed width $b$,
and the current panel being factored, whose width $c$ can expand to
fill the available memory.  Algorithm \ref{Alg:PFDGEQRF:outline} gives
an outline of the code, without cleanup for cases in which $c$ does
not evenly divide $n$ or $b$ does not evenly divide the current column
index minus one.  Algorithm \ref{Alg:PFDGEQRF:detail} near the end of
this section illustrates this ``border cleanup'' in detail, though we
do not need this level of detail in order to model the performance of
\lstinline!PFDGEQRF!.

\subsubsection{Communication pattern}\label{SSS:seqLL:common}

Algorithm \ref{Alg:seqLL:outline} shares a common communication
pattern with many variants of the QR factorization.  All these
factorization variants keep two panels in fast memory: a left panel
and a current panel.  For each current panel, the methods sweep from
left to right over the collection of left panels, updating the current
panel with each left panel in turn.  They then factor the current
panel and continue.  Applying the $Q$ or $Q^T$ factor from the
factorization (as in Algorithm \ref{Alg:seqLL:outline:apply}) has a
similar communication pattern, except that the trailing matrix is
replaced with the $B$ matrix.  If we model this communication pattern
once, we can then get models for all such methods, just by filling in
floating-point operation counts for each.  Any QR factorization which
works ``column-by-column,'' such as (classical or modified)
Gram-Schmidt, may be used in place of Householder QR without changing
the communication model.  This is because both the left panel and the
current panel are in fast memory, so neither the current panel update
nor the current panel factorization contribute to communication.
(Hence, this unified model only applies in the sequential case, unless
each processor contains an entire panel.)

\subsubsection{Fast memory usage}

The factorization uses $(b+c)m$ words of fast memory at once, not
counting lower-order terms or any BLAS 3 optimizations that the panel
factorizations and updates may use.  In order to maximize the amount
of fast memory used, we choose $b$ and $c$ so that $(b+c)m = W$.  If
$2m > W$ then we cannot use this factorization algorithm without
modification, as in its current form, at least two columns of the
matrix must be able to fit in fast memory.  The parameter $b$ is
typically a small constant chosen to increase the BLAS 3 efficiency,
whereas one generally chooses $c$ to fill the available fast memory,
for reasons which will be shown below.  Thus, in many cases we will
simplify the analysis by taking $b = 1$ and $c \approx W/m$.

An important quantity to consider is $mn/W$.  This is the theoretical
lower bound on the number of reads from slow memory (it is a latency
term).  It also bounds from below the number of slow memory writes,
assuming that we use the usual representation of the $Q$ factor as a
collection of dense Householder reflectors, and the usual
representation of the $R$ factor as a dense upper triangular matrix.

\subsubsection{Number of words transferred}

Algorithm \ref{Alg:seqLL:outline} transfers about
\begin{multline}
\label{eq:pfdgeqrf:bw:bc}
\sum_{j = 1}^{\frac{n}{c}} \left( 
  2c \left( 
    m - cj + 1
  \right) +
  \sum_{k=1}^{\frac{c(j-1)}{b}} b \left(
    m - bk + 1
  \right)
\right) = \\
\left( 
    \frac{3mn}{2}
    - \frac{3n^2}{4}
    + \frac{3n}{2}
\right)
+ \frac{bn}{4}
- \frac{13 cn}{12}
+ \frac{1}{c} \left(
    \frac{mn^2}{2} 
    - \frac{n^3}{6}
    + \frac{n^2}{2} 
    - \frac{bn^2}{4}
\right)
\end{multline}
floating-point words between slow and fast memory.  Our goal is to
minimize this expression as a function of $b$ and $c$.  Though it is a
complicated expression, one can intuitively guess that it is minimized
when $b = 1$ and $c$ is as large as possible, because if the current
panel width is larger, then one needs to iterate over all left panels
fewer times ($n/c$ times, specifically), whereas the number of words
read and written for all current panels is $\Theta(mn)$ regardless of
the values of $b$ and $c$.  If we take $b = 1$ and $c \approx W/m$,
the most significant terms of the above sum are
\[
\frac{3mn}{2} 
- \frac{3n^2}{4} 
+ \frac{m n}{W} \left(
    \frac{m n}{2} 
    - \frac{n^2}{6}.
\right)
\]
In the case of sequential TSQR, the lower bound on the number of words
transferred between slow and fast memory is $2mn$ (see Section
\ref{SS:lowerbounds:1d}); the third and fourth terms above, taken
together, can be arbitrarily larger when $mn \ll W$.  In the case of
sequential CAQR (see Section \ref{sec:SeqCAQR}), the lower bound is
\[
\frac{ \frac{mn^2}{4} - \frac{n^2}{8} \left( \frac{n}{2} + 1
  \right)}{(8W)^{1/2}} - W.
\]
If we assume $m \geq n$, then \lstinline!PFDGEQRF! transfers
$\Omega(m / \sqrt{W})$ times more words between fast and slow memory
than the lower bound.

\subsubsection{Number of slow memory accesses}

The total number of slow memory reads and writes performed by
Algorithm \ref{Alg:seqLL:outline} is about
\begin{equation}
\label{eq:pfdgeqrf:lat:bc}
\text{Messages}_{\text{\texttt{PFDGEQRF}}}(m, n, b, c) =
\sum_{j = 1}^{\frac{n}{c}} \left( 
  2 +
  \sum_{k=1}^{\frac{c(j-1)}{b}} 1
\right) =
\frac{n}{2bc} +
2\frac{n}{c} -
\frac{n}{2b}.
\end{equation}
This quantity is always minimized by taking $c$ as large as possible,
and therefore $b$ as small as possible.  If we let $b = 1$ and
approximate $c$ by $W/m$, the number of slow memory reads and writes
comes out to
\begin{equation}
\label{eq:pfdgeqrf:lat}
\text{Messages}_{\text{\texttt{PFDGEQRF}}}(m, n, W) =
\frac{2mn}{W} + \frac{mn^2}{2W} - \frac{n}{2}
\end{equation}
accesses.  

In the case of sequential TSQR (see Section \ref{SS:lowerbounds:1d}),
the lower bound on the number of messages between slow and fast memory
is $mn/W$ (or $2mn/W$ if we need to write the result to slow memory).
Thus, the above communication pattern is a factor of $n/4$ away from
optimality with respect to the latency term.  In the case of
sequential CAQR (see Section \ref{sec:SeqCAQR}, the lower bound is
\[
\frac{\frac{mn^2}{4}-\frac{n^2}{8}(\frac{n}{2}+1)}{(8W^3)^{1/2}} - 1
\]
messages.  Thus, \lstinline!PFDGEQRF! uses $\Theta(\sqrt{W})$ more
messages than the lower bound.

\subsubsection{Floating-point operations}

\lstinline!PFDGEQRF!, when run sequentially, performs exactly the same
floating-point operations as standard Householder QR, but in a
different order.  Thus, the flop count is the same, namely
\[
2mn^2 - \frac{2n^3}{3} + O(mn).
\]
This is independent of the parameters $b$ and $c$.  If BLAS 3
optimizations are used in the panel factorizations and updates, as is
likely in a practical implementation, then the flop count may go up by
a small constant factor.

\subsubsection{Runtime}

The estimated runtime of \lstinline!PFDGEQRF! as a function of the
parameters $b$ and $c$ is
\begin{multline}\label{eq:PFDGEQRF:runtime:bc}
T_{\text{\texttt{PFDGEQRF}}}(m, n, b, c) = 
\alpha \left[
    \frac{n}{2bc} +
    2\frac{n}{c} -
    \frac{n}{2b}
\right] + \\
\beta \left[
    \left( 
        \frac{3mn}{2}
        - \frac{3n^2}{4}
        + \frac{3n}{2}
    \right)
    + \frac{bn}{4}
    - \frac{13 cn}{12}
    + \frac{1}{c} \left(
        \frac{mn^2}{2} 
        - \frac{n^3}{6}
        + \frac{n^2}{2} 
        - \frac{bn^2}{4}
    \right)
\right] + \\
\gamma \left[
    2mn^2 - \frac{2n^3}{3} + O(mn) 
\right].
\end{multline}
When we choose $b$ and $c$ optimally according to the discussion
above, the estimated runtime is
\begin{multline}\label{eq:PFDGEQRF:runtime:W}
T_{\text{\texttt{PFDGEQRF}}}(m, n, W) = 
\alpha \left[
    \frac{2mn}{W} + \frac{mn^2}{2W} - \frac{n}{2}
\right] + \\
\beta \left[
    \frac{3mn}{2} - \frac{3n^2}{4} 
    + \frac{m n}{W} \left(
        \frac{m n}{2} 
        - \frac{n^2}{6}.
    \right)
\right] + \\
\gamma \left[
    2mn^2 - \frac{2n^3}{3} + O(mn) 
\right].
\end{multline}

\subsection{Applying $Q^T$}\label{SS:seqLL:apply}

\begin{algorithm}[h]
  \caption{Applying $Q^T$ from Algorithm \ref{Alg:PFDGEQRF:outline} to
    an $n \times r$ matrix $B$}
\label{Alg:seqLL:outline:apply}
\begin{algorithmic}[1]
\For{$j = 1$ to $r - c$ step $c$}
  \State{Read current panel of $B$ (columns $j : j + c - 1$)}
  \For{$k = 1$ to $j - 1$ step $b$}
    \State{Read left panel of $Q$ (columns $k : k + b - 1$)}
    \State{Apply left panel to current panel}
  \EndFor
  \State{Write current panel of $B$}
\EndFor
\end{algorithmic}[1]
\end{algorithm}

Algorithm \ref{Alg:seqLL:outline:apply} outlines the communication
pattern for applying the full but implicitly represented $Q^T$ factor
from Algorithm \ref{Alg:PFDGEQRF:outline} to an $n \times r$ dense
matrix $B$.  We assume that the matrix $B$ might be large enough not
to fit entirely in fast memory.  Applying $Q$ instead of $Q^T$ merely
changes the order of iteration over the Householder reflectors, and
does not change the performance model, so we can restrict ourselves
without loss of generality to computing $Q^T B$.

\subsubsection{Fast memory usage}

Algorithm \ref{Alg:seqLL:outline:apply} uses at most $bm + cn$ words
of fast memory at once, not counting lower-order terms or any BLAS 3
optimizations that the panel applications may use.  In order to
maximize the amount of fast memory used, we take $bm + cn = W$.  Note
that this is different from the factorization, in case $m \neq n$.  If
$m+n > W$ then we cannot use this algorithm at all, even if $b = 1$
and $c = 1$.  


\subsubsection{Number of words transferred}

Clearly, each panel of $B$ is read and written exactly once, so the
matrix $B$ contributes $2mr$ words to the total.  If we count the
number of times that each word of the $Q$ factor is read, we obtain
about 
\[
\sum_{j = 1}^{\frac{r}{c}} \sum_{k = 1}^{\frac{n}{b}} \left( bm -
  \frac{b(b-1)}{2} \right) = 
\frac{nr}{2c}
- \frac{bnr}{2c}
+ \frac{mnr}{c}
\]
The formula is clearly minimized when $c$ is maximized, namely when $b
= 1$ and $c \approx W/m$.  In that case, we obtain a total number of
words transferred of about
\[
mr 
+ \frac{m^2 n r}{W}\, \text{words.}
\]

\subsubsection{Number of slow memory accesses}

The total number of data transfers between slow and fast memory is about
\[
\sum_{j=1}^{r/c} \left( 2 + \sum_{k=1}^{n/b} 1 \right) = 
\frac{2r}{c} + \frac{nr}{bc} = 
\frac{r}{c} \left( 2 + \frac{n}{b} \right). 
\]
If we choose $b = 1$ and $c \approx W/m$, we get
\[
\frac{2mr}{W} + \frac{2mnr}{W},
\]
whereas if we choose $c = 1$ and $b \approx W/m$, then we get
\[
2r + \frac{mnr}{W}.
\]
It's clear that maximizing $b$ minimizes the latency term.

\subsubsection{Floating-point operations}

Applying the $Q$ factor costs 
\[
4mnr - 2n^2 r + O(mr)\, \text{flops.}
\]
This is independent of the parameters $b$ and $c$.  If BLAS 3
optimizations are used in the panel operations, as is likely in a
practical implementation, then the flop count may go up by a small
constant factor.

\begin{algorithm}[h]
\caption{More detailed outline of ScaLAPACK out-of-DRAM Householder QR
  factorization (PFDGEQRF), with border cleanup}
\label{Alg:PFDGEQRF:detail}
\begin{algorithmic}[1]
\For{$j = 1$ to $\left( \lfloor \frac{n}{c} \rfloor - 1 \right) c + 1$
  step $c$}
  \State{Read current panel (columns $j : j+c-1$, rows $1 : m$)}
  \For{$k = 1$ to $\left( \lfloor \frac{j-1}{b} \rfloor - 1 \right) b
    + 1$, step $b$} 
    \State{Read left panel (columns $k : k+b-1$, lower trapezoid,
       starting at row $k$)}
    \State{Apply left panel to rows $k : m$ of current panel}
  \EndFor
  \State{$k := \lfloor \frac{j-1}{b} \rfloor b + 1$}        
  \State{Read left panel (columns $k : j-1$, lower trapezoid, starting
    at row $k$)}
  \State{Apply left panel to rows $k : m$ of current panel}
  \State{Factor current panel (rows $1 : m$)}
  \State{Write current panel (rows $1 : m$)}
\EndFor
\State{$j := \lfloor \frac{n}{c} \rfloor c + 1$}
\State{Read current panel (columns $j : n$, rows $1 : m$)}
\For{$k = 1$ to $\left( \lfloor \frac{j-1}{b} \rfloor - 1 \right) b +
  1$, step $b$}
  \State{Read left panel (columns $k : k+b-1$, lower trapezoid,
    starting at row $k$)}
  \State{Apply left panel to rows $k : m$ of current panel}
\EndFor
\State{$k := \lfloor \frac{j-1}{b} \rfloor b + 1$}
\State{Read left panel (columns $k : j-1$, lower trapezoid, starting
  at row $k$)}
\State{Apply left panel to current panel}
\State{Factor current panel (rows $1 : m$)}
\State{Write current panel (rows $1 : m$)}
\end{algorithmic}
\end{algorithm}




\end{comment}
\newcommand{\nunder}{{\underline{n}}}
\newcommand{\munder}{{\underline{m}}}

\section{Communication Lower Bounds from Calculus}
\label{S:CommLowerBoundsFromCalculus}
 
\subsection{Summary}
In this section we address communication lower bounds for TSQR
needed in subsection~\ref{sec:par_TSQR_comm_bounds},
asking how much data two (or more) processors have to communicate in order
to compute the QR decomposition of a matrix whose rows are
distributed across them.
We analyze this in a way that applies to more general situations:
Suppose processor 1 and processor 2 each own some of the arguments
of a function $f$ that processor 1 wants to compute. What is the least
volume of communication required to compute the function?
We are interested in smooth functions of real or complex arguments,
and so will use techniques from calculus rather than modeling
the arguments as bit strings.

In this way, we will derive necessary conditions on the function $f$
for it to be evaluable by communicating fewer than all of its arguments
to one processor. We will apply these conditions to various linear
algebra operations to capture our intuition that it is in fact necessary 
to move all the arguments to one processor for correct evaluation of $f$.

\subsection{Communication lower bounds for one-way communication between
2 processors}
Suppose $x^{(m)} \in \RR^m$ is owned by processor 1 (P1) and
$y^{(n)} \in \RR^n$ is owned by P2; we use superscripts to
remind the reader of the dimension of each vector-valued variable or function.
Suppose P1 wants to compute  
$f^{(r)}(x^{(m)},y^{(n)}): \RR^{m} \times \RR^{n} \rightarrow \RR^r$.
We first ask how much information P2 has to send to P1, assuming
it is allowed to send one message, consisting of $\nunder \leq n$ real
numbers, which themselves could be functions of $y^{(n)}$.
In other words, we ask if functions $h^{(\nunder )} (y^{(n)}): \RR^{n} \rightarrow \RR^{\nunder}$
and $F^{(r)} (x^{(m)}, z^{(\nunder)})  : \RR^{m} \times \RR^{\nunder} \rightarrow \RR^r$, 
exist such that 
$f^{(r)}(x^{(m)},y^{(n)}) = F^{(r)} (x^{(m)}, h^{(\nunder)} (y^{(n)}))$.
When $\nunder = n$, the obvious choice is to send the original data $y^{(n)}$,
so that $h^{(\nunder)} (y^{(n)}) = y^{(n)}$ is the identity function and
$f^{(r)} = F^{(r)}$. The interesting question is whether we can send
less information, i.e. $\nunder < n$.

Unless we make further restrictions on the function $h$ we are allowed
to use, it is easy to see that we can always choose $\nunder =1$, i.e. send
the least possible amount of information: We do this by using a 
space-filling curve \cite{sagan1994space} to represent each $y^{(n)} \in \RR^{(n)}$ by
one of several preimages $\tilde{y} \in \RR$. In other words, 
$h^{(1)} (y^{(n)})$ 
maps $y^{(n)}$ to a scalar $\tilde{y}$ that P1 can map back to 
$y^{(n)}$ by a space filling curve. 
This is obviously unreasonable, since it implies we could try to 
losslessly compress $n$ 64-bit floating point numbers into one 64-bit 
floating point number. 
However, by placing some reasonable smoothness restrictions on  
the functions we use, since we can only hope to evaluate (piecewise) smooth
functions in a practical way anyway, we will see that we can draw useful 
conclusions about practical computations. 
To state our results, we use the notation $J_x f(x,y)$ to denote the
$r \times m$ Jacobian matrix of $f^{(r)}$ with respect to the arguments
$x^{(m)}$.  Using the above notation, we state

\lemma{
Suppose it is possible to compute 
$f^{(r)}(x^{(m)},y^{(n)})$ on P1 by communicating $\nunder < n$ words
$h^{(\nunder )} (y^{(n)})$ from P2 to P1, and evaluating
$f^{(r)}(x^{(m)},y^{(n)}) = F^{(r)} (x^{(m)}, h^{(\nunder )} (y^{(n)}))$.
Suppose $h^{(\nunder )}$ and $F^{(r)}$ are continuously differentiable
on open sets. Then necessary conditions for this to be possible
are as follows.
\begin{enumerate}
\item Given any fixed $y^{(n)}$ in the open set, then for all
$x^{(m)}$ in the open set, 
the rows of $J_y f(x,y)$ must lie
in a fixed subspace of $\RR^n$
of dimension at most $\nunder < n$.
\item Given any fixed $\tilde{y}^{(\nunder)} \in \RR^{\nunder}$ satisfying
$\tilde{y}^{(\nunder)} = h^{(\nunder)} (y^{(n)})$ for some $y^{(n)}$ in
the interior of the open set, there is 
a set $C \subset \RR^{n}$ containing $y^{(n)}$,
of dimension at least $n-\nunder$,
such that for each 
$x$, $f(x,y)$ is constant for $y \in C$.
\item If $r=n$, and for each fixed $x$, $f^{(r)}(x,y^{(n)})$ is a bijection,
then it is necessary and sufficient to send $n$ words from P2 to P1 
to evaluate $f$.
\end{enumerate}
}

\begin{proof}
Part 1 is proved simply by differentiating, using the chain rule,
and noting the dimensions of the Jacobians being multiplied:
\[
J_y^{(r \times n)} f^{(r)}(x,y) = J_h^{(r \times \nunder)} F^{(r)}(x,h) 
\cdot J_y^{(\nunder \times n)} h^{(\nunder)} (y)
\]
implying that for all $x$, each row of $J_y^{(r \times m)} f^{(r)} (x,y)$
lies in the space spanned by the $\nunder$ rows of 
$J_y^{(\nunder \times n)} h^{(\nunder)} (y)$.

Part 2 is a consequence of the implicit function theorem.
Part 3 follows from part 2, since if the function is a bijection,
then there is no set $C$ along which $f$ is constant.
\end{proof} \rm

Either part of the lemma can be used to derive 
lower bounds on the volume of communication needed to compute $f(x,y)$,
for example
by choosing an $\nunder$ equal to the lower bound minus 1, and 
confirming that either necessary condition in the Lemma is
violated, at least in some open set. 

We illustrate this for a simple matrix factorization problem.

\corollary{
Suppose P1 owns the $r_1 \times c$ matrix $A_1$, and
P2 owns the $r_2 \times c$ matrix $A_2$, with
$r_2 \geq c$. Suppose P1 wants to compute the
$c \times c$ Cholesky factor $R$ of
$R^T \cdot R = A_1^T \cdot A_1 + A_2^T \cdot A_2$,
or equivalently the $R$ factor in the $QR$ decomposition
of $\bmat{c} A_1 \\ A_2 \emat$. Then P2 has to communicate
at least $c(c+1)/2$ words to P1, and it is possible to
communicate this few, namely either the entries
on and above the diagonal of the symmetric $c \times c$ matrix $A_2^T \cdot A_2$, 
or the entries of its
Cholesky factor $R$, so that $R^T \cdot R = A_2^T \cdot A_2$ 
(equivalently, the $R$ factor of the $QR$ factorization of $A_2$).
}

\begin{proof}
That it is sufficient to communicate the $c(c+1)/2$ entries described
above is evident. We use  Corollary~1 to prove that these many words are
necessary. We use the fact that mapping between 
the entries on and above the diagonal of the symmetric positive definite
matrix and its Cholesky factor is a bijection
(assuming positive diagonal entries of the Cholesky factor).
To see that for any fixed $A_1$, $f(A_1,R) = $ the Cholesky factor
of $A_1^T \cdot A_1 + R^T \cdot R$ is a bijection,
note that it is a composition of three bijections:
the mapping from $R$ to the entries on and above the 
diagonal of $Y = A_2^T \cdot A_2$, the entries on and
above the diagonal of $Y$ and those on and above the diagonal
of $X = A_1^T \cdot A_1 + Y$, and the mapping between the entries on
and above the diagonal of $X$ and its Cholesky factor $f(A_1,R)$.
\end{proof} \rm

\subsection{Reduction operations}

We can extend this result slightly to make it apply to the case of
more general reduction operations, where one processor P1 is trying to
compute a function of data initially stored on multiple other processors
P2 through P$s$. We suppose that there is a tree
of messages leading from these processors eventually reaching P1.
Suppose each P$i$ only sends data up the tree, so that the communication
pattern forms a DAG (directed acylic graph) with all paths ending at P1. 
Let P$i$'s data be denoted $y^{(n)}$.
Let all the variables on P1 be denoted $x^{(m)}$,
and treat all the other variables on the other processors as constants.
Then exactly the same analysis as above applies, and we can conclude that
{\em every} message along the unique path from P$i$ to P1 has the same 
lower bound on its size, as determined by Lemma~1.
This means Corollary~1 extends to include reduction operations where
each operation is a bijection between one input (the other being fixed)
and the output. In particular, it applies to TSQR.

We emphasize again that using a real number model to draw conclusions about
finite precision computations must be done with care. For example,
a bijective function depending on many variables could hypothetically 
round to the same floating point output for all floating point inputs,
eliminating the need for any communication or computation 
for its evaluation. But this is not the
case for the functions we are interested in.

Finally, we note that the counting must be done slightly 
differently for the QR decomposition of complex data,
because the diagonal entries $R_{i,i}$ are generally
taken to be real. Alternatively, there is a degree of
freedom in choosing each row of $R$, which can be
multiplied by an arbitrary complex number of absolute
value 1.

\subsection{Extensions to two-way communication}

While the result of the previous subsection is adequate for the results 
of this paper,
we note that it may be extended as follows. For motivation, suppose that
P1 owns the scalar $x$, and wants to evaluate the polynomial
$\sum_{i=1}^{n} y_i x^{i-1}$, where P2 owns the vector $y^{(n)}$.
The above results can be used to show that P2 needs to send $n$
words to P1 (all the coefficients of the polynomial, for example).
But there is an obvious way to communicate just 2 words:
(1) P1 sends $x$ to P2, (2) P2 evaluates the polynomial, and
(3) P2 sends the value of the polynomial back to P1.

More generally, one can imagine $k$ phases, during each of which
P1 sends one message to P2 and then P2 sends one message to P1.
The contents of each message can be any smooth functions of all 
the data available to the sending processor, either originally 
or from prior messages. At the end of the $k$-th phase, P1 then
computes $f(x,y)$.

More specifically, the computation and communication proceeds as
follows:
\begin{itemize}
\item In Phase 1, P1 sends $g_1^{(m_1)} (x^{(m)})$ to P2
\item In Phase 1, P2 sends $h_1^{(n_1)} (y^{(n)}, g_1^{(m_1)} (x^{(m)}))$ to P1
\item In Phase 2, P1 sends $g_2^{(m_2)} (x^{(m)}, h_1^{(n_1)} (y^{(n)}, g_1^{(m_1)} (x^{(m)})))$ to P2
\item In Phase 2, P2 sends 
$h_2^{(n_2)} (y^{(n)}, 
g_1^{(m_1)} (x^{(m)}),
g_2^{(m_2)} (x^{(m)}, h_1^{(n_1)} ( y^{(n)}, g_1^{(m_1)} (x^{(m)}))))$ to P1
\item $\dots$
\item In Phase $k$, P1 sends
$g_k^{(m_k)} (x^{(m)}, h_1^{(n_1)} (\dots), h_2^{(n_2)} (\dots) , \dots , 
h_{k-1}^{(n_{k-1})} (\dots) )$ to P2
\item In Phase $k$, P2 sends
$h_k^{(n_k)} (y^{(n)}, g_1^{(m_1)} (\dots), g_2^{(m_2)} (\dots) , \dots , 
g_{k}^{(m_k)} (\dots) )$ to P1
\item P1 computes 
\begin{eqnarray*}
f^{(r)} (x^{(m)} , y^{(n)} ) & = &
F^{(r)} ( x^{(m)}, h_1^{(n_1)} ( y^{(n)}, g_1^{(m_1)} (x^{(m)})), 
\\ & & 
h_2^{(n_2)} ( y^{(n)}, 
g_1^{(m_1)} (x^{(m)}), 
g_2^{(m_2)} (x^{(m)}, h_1^{(n_1)} ( y^{(n)}, g_1^{(m_1)} (x^{(m)})))),
\\ & & \dots
\\ & & h_k^{(n_k)} (y^{(n)}, g_1^{(m_1)} (\dots), g_2^{(m_2)} (\dots) , \dots , 
g_{k}^{(m_k)} (\dots) ))
\end{eqnarray*}
\end{itemize}

\lemma{
Suppose it is possible to compute 
$f^{(r)}(x^{(m)},y^{(n)})$ on P1 by the scheme described above.
Suppose all the functions involved are continuously differentiable
on open sets. Let $\nunder = \sum_{i=1}^k n_i$ and
$\munder = \sum_{i=1}^k m_i$.
Then necessary conditions for this to be possible
are as follows.
\begin{enumerate}
\item 
Suppose $\nunder < n$ and $\munder \leq m$, ie. P2 cannot communicate
all its information to P1, but P1 can potentially send its information
to P2.
Then there is a set $C_x \subset \RR^m$ 
of dimension at least $m-\munder$ 
and a set $C_y \subset \RR^n$ 
of dimension at least $n-\nunder$ such that
for $(x,y) \in C = C_x \times C_y$, 
the value of $f(x,y)$ is independent of $y$.
\item 
If $r=n=m$, and for each fixed $x$ or fixed $y$, 
$f^{(r)}(x^{(m)},y^{(n)})$ is a bijection,
then it is necessary and sufficient to send $n$ words from P2 to P1 
to evaluate $f$.
\end{enumerate}
}

\begin{proof}
We define the sets $C_x$ and $C_y$ by the following constraint equations,
one for each communication step in the algorithm:
\begin{itemize}
\item 
$\tilde{g}_1^{(m_1)} = g_1^{(m_1)} (x^{(m)})$ is a fixed constant, 
placing $m_1$ smooth constraints on $x^{(m)}$.
\item 
In addition to the previous constraint,
$\tilde{h}_1^{(n_1)} = h_1^{(n_1)} (y^{(n)}$, $g_1^{(m_1)} (x^{(m)}))$ 
is a fixed constant,
placing $n_1$ smooth constraints on $y^{(n)}$.
\item 
In addition to the previous constraints, \linebreak
$\tilde{g}_2^{(m_2)} = g_2^{(m_2)} (x^{(m)}, h_1^{(n_1)} (y^{(n)}, g_1^{(m_1)} (x^{(m)})))$ 
is a fixed constant, placing $m_2$ more smooth constraints on $x^{(m)}$.
\item 
In addition to the previous constraints, \linebreak
$\tilde{h}_2^{(n_2)} = h_2^{(n_2)} (y^{(n)}, 
g_1^{(m_1)} (x^{(m)}),
g_2^{(m_2)} (x^{(m)}, h_1^{(n_1)} ( y^{(n)}, g_1^{(m_1)} (x^{(m)}))))$ 
is a fixed constant,
placing $n_2$ more smooth constraints on $y^{(n)}$.
\item 
\dots
\item 
In addition to the previous constraints, \linebreak
$\tilde{g}_k^{(m_k)} = g_k^{(m_k)} (x^{(m)}, h_1^{(n_1)} (\dots), h_2^{(n_2)} (\dots) , \dots , 
h_{k-1}^{(n_{k-1})} (\dots) )$ is a fixed constant,
placing $m_k$ more smooth constraints on $x^{(m)}$.
\item 
In addition to the previous constraints, \linebreak
$\tilde{h}_k^{(n_k)} = h_k^{(n_k)} (y^{(n)}, g_1^{(m_1)} (\dots), g_2^{(m_2)} (\dots) , \dots , 
g_{k}^{(m_k)} (\dots) )$ is a fixed constant,
placing $n_k$ more smooth constraints on $y^{(n)}$.
\end{itemize}
Altogether, we have placed 
$\nunder = \sum_{i=1}^k n_i < n$ smooth constraints on $y^{(n)}$ and
$\munder = \sum_{i=1}^k m_i \leq m$ smooth constraints on $x^{(m)}$,
which by the implicit function theorem define surfaces 
$C_y (\tilde{h}_1^{(n_1)} , \dots , \tilde{h}_k^{(n_k)} )$
and 
$C_x (\tilde{g}_1^{(m_1)} , \dots , \tilde{g}_k^{(m_k)} )$,
of dimensions at least
$n - \nunder > 0$ and $m - \munder \geq 0$, respectively,
and parameterized by
$\{\tilde{h}_1^{(n_1)} , \dots , \tilde{h}_k^{(n_k)} \}$ and
$\{\tilde{g}_1^{(m_1)} , \dots , \tilde{g}_k^{(m_k)} \}$,
respectively.
For $x \in C_x$ and $y \in C_y$, the values communicated by 
P1 and P2 are therefore constant. Therefore, for $x \in C_x$
and $y \in C_y$, $f(x,y) = F(x,h_1, \dots, h_k)$ depends only on $x$,
not on $y$. This completes the first part of the proof.

For the second part, we know that if $f(x,y)$ is a bijection
in $y$ for each fixed $x$, then by the first part
we cannot have $\nunder < n$, because otherwise 
$f(x,y)$ does not depend on $y$ for certain values of $x$,
violating bijectivity. But if we can send $\nunder = n$ words from
P2 to P1, then it is clearly possible to compute $f(x,y)$ by
simply sending every component of $y^{(n)}$ from P2 to P1 explicitly.
\end{proof}

\corollary{
Suppose P1 owns the $c$-by-$c$ upper triangular matrix $R_1$,
and P2 owns the $c$-by-$c$ upper triangular matrix $R_2$, and
P1 wants to compute the R factor in the QR decomposition of
$\bmat{c} R_1 \\ R_2 \emat$. Then it is neccessary and sufficient
to communicate $c(c+1)/2$ words from P2 to P1 (in particular, the
entries of $R_2$ and sufficient).
}

We leave extensions to general communication patterns among
multiple processors to the reader.


\bibliographystyle{siam}
\bibliography{qr}
\end{document}